\documentclass[10pt]{amsart}
\usepackage[T1]{fontenc}
\usepackage{lmodern}
\usepackage{geometry}
\usepackage[latin1] {inputenc}
\usepackage{amsmath}
\usepackage{amsfonts, amssymb, textcomp}
\usepackage[colorlinks=flase, linkcolor=red,urlcolor=green, citecolor=blue]{hyperref}
\usepackage{subeqnarray}

\usepackage{latexsym}
\usepackage{fancyhdr}
\usepackage{longtable}
\usepackage{amsmath, amssymb}
\usepackage{graphicx}

\numberwithin{equation}{section}

\theoremstyle{plain}
\newtheorem{theorem}{Theorem}[section]
\newtheorem{lemma}[theorem]{Lemma}
\newtheorem{proposition}[theorem]{Proposition}
\newtheorem{corollary}[theorem]{Corollary}
\theoremstyle{definition}
\newtheorem{definition}[theorem]{Definition}

\newtheorem{remark}[theorem]{Remark}

\newcommand{\RR}{\mathbb{R}}
\newcommand{\CC}{\mathbb{C}}
\newcommand{\NN}{\mathbb{N}}

\newcommand{\supp}{\operatorname{supp}}

\newcommand{\ess}{\operatorname{ess}}
\let\on=\operatorname

\newenvironment{demo}[1]{\par\smallskip\noindent{\bf #1.}}{\par\smallskip}

\makeatletter
\@namedef{subjclassname@2020}{%
  \textup{2020} Mathematics Subject Classification}
\makeatother

\title[On inclusion relations of weighted $L^p$-type spaces]
{On inclusion relations of weighted $L^p$-type spaces defined in terms of weight function matrices}

\author[G.~Schindl]{Gerhard Schindl}

\address{G.~Schindl: Fakult\"at f\"ur Mathematik, Universit\"at Wien, Oskar-Morgenstern-Platz~1, A-1090 Wien, Austria.}
\email{gerhard.schindl@univie.ac.at}

\begin{document}

\begin{abstract}
We introduce new weighted $L^p$-type spaces defined in terms of weight function matrices and characterize the inclusion relations in terms of the defining matrices. Moreover, we provide a detailed study concerning the coincidence with the common (non-weighted) $L^p$-spaces, the (non-)triviality of such weighted spaces and investigate their translation invariance. The obtained results are then applied to particular weight function matrices which are expressed in terms of one single weight function and a positive real parameter. Also variations of this new weighted setting are discussed; more precisely weighted Banach (sub-)spaces of $L^p$ and when weighting the Fourier image of appropriate Banach spaces of functions. The general framework allows to describe the known ultradifferentiable weight function setting by Beurling-Bj\"{o}rck which is more original than the approach presented by Braun, Meise and Taylor. When applying the characterization of the inclusion relations to Beurling-Bj\"{o}rck-type spaces we are able to emphasize the difference between both ultradifferentiable weight function settings: We construct a technical (counter-)example which is a weight in the sense of Beurling-Bj\"{o}rck but violates the standard and crucial convexity condition needed in the Braun-Meise-Taylor setting.
\end{abstract}

\thanks{This research was funded in whole or in part by the Austrian Science Fund (FWF) 10.55776/PAT9445424}
\date{\today}

\maketitle

\section{Introduction}\label{Intro}
When dealing with ultradifferentiable function classes and analogously defined weighted spaces, in the literature ``classically'' one can find two settings: either using a \emph{weight sequence} $\mathbf{M}$ or a \emph{weight function} $\omega$. The ultradifferentiable weight sequence case is more original, we refer to \cite{Komatsu73} and \cite{mandelbrojtbook}. On the other hand, in the weight function setting there exist two methods: First the growth of the Fourier transform $\widehat{f}$ (of compactly supported functions $f$) has been measured w.r.t. weighted $L^1$ or $L^{\infty}$-based norms; see \cite{Beurling61}, \cite{Bjorck66}, \cite{PetzscheVogt}, the explanations in \cite{BraunMeiseTaylor90} and Section \ref{examplesect} for precise definitions and a summary. More precisely, the boundedness of $\int_{\RR^d}|\widehat{f}(x)|e^{\ell\omega(x)}dx$ or $\sup_{x\in\RR^d}|\widehat{f}(x)|e^{\ell\omega(x)}$ for some $\ell>0$ \emph{(``Roumieu-type'')} resp. for all $\ell>0$ \emph{(``Beurling-type'')} is required and so the admissible growth is controlled by $\omega$ directly \emph{(``Beurling-Bj\"{o}rck setting'').}

Later, in \cite{BraunMeiseTaylor90} the nowadays frequently used definition has been given: Let $\varphi_{\omega}: t\mapsto\omega(e^t)$ and let $\varphi^{*}_{\omega}$ be the \emph{Legendre-Fenchel-Young conjugate} given by
\begin{equation}\label{LFYconj}
\varphi^{*}_{\omega}(x):=\sup_{y\in\RR}\{xy-\varphi_{\omega}(y)\},\;x\in\RR.
\end{equation}
Then the derivatives of the functions under investigation are weighted analogously to the weight sequence framework and the weights are given by $W^{(\ell)}_j:=\exp(\frac{1}{\ell}\varphi^{*}_{\omega}(\ell j))$, $j\in\NN$ \emph{(``Braun-Meise-Taylor setting'').} Here, $\ell>0$ denotes the crucial parameter, hence analogously one has the \emph{Roumieu-type} $\mathcal{E}_{\{\omega\}}$ and the \emph{Beurling-type} $\mathcal{E}_{(\omega)}$. For obtaining the coincidence of both weight function settings a crucial \emph{Paley-Wiener-type theorem} has been shown, see \cite[Lemma 3.3, Prop. 3.4]{BraunMeiseTaylor90}.\vspace{6pt}

In each approach basic growth and regularity assumptions on the weights $\omega$ are required and, as it has been mentioned in the introduction in \cite{BraunMeiseTaylor90}, several crucial differences between these two settings can be recognized. For the convenience of the reader we recall now some definitions; see also Sections \ref{growthpropsection}, \ref{weightfctnotionsect} and \ref{examplesect} for more details and note that the abbreviations of the conditions have already been used in \cite{dissertation}.

\begin{definition}\label{defweightfct}
$\omega:[0,+\infty)\rightarrow[0,+\infty)$ is called a \emph{weight function} in the terminology of \cite[Sect. 2.1]{index}, \cite[Sect. 2.2]{sectorialextensions} and \cite[Sect. 2.2]{sectorialextensions1}, if $\omega$ is
\begin{itemize}
\item[$(*)$] continuous,

\item[$(*)$] non-decreasing,

\item[$(*)$] $\omega(0)=0$, and

\item[$(*)$] $\lim_{t\rightarrow+\infty}\omega(t)=+\infty$.
\end{itemize}
If $\omega$ satisfies in addition $\omega(t)=0$ for all $t\in[0,1]$, then we call $\omega$ \emph{normalized.}
\end{definition}

\emph{Note:} For some considerations an even more general setting is sufficient; occasionally in \cite{index} we have even worked without the assumptions continuity and $\omega(0)=0$.

\begin{definition}\label{defBBweightfct}
According to \cite[Def. 1.3.22]{Bjorck66} the function $\omega:[0,+\infty)\rightarrow[0,+\infty)$ is called a weight in the sense of \emph{Beurling-Bj\"{o}rk,} and write \emph{BB-weight} for short, if $\omega$ is
\begin{itemize}
\item[$(i)$] continuous,

\item[$(ii)$] $\omega(0)=0$,

\item[$(iii)$] $\omega$ satisfies \hyperlink{sub}{$(\omega_{\on{sub}})$}; i.e. $\omega(s+t)\le\omega(s)+\omega(t)$ for all $s,t\ge 0$.

\item[$(iv)$] $\omega$ satisfies \hyperlink{omnq}{$(\omega_{\on{nq}})$}; i.e. $\int_1^{\infty}\frac{\omega(t)}{t^2}dt<+\infty$.

\item[$(v)$] $\omega$ satisfies \hyperlink{om3w}{$(\omega_{3,w})$}; i.e. $\log(t)=O(\omega(t))$ as $t\rightarrow+\infty$.
\end{itemize}
\end{definition}

Closely related to BB-weights are weights in the sense by H.-J. Petzsche and D. Vogt in \cite{PetzscheVogt}:

\begin{definition}\label{defPVweightfct}
According to \cite{PetzscheVogt} the function $\omega:[0,+\infty)\rightarrow[0,+\infty)$ is called a weight in the sense of \emph{Petzsche-Vogt,} and write \emph{PV-weight} for short, if

\begin{itemize}
\item[$(1)$] $\omega$ is non-decreasing,

\item[$(2)$] $\omega$ satisfies \hyperlink{sub}{$(\omega_{\on{sub}})$} (abbreviated by $(\alpha)$ in \cite{PetzscheVogt}),

\item[$(3)$] $\omega$ satisfies \hyperlink{omnq}{$(\omega_{\on{nq}})$} (abbreviated by $(\beta)$),

\item[$(4)$] $\omega$ satisfies \hyperlink{om3w}{$(\omega_{3,w})$} (abbreviated by $(\gamma)$).
\end{itemize}
\emph{Note:} The differences between BB- and PV-weights are subtle.
\end{definition}

On the other hand in the Braun-Meise-Taylor setting the following set of functions is standard:

\begin{definition}\label{defBMTweightfct}
$\omega:[0,+\infty)\rightarrow[0,+\infty)$ is called a weight function in the sense of \emph{Braun-Meise-Taylor,} see \cite[Def. 1.1]{BraunMeiseTaylor90} and write \emph{BMT-weight} for short, if
\begin{itemize}
\item[$(*)$] $\omega$ is a weight function in the sense of Definition \ref{defweightfct},

\item[$(*)$] $\omega$ has \hyperlink{om1}{$(\omega_1)$}; i.e. $\omega(2t)=O(\omega(t))$ as $t\rightarrow+\infty$,

\item[$(*)$] $\omega$ has \hyperlink{om3}{$(\omega_3)$}; i.e. $\log(t)=o(\omega(t))$ as $t\rightarrow+\infty$,

\item[$(*)$] $\omega$ has \hyperlink{om4}{$(\omega_4)$}; i.e. $\varphi_{\omega}:t\mapsto\omega(e^t)$ is a convex function (on $\RR$).
\end{itemize}
\end{definition}

Summarizing, when comparing the basic requirements on the weights it follows that \emph{subadditivity} for $\omega$ in \cite{Beurling61}, \cite{Bjorck66} and \cite{PetzscheVogt} (Beurling-Bj\"{o}rck setting) is replaced resp. relaxed in the Braun-Meise-Taylor setting by $\omega(2t)=O(\omega(t))$ which is frequently abbreviated by $(\alpha)$ in the literature. A further achievement is that the definition of weighted spaces within the Braun-Meise-Taylor setting admits the possibility to introduce and investigate the notion of \emph{quasianalyticity} as well analogously to the weight sequence setting \emph{(``Denjoy-Carleman-Theorem'')} whereas within the Beurling-Bj\"{o}rck setting purely (non-trivial) test function spaces are studied; see Section \ref{examplesect} for precise definitions.\vspace{6pt}

On the other hand, for establishing the equivalent description in \cite[Lemma 3.3, Prop. 3.4]{BraunMeiseTaylor90} crucially it has been used that $\varphi_{\omega}$ is \emph{convex;} see \hyperlink{om4}{$(\omega_4)$} and this condition is frequently abbreviated by $(\delta)$ in the literature. \hyperlink{om4}{$(\omega_4)$} has become a standard assumption on $\omega$ also in related weighted settings since then; e.g. when considering globally defined functions of \emph{Gelfand-Shilov-type.} This convexity condition is unavoidable for obtaining the coincidence of both settings since the proof of the \emph{Paley-Wiener-type theorem} involves the double conjugate $\varphi^{**}_{\omega}$ and the identity $\varphi^{**}_{\omega}=\varphi_{\omega}$ which holds if and only if $\varphi_{\omega}$ is convex. In general, one only has $\varphi^{**}_{\omega}\le\varphi_{\omega}$ and for this note that $\varphi^{**}_{\omega}$ and $\varphi^{*}_{\omega}$ both are automatically convex; see the comments in \cite[Sect. 3.1]{dissertation} and more details are given in Section \ref{LFYconjgsection}.

In \cite{dissertation} and \cite{compositionpaper} we have shown that ultradifferentiable classes in the Braun-Meise-Taylor setting can equivalently be described by using the \emph{associated weight matrix} $\mathcal{M}_{\omega}:=\{\mathbf{W}^{(\ell)}=(W^{(\ell)}_j)_{j\in\NN}: \ell>0\}$, $W^{(\ell)}_j:=\exp(\frac{1}{\ell}\varphi^{*}_{\omega}(\ell j))$; we refer to the main statement \cite[Thm. 5.14]{compositionpaper} and for this alternative description \hyperlink{om1}{$(\omega_1)$} is crucial. It is immediate by definition that properties for $\omega$ can be transferred to $\mathcal{M}_{\omega}$ but, however, for the converse and hence when establishing a characterization of growth conditions, then again \hyperlink{om4}{$(\omega_4)$} is becoming relevant in order to treat the double conjugate. Similarly, this comment applies when properties for the classes $\mathcal{E}_{\{\omega\}}$, $\mathcal{E}_{(\omega)}$ are transferred to $\omega$ (resp. $\varphi_{\omega}$); see e.g. the proof of \cite[Cor. 5.17]{compositionpaper}. Another crucial property for which \hyperlink{om4}{$(\omega_4)$} is unavoidable is the equivalence between $\omega$ and each \emph{associated weight function} $\omega_{\mathbf{W}^{(\ell)}}$; see \cite[Lemma 5.7]{compositionpaper}, \cite[Theorem 4.0.3, Lemma 5.1.3]{dissertation}, and also \cite[Lemma 2.5]{sectorialextensions}. On the other hand, for the sake of completeness we mention that basically $\mathcal{M}_{\omega}$ can be introduced for more general $\omega$. Indeed, it suffices to assume that $\omega$ is non-decreasing and $\log(t)=o(\omega(t))$ as $t\rightarrow+\infty$ (i.e. \hyperlink{om3}{$(\omega_3)$}) in order to ensure well-definedness of $\mathcal{M}_{\omega}$ \emph{(``matrix admissible weights'');} see Section \ref{LFYconjgsection} for more details.\vspace{6pt}

The main aim of this article is to go one step further, to define new weighted $L^p$-type spaces both of Roumieu- and Beurling-type in terms of abstractly given \emph{weight function matrices} $\mathcal{W}:=\{\omega^{\ell}: \ell>0\}$ with matrix parameter $\ell>0$ (see Definition \ref{functionmatrixdef}) and to answer the following question; see Sections \ref{generalclasssect}, \ref{charactinclsection}, \ref{specialcasesect} and \ref{BBclasssect}:
\begin{itemize}
\item[$(Q1)$] Characterize the inclusion relations of weighted $L^p$-type spaces, and hence the equality (as l.c.v.s.), in terms of the defining matrices.
\end{itemize}
We answer $(Q1)$ in the main statements Theorem \ref{maincharthmbeur} for the Beurling-type and Theorem \ref{maincharthmroum} for the Roumieu-type and work in a very general setting. First, we deal with abstractly given matrices $\mathcal{W}$ and, second, in Sections \ref{generalclasssect}, \ref{charactinclsection} we are concerned with weighted subspaces of $L^{p}(\RR^d)$, $1\le p\le\infty$, subject to $\mathcal{W}$. Indeed, the weights are given by $x\mapsto\exp(\omega^{\ell}(x))$ and so involving the exponential function. Section \ref{specialcasesect} summarizes the situation when we focus on two special but natural weight function matrices for which the crucial information can be expressed in terms of a single weight function $\omega$: matrices of \emph{``exponential-type''} and of \emph{``dilatation-type''.} The special structure of these matrices allows for a more compact characterization. In Section \ref{BBclasssect} we apply the information to related and analogous situations: First, in Section \ref{subspacesect} we consider the setting when weighting a fixed \emph{translation invariant} subspace $\mathcal{A}^p\subseteq L^p(\RR^d)$. In Section \ref{weightedFouriersection} we focus in detail on the case when weighting the ``Fourier image'' of a given (fixed) \emph{modulation invariant} function space $\mathcal{A}$; i.e. not the growth of $f$ directly is controlled but of its Fourier transform $\widehat{f}$. Of course, for this case it is crucial to ensure that the Fourier transform is well-defined on $\mathcal{A}$ and its image is non-trivial. Finally, in each instance we also study the \emph{anisotropic setting;} i.e. when the weight functions are not radially extended to $\RR^d$.\vspace{6pt}

Summarizing, this general framework allows to treat in particular the known Beurling-Bj\"{o}rck setting (for both types) which is  corresponding to the special weight matrix $\omega^{\ell}:=\ell\omega$, $\ell>0$, where $\omega$ is given and fixed \emph{(``exponential type''),} see Section \ref{specialcasesect1}, and when weighting the Fourier image of fixed subspaces (of functions having compact support) in $L^1(\RR^d)$. We refer to Section \ref{examplesect} for more details. Recall that $(Q1)$ for the Braun-Meise-Taylor setting has been solved in \cite[Lemma 5.16 \& Cor. 5.17]{compositionpaper} and the crucial growth relations are expressed by an O-growth resp. $o$-growth restriction and thus equivalence of weight functions characterizes the equality (as l.c.v.s.); see Section \ref{growthrelsection} for the relevant definitions. As expected, for the Beurling-Bj\"{o}rck setting the same growth relations are becoming relevant.

In addition, we provide a detailed study concerning the coincidence of the new weighted spaces with the common (non-weighted) $L^p$-spaces (see Section \ref{unboundedsect}), of their (non-)triviality (see Section \ref{nontrivialsect}) and investigate their translation invariance in Section \ref{stabtranslationsect}.\vspace{6pt}

By taking into account the information from $(Q1)$ and the fact that \hyperlink{om4}{$(\omega_4)$} is in general not preserved under equivalence of weight functions, we can see the second main question:

\begin{itemize}
\item[$(Q2)$] Illustrate the difference between the Beurling-Bj\"{o}rck and the Braun-Meise-Taylor setting in view of the convexity for $\varphi_{\omega}$; more precisely we ask: How restrictive is \hyperlink{om4}{$(\omega_4)$}? Does for any weight function (in the sense of Beurling-Bj\"{o}rck) an equivalent function exist which satisfies \hyperlink{om4}{$(\omega_4)$}?
\end{itemize}

Without \hyperlink{om4}{$(\omega_4)$} the main Paley-Wiener-type result from \cite{BraunMeiseTaylor90} and hence the equality between both settings is becoming unclear. Therefore, when treating the pure Beurling-Bj\"{o}rck framework in view of $(Q1)$ the specific role resp. the loss of \hyperlink{om4}{$(\omega_4)$} has to be investigated. This question has been motivated by the explanations above and it is based on the comments in \cite[Sect. 3.1]{dissertation}. But also in the recent article \cite{NeytVindas23} in the introduction this convexity condition, its possible failure and its role in proofs within the Beurling-Bj\"{o}rck are discussed. (\hyperlink{om4}{$(\omega_4)$} is also denoted by $(\delta)$ in \cite{NeytVindas23}.)\vspace{6pt}

The last Section \ref{failuresection} is dedicated to $(Q2)$ and to the construction of a (counter-)example showing that there exist weights in the sense of Beurling-Bj\"{o}rck (and of Petzsche-Vogt) such that none equivalent (weight) function satisfies \hyperlink{om4}{$(\omega_4)$}. Indeed, we even verify that there exist ``many'' of such weight functions; i.e. uncountable infinitely many different equivalence classes. The geometric construction also involves recent knowledge about the growth index $\gamma(\omega)$ for weight functions $\omega$ introduced in \cite{index}; see Section \ref{growthindexsection}. An immediate consequence of this example is to see the limitation of the weight sequence setting in this approach; i.e. when considering \emph{associated weight functions} $\omega_{\mathbf{M}}$ with $\mathbf{M}$ being a sequence satisfying basic growth and regularity properties. This is due to the fact that $\omega_{\mathbf{M}}$ satisfies automatically \hyperlink{om4}{$(\omega_4)$}; see Remark \ref{nonassorem} for more details. On the other hand $\omega_{\mathbf{M}}$ appears frequently in the Braun-Meise-Taylor setting as a (counter-)example and in the comparison results from \cite{BonetMeiseMelikhov07}; this again illustrates the difference between both settings. Due to the failure of \hyperlink{om4}{$(\omega_4)$} this (counter-)example also implies the fact that the information concerning the characterization of inclusion relations applies to spaces for which the corresponding BMT-setting results \cite[Lemma 5.16 \& Cor. 5.17]{compositionpaper} cannot be used. Finally, we emphasize that the constructed $\omega$ in Section \ref{failuresection} is matrix admissible.\vspace{6pt}

In the first Section \ref{weightfctsect} we list, summarize and compare in detail the numerous growth conditions and relations on weight functions in the different frameworks. One can expect that this content is also interesting for other different weighted settings since we exclusively deal with growth and regularity properties for functions. The same comment applies to the construction provided in Section \ref{failuresection}.\vspace{6pt}

\textbf{Acknowledgements.} The author of this article thanks the two anonymous referees for their efforts, the careful reading and the valuable suggestions which have improved the presentation of the results.

\section{Weight functions and their growth conditions}\label{weightfctsect}

\subsection{General notation}
We write $\NN:=\{0,1,2,\dots\}$, $\NN_{>0}:=\{1,2,\dots\}$ and for $x\in\RR^d$ the expression $|x|$ denotes the usual Euclidean norm on $\RR^d$.

\subsection{Relevant growth properties}\label{growthpropsection}
We list (known) conditions for functions $\omega:[0,+\infty)\rightarrow[0,+\infty)$; the abbreviations are also appearing in \cite{dissertation}. Some of these properties appear exclusively for weights in the sense of Beurling-Bj\"{o}rck (see \cite{Bjorck66}), some of them only for Braun-Meise-Taylor weights (see \cite{BraunMeiseTaylor90}) and some of them are crucial in both settings.\vspace{6pt}

We call a function $\omega$ \emph{normalized} if $\omega(t)=0$ for all $0\le t\le 1$. Moreover, we shall assume $\omega\neq 0$ except stated explicitly otherwise. Inspired by the notation used in \cite{dissertation} we consider:\vspace{6pt}

\begin{itemize}
\item[\hypertarget{om1}{$(\omega_1)}$] $\omega(2t)=O(\omega(t))$ as $t\rightarrow+\infty$.
	
\item[\hypertarget{om2}{$(\omega_2)$}] $\omega(t)=O(t)$ as $t\rightarrow+\infty$.
	
\item[\hypertarget{om3}{$(\omega_3)$}] $\log(t)=o(\omega(t))$ as $t\rightarrow+\infty$.
	
\item[\hypertarget{om4}{$(\omega_4)$}] $\varphi_{\omega}:t\mapsto\omega(e^t)$ is a convex function (on $\RR$).
	
\item[\hypertarget{om5}{$(\omega_5)$}] $\omega(t)=o(t)$ as $t\rightarrow+\infty$.

\item[\hypertarget{om6}{$(\omega_6)$}] $\exists\;H\ge 1\;\forall\;t\ge 0:\;2\omega(t)\le\omega(H t)+H$.

\item[\hypertarget{omnq}{$(\omega_{\text{nq}})$}] $\int_1^{\infty}\frac{\omega(t)}{t^2}dt<+\infty$.
	
\item[\hypertarget{omsnq}{$(\omega_{\text{snq}})$}] $\exists\;C>0\;\forall\;y>0: \int_1^{\infty}\frac{\omega(y t)}{t^2}dt\le C\omega(y)+C$.
\end{itemize}

The next condition is crucial for characterizing desired stability properties for Braun-Meise-Taylor classes $\mathcal{E}_{\{\omega\}}$ and $\mathcal{E}_{(\omega)}$ (of both types), see \cite{compositionpaper}, \cite{characterizationstabilitypaper} and the citations therein:
\begin{itemize}
\item[\hypertarget{alpha0}{$(\alpha_0)$}] $\exists\;C\ge 1\;\exists\;t_0\ge 0\;\forall\;\lambda\ge 1\;\forall\;t\ge t_0:\;\;\;\omega(\lambda t)\le C\lambda\omega(t).$
\end{itemize}
In \cite[$(3.8.4)$]{dissertation} this requirement has been denoted by $(\omega_{1^{'}})$ but frequently in the literature it is abbreviated by $(\alpha_0)$. Next, let us introduce the following two requirements appearing in the Beurling-Bj\"{o}rck framework (see \cite{Bjorck66}):

\begin{itemize}
\item[\hypertarget{om3w}{$(\omega_{3,w})$}] $\log(t)=O(\omega(t))$ as $t\rightarrow+\infty$.

\item[\hypertarget{sub}{$(\omega_{\text{sub}})$}] $\omega(s+t)\le\omega(s)+\omega(t)$ for all $s,t\ge 0$.
\end{itemize}
Finally, recall that the \emph{Legendre-Fenchel-Young conjugate (of $\varphi_{\omega}$)} is defined as follows (see \eqref{LFYconj}):
\begin{equation*}
\varphi^{*}_{\omega}(x):=\sup_{y\in\RR}\{xy-\varphi_{\omega}(y)\},\;x\in\RR.
\end{equation*}
This transform is crucial when defining weighted classes in the sense of Braun-Meise-Taylor and it is also basic for introducing the \emph{associated weight matrix} $\mathcal{M}_{\omega}:=\{\mathbf{W}^{(\ell)}=(W^{(\ell)}_j)_{j\in\NN}: \ell>0\}$, $W^{(\ell)}_j:=\exp\left(\frac{1}{\ell}\varphi^{*}_{\omega}(\ell j)\right)$. We investigate this transform in more detail in Section \ref{LFYconjgsection}.

\subsection{Relevant growth relations between weight functions}\label{growthrelsection}
Let $\sigma,\tau: [0,+\infty)\rightarrow[0,+\infty)$ be arbitrary functions and introduce the following growth relations between $\sigma$ and $\tau$:

\begin{itemize}
\item[$(*)$] Write $\sigma\le\tau$ if $\sigma(t)\le\tau(t)$ for all $t\in[0,+\infty)$.

\item[$(*)$] Write $\sigma\hypertarget{ompreceq}{\preceq}\tau$ if
$$\exists\;C\ge 1\;\forall\;t\ge 0:\;\;\;\tau(t)\le C\sigma(t)+C,$$
and $\sigma$ and $\tau$ are called \emph{equivalent,} written $\sigma\hypertarget{sim}{\sim}\tau$, if $\sigma\hyperlink{ompreceq}{\preceq}\tau$ and $\tau\hyperlink{ompreceq}{\preceq}\sigma$. This is the same relation as it has been considered in \cite[Sect. 2.1]{index}.

\item[$(*)$] Write $\sigma\hypertarget{omtriangle}{\vartriangleleft}\tau$, if
$$\forall\;\epsilon>0\;\exists\;C\ge 1\;\forall\;t\ge 0:\;\;\;\tau(t)\le\epsilon\sigma(t)+C.$$
\end{itemize}

Obviously, $\sigma\hyperlink{omtriangle}{\vartriangleleft}\tau$ implies $\sigma\hyperlink{ompreceq}{\preceq}\tau$ but \hyperlink{omtriangle}{$\vartriangleleft$} is in general not reflexive: $\sigma\hyperlink{omtriangle}{\vartriangleleft}\sigma$ holds if and only if $\sup_{t\ge 0}\sigma(t)<+\infty$.

\emph{Note:}
\begin{itemize}
\item[$(i)$] In the literature, frequently the crucial growth relations for weight functions read as follows: $\sigma\preceq\tau$ resp. $\sigma\vartriangleleft\tau$ precisely means
\begin{equation}\label{bigOrelation}
\tau(t)=O(\sigma(t))\;\text{as}\;t\rightarrow+\infty,\;\;\;\text{resp.}\;\;\;\tau(t)=o(\sigma(t))\;\text{as}\;t\rightarrow+\infty.
\end{equation}
For these, in general stronger requirements, we also refer e.g. to \cite[Rem. 2.6, $(5)$]{index}.

\item[$(ii)$] If $\sigma$ satisfies $\liminf_{t\rightarrow+\infty}\sigma(t)>0$; i.e.
\begin{equation}\label{liminfweighfct}
\exists\;t_0\ge 0\;\exists\;\epsilon>0\;\forall\;t\ge t_0:\;\;\;\sigma(t)\ge\epsilon,
\end{equation}
then $\sigma\hyperlink{ompreceq}{\preceq}\tau$ if and only $\tau(t)=O(\sigma(t))$ as $t\rightarrow+\infty$ and $\sigma\hyperlink{omtriangle}{\vartriangleleft}\tau$ if and only if $\tau(t)=o(\sigma(t))$ as $t\rightarrow+\infty$.

\item[$(iii)$] \eqref{liminfweighfct} is a mild requirement; e.g. it holds if $\sigma\neq 0$ is assumed to be non-decreasing. Moreover, it holds if $\lim_{t\rightarrow+\infty}\sigma(t)=+\infty$ and this has been a basic and natural growth requirement when studying weight functions even in a very general setting in \cite{index}. And, finally, we mention that to each arbitrary function $\sigma:[0,+\infty)\rightarrow[0,+\infty)$ and $\epsilon>0$ there exists $\sigma_{\epsilon}$ satisfying \eqref{liminfweighfct} for this $\epsilon>0$ and with $t_0=0$ and such that $\sigma\hyperlink{sim}{\sim}\sigma_{\epsilon}$: We simply put $\sigma_{\epsilon}(t):=\sigma(t)+\epsilon$ and get $\sigma_{\epsilon}(t)\ge\epsilon$, $\sigma(t)\le\sigma_{\epsilon}(t)\le C\sigma(t)+C$ for all $t\ge 0$ with $C:=\max\{1,\epsilon\}$.

\item[$(iv)$] If $\sigma:[0,+\infty)\rightarrow(0,+\infty)$, so if $\sigma(t)>0$ for all $t\ge 0$, then $\tau(t)=O(\sigma(t))$ as $t\rightarrow+\infty$ precisely means $\tau(t)\le C\sigma(t)$ for some $C\ge 1$ and all $t\ge 0$; see e.g. the weighted entire setting treated in \cite{weightedentireinclusion1}, \cite{weightedentireinclusion2}.

\item[$(v)$] In view of $(ii)$, for any $\sigma$ satisfying \eqref{liminfweighfct} we get that \hyperlink{om1}{$(\omega_1)$} is equivalent to $\omega(2t)\le L\omega(t)+L$ for some $L\ge 1$ and all $t\ge 0$ and this estimate is appropriate even if \eqref{liminfweighfct} fails.

    On the other hand, for \hyperlink{om2}{$(\omega_2)$} and \hyperlink{om5}{$(\omega_5)$} this technical observation is not required since $t\mapsto t$ is obviously satisfying \eqref{liminfweighfct}.
\end{itemize}

All properties listed in Section \ref{growthpropsection} except \hyperlink{sub}{$(\omega_{\on{sub}})$} and the convexity property \hyperlink{om4}{$(\omega_4)$} are automatically preserved under \hyperlink{sim}{$\sim$}; see also \cite[Sect. 3.2]{dissertation}. Indeed, the effects on \hyperlink{sub}{$(\omega_{\on{sub}})$} when passing to an equivalent weight and the relation to \hyperlink{alpha0}{$(\alpha_0)$} are investigated and revisited in Section \ref{alphazerosubsection} whereas the failure of \hyperlink{om4}{$(\omega_4)$} is studied in detail in Section \ref{failuresection} where a (counter-)example is constructed.\vspace{6pt}

Naturally, one can also consider the following growth relations:

\begin{itemize}
\item[$(*)$] Write $\sigma\hypertarget{ompreceqc}{\preceq_{\mathfrak{c}}}\tau$ if
$$\exists\;C_1,C_2\ge 1\;\forall\;t\ge 0:\;\;\;\tau(t)\le\sigma(C_1t)+C_2,$$
and write $\sigma\hypertarget{simc}{\sim_{\mathfrak{c}}}\tau$, if $\sigma\hyperlink{ompreceqc}{\preceq_{\mathfrak{c}}}\tau$ and $\tau\hyperlink{ompreceqc}{\preceq_{\mathfrak{c}}}\sigma$.

\item[$(*)$] Write $\sigma\hypertarget{omtrianglec}{\vartriangleleft_{\mathfrak{c}}}\tau$, if
$$\forall\;\epsilon>0\;\exists\;C\ge 1\;\forall\;t\ge 0:\;\;\;\tau(t)\le\sigma(\epsilon t)+C.$$
\end{itemize}
If $\sigma\hyperlink{ompreceqc}{\preceq_{\mathfrak{c}}}\tau$ and either $\tau$ or $\sigma$ is non-decreasing, then in this relation we can assume w.l.o.g. $C_1=C_2$ by taking $C:=\max\{C_1,C_2\}$. \hyperlink{omvartrianglec}{$\vartriangleleft_{\mathfrak{c}}$} implies \hyperlink{ompreceqc}{$\preceq_{\mathfrak{c}}$} but \hyperlink{omvartrianglec}{$\vartriangleleft_{\mathfrak{c}}$} is in general not reflexive: If $\sigma$ is non-decreasing and $\lim_{t\rightarrow+\infty}\sigma(t)=+\infty$, then $\sigma\hyperlink{omvartrianglec}{\vartriangleleft_{\mathfrak{c}}}\sigma$ implies
\begin{equation}\label{slowlyvarying}
\forall\;u>0:\;\;\;\lim_{t\rightarrow+\infty}\frac{\sigma(tu)}{\sigma(t)}=1;
\end{equation}
i.e. $\sigma$ has to be \emph{slowly varying,} see \cite[$(1.2.1)$]{regularvariation} and \cite[Sect. 2.2]{index}: $\sigma\hyperlink{omvartrianglec}{\vartriangleleft_{\mathfrak{c}}}\sigma$ yields $\limsup_{t\rightarrow+\infty}\frac{\sigma(tu)}{\sigma(t)}\le 1$ for all $u>0$ and since $\sigma$ is non-decreasing $\lim_{t\rightarrow+\infty}\frac{\sigma(tu)}{\sigma(t)}=1$ follows for any $u\ge 1$. And for $0<u<1$ set $s:=ut$ and get $\lim_{t\rightarrow+\infty}\frac{\sigma(tu)}{\sigma(t)}=\lim_{s\rightarrow+\infty}\frac{\sigma(s)}{\sigma(s/u)}=1$. (This argument implies that it suffices to check \eqref{slowlyvarying} for all $u>1$.) Indeed, we recognize that $\sigma\hyperlink{omtrianglec}{\vartriangleleft_{\mathfrak{c}}}\tau$ is precisely \cite[Sect. 4.2 $(B)$]{genLegendreconj} between $\sigma$ and $\tau$; i.e. \cite[$(4.4)$]{genLegendreconj} with $t_0=+\infty$ and hence this relation is becoming (also) relevant for ensuring well-definedness of the operation $\widehat{\star}$ introduced and studied in \cite[Sect. 4]{genLegendreconj}: $\tau\hyperlink{omvartrianglec}{\vartriangleleft_{\mathfrak{c}}}\sigma$ if and only if $\sigma\widehat{\star}\tau$ is well-defined on $[0,+\infty)$; and the above arguments concerning \eqref{slowlyvarying} are also contained in \cite[Rem. 4.4]{genLegendreconj}.

\subsection{On the comparison of crucial growth properties}\label{relevantgrowthsect}
Let $\omega:[0,+\infty)\rightarrow[0,+\infty)$ be given. We gather now immediate observations for the conditions listed in Section \ref{growthpropsection} and study relations between them.

\begin{itemize}
\item[$(a)$] \hyperlink{om3}{$(\omega_3)$} is equivalent to $\lim_{t\rightarrow+\infty}\frac{\varphi_{\omega}(t)}{t}=+\infty$ and it clearly implies $\lim_{t\rightarrow+\infty}\omega(t)=+\infty$. Moreover, \hyperlink{om3}{$(\omega_3)$} is equivalent to $\log(a+t)=o(\omega(t))$ as $t\rightarrow+\infty$ and \hyperlink{om3w}{$(\omega_{3,w})$} is equivalent to $\log(a+t)=O(\omega(t))$ as $t\rightarrow+\infty$, where $a\ge 1$ is an arbitrary real parameter.

\item[$(b)$] \hyperlink{om3}{$(\omega_3)$} implies \hyperlink{om3w}{$(\omega_{3,w})$} and \hyperlink{om3w}{$(\omega_{3,w})$} implies $\lim_{t\rightarrow+\infty}\omega(t)=+\infty$. The crucial difference between these two conditions is that \hyperlink{om3w}{$(\omega_{3,w})$} formally allows to consider the ``limiting weight'' $\omega(t):=\log(1+t)$ whereas \hyperlink{om3}{$(\omega_3)$} excludes this special case. In the literature occasionally $t\mapsto\log(1+t)$ is treated for the  Beurling-type and which gives back indeed the non-weighted situation.

\item[$(c)$] \hyperlink{omsnq}{$(\omega_{\on{snq}})$} implies \hyperlink{omnq}{$(\omega_{\on{nq}})$}.

\item[$(d)$] If $\omega$ satisfies \eqref{liminfweighfct}, then \hyperlink{alpha0}{$(\alpha_0)$} implies \hyperlink{om1}{$(\omega_1)$} (choose $\lambda=2$).

\item[$(e)$] If $\omega$ satisfies \eqref{liminfweighfct}, then also \hyperlink{sub}{$(\omega_{\on{sub}})$} implies \hyperlink{om1}{$(\omega_1)$} (take $x=y$) and with $x=y=0$ we get $\omega(0)\le 2\omega(0)$ which yields $\omega(0)\ge 0$.

\item[$(f)$] In the literature sometimes instead of \hyperlink{om1}{$(\omega_1)$} the following condition is used (see e.g. \cite[Lemma 1.2]{BraunMeiseTaylor90}):
    \begin{equation}\label{om1reform}
    \exists\;L\ge 1\;\forall\;s,t\ge 0:\;\;\;\omega(s+t)\le L(\omega(s)+\omega(t))+L.
    \end{equation}
     When \eqref{liminfweighfct} holds, then set $s=t$ and get immediately \hyperlink{om1}{$(\omega_1)$} when choosing the constant $2L$. If $\omega$ is non-decreasing, then \eqref{om1reform} is equivalent to \hyperlink{om1}{$(\omega_1)$}: Indeed, if $\omega$ has \hyperlink{om1}{$(\omega_1)$} then $\omega(s+t)\le\omega(2\max\{s,t\})\le L\omega(\max\{s,t\})+L\le L(\omega(s)+\omega(t))+L$ and so \eqref{om1reform} holds with the same $L$.

\item[$(g)$] \hyperlink{om5}{$(\omega_5)$} implies \hyperlink{om2}{$(\omega_2)$}.

\item[$(h)$] If $\omega$ is non-decreasing, then \hyperlink{omnq}{$(\omega_{\on{nq}})$} implies \hyperlink{om5}{$(\omega_5)$} because $\int_s^{\infty}\frac{\omega(t)}{t^2}dt\ge\omega(s)\int_s^{\infty}\frac{1}{t^2}dt=\frac{\omega(s)}{s}$ for all $s\ge 1$.

\item[$(i)$] Let $\omega$ be non-decreasing (and $\omega\neq 0$) and assume that there exists some $t_0>0$ such that $\omega(t)=0$ for all $t\in[0,t_0]$. In particular, this last assumption holds when $\omega$ is normalized; i.e. $\omega(t)=0$ for $t\in[0,1]$. Then several growth conditions are violated:

    \begin{itemize}
    \item[$(a)$] Concavity fails: This condition means $\omega(\lambda t+(1-\lambda)s)\ge\lambda\omega(t)+(1-\lambda)\omega(s)$ for all $s,t\ge 0$ and $\lambda\in[0,1]$. Choose $s\in[0,t_0)$, $t\in(t_0,+\infty)$ with $\omega(t)>0$ and such $t$ exists since $\omega$ is non-decreasing and $\omega\neq 0$. Then there exists $\lambda_0\in(0,1)$ such that $\lambda_0 s+(1-\lambda_0)t=t_0$ and thus $0=\omega(t_0)\ge(1-\lambda_0)\omega(t)>0$ follows, a contradiction.

\item[$(b)$] \hyperlink{alpha0}{$(\alpha_0)$} fails (see also \cite[Sect. 4.1]{subaddlike} and the references therein): $C\lambda\omega(t)=0$ holds for all $C,\lambda\ge 1$ and all $0<t\le t_0$ but, since $\omega\neq 0$ and $\omega$ is non-decreasing, we have $\omega(\lambda t)>0$ as $\lambda\rightarrow+\infty$ for any $t>0$ fixed.

\item[$(c)$] \hyperlink{sub}{$(\omega_{\on{sub}})$} fails: For all $0<s\le t_0$ and all $t\ge 0$ we get $\omega(t)\le\omega(s+t)\le\omega(s)+\omega(t)=\omega(t)$ which implies $\omega=0$; again a contradiction.
\end{itemize}
Indeed, $(a)$ and $(b)$ are even valid if \eqref{liminfweighfct} and $\omega(t)=0$ for all $t\in[0,t_0]$.
\end{itemize}

We close this section by investigating conditions \hyperlink{om1}{$(\omega_1)$} and \hyperlink{om6}{$(\omega_6)$}. First, we state a technical result.

\begin{lemma}\label{om1om6itlemma}
Let $\omega: [0,+\infty)\rightarrow[0,+\infty)$ be non-decreasing.
\begin{itemize}
\item[$(I)$] The following are equivalent:
\begin{itemize}
\item[$(i)$] $\omega$ satisfies \hyperlink{om1}{$(\omega_1)$}.

\item[$(ii)$] $\omega$ satisfies
$$\exists\;a>1\;\exists\;L\ge 1\;\forall\;t\ge 0:\;\;\;\omega(at)\le L\omega(t)+L.$$

\item[$(iii)$] $\omega$ satisfies
$$\forall\;a>1\;\exists\;L\ge 1\;\forall\;t\ge 0:\;\;\;\omega(at)\le L\omega(t)+L.$$
\end{itemize}

\item[$(II)$] The following are equivalent:

\begin{itemize}
\item[$(i)$] $\omega$ satisfies \hyperlink{om6}{$(\omega_6)$}.

\item[$(ii)$] $\omega$ satisfies
$$\exists\;a>1\;\exists\;H\ge 1\;\forall\;t\ge 0:\;\;\;a\omega(t)\le\omega(Ht)+H.$$

\item[$(iii)$] $\omega$ satisfies
$$\forall\;a>1\;\exists\;H\ge 1\;\forall\;t\ge 0:\;\;\;a\omega(t)\le\omega(Ht)+H.$$
\end{itemize}
\end{itemize}
\end{lemma}

Note that assumption \cite[$(\alpha)$]{Franken94} precisely corresponds to $(I)(ii)$.

\demo{Proof}
The particular implications $(i)\Rightarrow(ii)$ and $(iii)\Rightarrow(i)$ are clear (set $a:=2$) and $(ii)\Rightarrow(iii)$ follows by iteration: Let $b>1$ be arbitrary, then choose $n\in\NN_{>0}$ such that $a^n\ge b$ with $a$ appearing in $(ii)$, and iterate the corresponding estimate in $(ii)$ $n$-times. (Concerning \hyperlink{om6}{$(\omega_6)$} see also \cite[Rem. 2.2]{modgrowthstrangeII} but, however, the assumption $\lim_{t\rightarrow+\infty}\omega(t)=+\infty$ made there is not required necessarily for the proof of these equivalences.)
\qed\enddemo

By taking into account this result, conditions \hyperlink{om1}{$(\omega_1)$} and \hyperlink{om6}{$(\omega_6)$} give the possibility to compare the relations \hyperlink{ompreceq}{$\preceq$} and \hyperlink{ompreceqc}{$\preceq_{\mathfrak{c}}$} resp. \hyperlink{omtriangle}{$\vartriangleleft$} and \hyperlink{omtrianglec}{$\vartriangleleft_{\mathfrak{c}}$}; confirm also \cite[Prop. 5.1]{weightedentireinclusion2} in the weighted entire setting.

\begin{proposition}\label{weightcomparisonprop}
Let $\sigma,\tau:[0,+\infty)\rightarrow[0,+\infty)$ be non-decreasing.
\begin{itemize}
\item[$(i)$] If either $\sigma$ or $\tau$ satisfies \hyperlink{om6}{$(\omega_6)$}, then $\sigma\hyperlink{ompreceq}{\preceq}\tau$ implies $\sigma\hyperlink{ompreceqc}{\preceq_{\mathfrak{c}}}\tau$.

\item[$(ii)$] If either $\sigma$ or $\tau$ satisfies \hyperlink{om1}{$(\omega_1)$}, then $\sigma\hyperlink{ompreceqc}{\preceq_{\mathfrak{c}}}\tau$ implies $\sigma\hyperlink{ompreceq}{\preceq}\tau$.

\item[$(iii)$] If either $\sigma$ or $\tau$ satisfies \hyperlink{om1}{$(\omega_1)$}, then $\sigma\hyperlink{omtriangle}{\vartriangleleft}\tau$ implies $\sigma\hyperlink{omtrianglec}{\vartriangleleft_{\mathfrak{c}}}\tau$.

\item[$(iv)$] If either $\sigma$ or $\tau$ satisfies \hyperlink{om6}{$(\omega_6)$}, then $\sigma\hyperlink{omtrianglec}{\vartriangleleft_{\mathfrak{c}}}\tau$ implies $\sigma\hyperlink{omtriangle}{\vartriangleleft}\tau$.
\end{itemize}
\end{proposition}

\demo{Proof}
$(i)$ Let $\sigma\hyperlink{ompreceq}{\preceq}\tau$, then $\tau(t)\le a\sigma(t)+a$ for some $a\ge 1$ and all $t\ge 0$. If $\sigma$ satisfies \hyperlink{om6}{$(\omega_6)$}, then $(II)$ in Lemma \ref{om1om6itlemma} yields the conclusion. If $\tau$ satisfies \hyperlink{om6}{$(\omega_6)$}, then we have the estimate $\tau(t)\le\frac{1}{a}\tau(Ht)+\frac{H}{a}\le\sigma(Ht)+\frac{H}{a}+1$ for some $H\ge 1$ (depending on $a$) and all $t\ge 0$, hence again $\sigma\hyperlink{ompreceqc}{\preceq_{\mathfrak{c}}}\tau$.

$(ii)$ Let $\sigma\hyperlink{ompreceqc}{\preceq_{\mathfrak{c}}}\tau$, then $\tau(t)\le\sigma(at)+a$ for some $a\ge 1$ and all $t\ge 0$. If $\sigma$ satisfies \hyperlink{om1}{$(\omega_1)$}, then $(I)$ in Lemma \ref{om1om6itlemma} yields the conclusion. If $\tau$ satisfies \hyperlink{om1}{$(\omega_1)$}, then we have the estimate $\tau(t)\le L\tau(a^{-1}t)+L\le L\sigma(t)+La+L$ for some $L\ge 1$ depending on $a$ and all $t\ge 0$.\vspace{6pt}

$(iii)$ Let $\sigma\hyperlink{omtriangle}{\vartriangleleft}\tau$ and let $\epsilon>0$ be given (small). If $\sigma$ satisfies \hyperlink{om1}{$(\omega_1)$}, then $\sigma(t)\le L\sigma(\epsilon t)+L$ for some $L\ge 1$ depending on given $\epsilon$ and all $t\ge 0$. Then apply relation $\sigma\hyperlink{omtriangle}{\vartriangleleft}\tau$ to $L^{-1}$ (small) and get $\tau(t)\le L^{-1}\sigma(t)+C\le\sigma(\epsilon t)+1+C$ for some $C$ depending on $\epsilon$ (via $L$) and all $t\ge 0$. If $\tau$ satisfies \hyperlink{om1}{$(\omega_1)$}, then $\tau(t)\le L\tau(\epsilon t)+L\le\sigma(\epsilon t)+CL+L$ again for some $C$ and $L$ both depending on $\epsilon$ and all $t\ge 0$ when applying $\sigma\hyperlink{omtriangle}{\vartriangleleft}\tau$ to $L^{-1}$. Therefore, in both cases $\sigma\hyperlink{omtrianglec}{\vartriangleleft_{\mathfrak{c}}}\tau$ is valid.\vspace{6pt}

$(iv)$ Let $\sigma\hyperlink{omtrianglec}{\vartriangleleft_{\mathfrak{c}}}\tau$ be valid and let $\epsilon>0$ be given (small). If $\sigma$ satisfies \hyperlink{om6}{$(\omega_6)$}, then $\epsilon^{-1}\sigma(t)\le\sigma(Ht)+H$ for some $H\ge 1$ depending on given $\epsilon$ and all $t\ge 0$. Then apply $\sigma\hyperlink{omtrianglec}{\vartriangleleft_{\mathfrak{c}}}\tau$ to $H^{-1}$ (small) and get $\tau(t)\le\sigma(H^{-1}t)+C\le\epsilon\sigma(t)+\epsilon H+C$ for some $H,C\ge 1$ both depending on $\epsilon$ and all $t\ge 0$. If $\tau$ satisfies \hyperlink{om6}{$(\omega_6)$}, then $\epsilon^{-1}\tau(t)\le\tau(Ht)+H$ for some $H\ge 1$ depending on given $\epsilon$ and all $t\ge 0$ and so $\tau(t)\le\epsilon\tau(Ht)+\epsilon H\le\epsilon\sigma(t)+\epsilon C+\epsilon H$ for some $C,H\ge 1$ both depending on $\epsilon$ and all $t\ge 0$ when applying again $\sigma\hyperlink{omtrianglec}{\vartriangleleft_{\mathfrak{c}}}\tau$ to $H^{-1}$. Summarizing, in both cases $\sigma\hyperlink{omtriangle}{\vartriangleleft}\tau$ is verified.
\qed\enddemo

\subsection{On subadditivity-like conditions}\label{alphazerosubsection}
We recall the following equivalent reformulations of \hyperlink{alpha0}{$(\alpha_0)$}; see \cite[Lemma 1]{peetre}, \cite[Prop. 2.23]{index} and finally \cite[Lemma 4.1]{subaddlike}:

\begin{lemma}\label{alpha0equivalence}
Let $\omega:[0,+\infty)\rightarrow[0,+\infty)$ be non-decreasing. Then the following are equivalent:
\begin{itemize}
\item[$(i)$] $\omega$ is equivalent to a subadditive function $\sigma$.

\item[$(ii)$] $\omega$ satisfies \hyperlink{alpha0}{$(\alpha_0)$}.

\item[$(iii)$] $\omega$ satisfies
\begin{equation}\label{alpha0rewritten}
\exists\;C\ge 1\;\exists\;D\ge 1\;\forall\;\lambda\ge 1\;\forall\;t\ge 0:\;\;\;\omega(\lambda t)\le C\lambda\omega(t)+D\lambda.
\end{equation}

\item[$(iv)$] $\omega$ is equivalent to its least concave majorant $F_{\omega}$.
\end{itemize}
\end{lemma}

\emph{Note:}

\begin{itemize}
\item[$(*)$] Assumption $\omega(0)\ge 0$ in \cite[Lemma 4.1]{subaddlike} is superfluous.

\item[$(*)$] To be formally correct, in the proof of \cite[Lemma 4.1 $(ii)\Rightarrow(iii)$]{subaddlike} one shall assume $s>0$ and note that the second part of \cite[$(12)$]{index} is trivial when $s=0$.

\item[$(*)$] In \cite{subaddlike} equivalence is defined via O-growth relations like in \eqref{bigOrelation}, see \cite[Sect. 2.4, p. 406]{subaddlike}, and which coincides with the above given notion \hyperlink{sim}{$\sim$} since \eqref{liminfweighfct} holds by $\lim_{t\rightarrow+\infty}\omega(t)=+\infty$ as a basic requirement in this work.

    However, in \cite[Lemma 4.1]{subaddlike} it was only assumed that $\omega$ is non-decreasing but this is sufficient to get \eqref{liminfweighfct} when $\omega\neq 0$ and as seen before both notions of equivalence coincide in this case.

    And if $\omega=0$, then one immediately infers $F_{\omega}= 0$ too and the above result is trivial.
\end{itemize}

\hyperlink{alpha0}{$(\alpha_0)$} should be compared with $(\alpha_1)$ in \cite{PetzscheVogt} which reads as follows:
\begin{equation}\label{alpha1}
\sup_{\lambda\ge 1}\limsup_{t\rightarrow+\infty}\frac{\omega(\lambda t)}{\lambda\omega(t)}<+\infty.
\end{equation}
Note that in \cite[Prop. 1.1]{PetzscheVogt} it has been claimed that \eqref{alpha1} is equivalent to \hyperlink{alpha0}{$(\alpha_0)$} but $(\alpha_1)$ is weaker and the difference is subtle: \hyperlink{alpha0}{$(\alpha_0)$} requires for the quotient under consideration \emph{uniform boundedness} for all $t\ge t_0$ and for all $\lambda\ge 1$ and $t_0$ must not depend on $\lambda$ whereas for $(\alpha_1)$ it is sufficient to have boundedness for each $\lambda\ge 1$ separately but here the values for $t$ may depend on $\lambda$, i.e. having the estimate only for all $t\ge t_{\lambda}$, cf. \cite[$(5.1)$]{subaddlike}. We refer to \cite{subaddlike} for more comments and explanations and also to the discussion in \cite[Sect. 3.8, p. 26-27]{dissertation}: There the (counter-)example given by $\omega(t):=t\log(t)$ for $t>1$ and $\omega(t):=0$ for $t\in[0,1]$ is studied.\vspace{6pt}

In the next step we show that \eqref{alpha0rewritten} is also related to the \emph{largest convex minorant;} see the proofs of \cite[Prop. 2.24 \& Cor. 2.26]{index}.

For arbitrary $\omega:[0,+\infty)\rightarrow[0,+\infty)$ introduce $\omega^{\iota}:(0+\infty)\rightarrow[0,+\infty)$ given by $\omega^{\iota}(t):=\omega(1/t)$. So $\omega^{\iota}$ is non-increasing if and only if $\omega$ is non-decreasing, $\lim_{t\rightarrow 0}\omega^{\iota}(t)=\lim_{t\rightarrow+\infty}\omega(t)$ and $\lim_{t\rightarrow+\infty}\omega^{\iota}(t)=\lim_{t\rightarrow 0}\omega(t)$ if the limits exist. Moreover, $\omega^{\iota}=0$ if and only if $\omega=0$.

$\overline{\omega^{\iota}}^c$ shall denote the \emph{largest convex minorant of $\omega^{\iota}$} which can be represented as follows, see e.g. the proof of \cite[Prop. 2.24]{index}:
$$\overline{\omega^{\iota}}^c(t):=\inf\{\lambda_1\omega^{\iota}(t_1)+\lambda_2\omega^{\iota}(t_2):\;\lambda_1+\lambda_2=1,\;\lambda_1t_1+\lambda_2t_2=t,\;\lambda_i\ge 0\}.$$
Of course, this definition makes sense for any function $h:(0+\infty)\rightarrow[0,+\infty)$ and $\overline{h}^c(t)\le h(t)$ for all $t$ when choosing $\lambda_1=1$, $\lambda_2=0$ and $t_1=t$, $t_2\in(0,+\infty)$ arbitrary.

Now we are ready to prove the main result; one of the referees has modified in the report the original proof and was able to replace the stronger assumption \eqref{alpha0rewritten} in the preprint by \hyperlink{om1}{$(\omega_1)$}.

\begin{proposition}\label{alpha0impl}
Let $\omega:[0,+\infty)\rightarrow[0,+\infty)$ be non-decreasing. If $\omega$ satisfies \hyperlink{om1}{$(\omega_1)$}, then
$$\exists\;A\ge 1\;\forall\;t\in(0,+\infty):\;\;\;\overline{\omega^{\iota}}^c(t)\le\omega^{\iota}(t)\le A\overline{\omega^{\iota}}^c(t)+A;$$
hence $\omega^{\iota}$ is equivalent to its largest convex minorant.
\end{proposition}

\emph{Note:} If $\omega=0$, then the conclusion is trivial.

\demo{Proof}
We verify \cite[$(14)$]{index} for $\omega^{\iota}$, i.e.
\begin{equation}\label{equ14}
\exists\;C\ge 1\;\exists\;\beta>0\;\forall\;s,t\in(0,+\infty):\;\;\;\omega^{\iota}(s)+C\ge\frac{1}{C}\min\left\{1,\frac{t^{\beta}}{s^{\beta}}\right\}\omega^{\iota}(t),
\end{equation}
and then apply \cite[Prop. 2.24]{index}. Note that \eqref{equ14} is equivalent to the existence of universal constants $C_1,C,D\ge 1$ such that $C_1\omega^{\iota}(s)+D\ge\frac{1}{C_2}\min\left\{1,\frac{t^{\beta}}{s^{\beta}}\right\}\omega^{\iota}(s)$ holds for all $s,t\in(0,+\infty)$ (set $C:=\max\{D/C_1,C_1C_2\}$).

First, since $\omega$ is non-decreasing for any $0<s\le t$ we have $\omega^{\iota}(t)=\omega(1/t)\le\omega(1/s)=\omega^{\iota}(s)$ and $\min\{1,\frac{t^{\beta}}{s^{\beta}}\}=1$ for all $\beta>0$. Thus \eqref{equ14} holds in this situation with any $C\ge 1$ and $\beta>0$.

For the second case we follow the ideas presented by one of the referees in the report: By \hyperlink{om1}{$(\omega_1)$}, and since $\omega$ is non-decreasing, we get $\omega(2t)\le L\omega(t)+L$ for some $L>1$ and all $t\ge 0$. By iteration we obtain (see e.g. \cite[Sect. 3.3, (3.3.1)]{dissertation})
$$\exists\;L>1\;\forall\;n\in\NN_{>0}\;\forall\;t\ge 0:\;\;\;\omega(2^nt)\le L^n\omega(t)+\sum_{j=1}^nL^j=L^n\omega(t)+\frac{L^{n+1}-L}{L-1}\le L^n\omega(t)+L^{n+1}.$$
According to this given $L$ we set $\beta:=\log(L)/\log(2)$, i.e. $2^{\beta}=L$, and so $\beta>0$. Now let $\lambda\ge 1$ be given, choose $n\in\NN_{>0}$ such that $2^n\le\lambda<2^{n+1}$ and estimate as follows for all $t\ge 0$:
$$\omega(\lambda t)\le\omega(2^{n+1}t)\le L^{n+1}\omega(t)+L^{n+2}\le\lambda^{\beta}(2^{\beta}\omega(t)+4^{\beta}).$$
The second estimate is valid since $L^{n+1}=2^{\beta(n+1)}\le(2\lambda)^{\beta}\Leftrightarrow 2^n\le\lambda$ and $L^{n+2}=2^{\beta(n+2)}\le(4\lambda)^{\beta}\Leftrightarrow 2^n\le\lambda$. Thus,
with $t':=\frac{1}{\lambda t}$ we get
$$\forall\;\lambda\ge 1\;\forall\;t'>0:\;\;\;\omega^{\iota}(t')=\omega(\lambda t)\le\lambda^{\beta}(2^{\beta}\omega(t)+4^{\beta})=\lambda^{\beta}(2^{\beta}\omega^{\iota}(\lambda t')+4^{\beta}).$$
Let now $0<t<s$ be arbitrary and with $\lambda=\frac{s}{t}>1$ the above estimate gives
$$\left(\frac{s}{t}\right)^{\beta}(2^{\beta}\omega^{\iota}(s)+4^{\beta})=\lambda^{\beta}(2^{\beta}\omega^{\iota}(\lambda t)+4^{\beta})\ge\omega^{\iota}(t),$$
and so it follows that
$$\forall\;s>t>0:\;\;\;2^{\beta}\omega^{\iota}(s)+4^{\beta}\ge\left(\frac{t}{s}\right)^{\beta}\omega^{\iota}(t)=\min\left\{1,\frac{t^{\beta}}{s^{\beta}}\right\}\omega^{\iota}(t).$$
Consequently, \eqref{equ14} is shown and \cite[Prop. 2.24]{index} yields the assertion.
\qed\enddemo

On the other hand, if there exists a convex function $\sigma:(0,+\infty)\rightarrow[0,+\infty)$ which is equivalent to $\omega^{\iota}$, then there exists some $A\ge 1$ such that for all $t,s\in(0,+\infty)$ and $\lambda_1,\lambda_2\in[0,1]$ satisfying $\lambda_1+\lambda_2=1$:
\begin{align*}
\omega^{\iota}(\lambda_1t+\lambda_2s)&\le A\sigma(\lambda_1t+\lambda_2s)+A\le A\lambda_1\sigma(t)+A\lambda_2\sigma(s)+A
\\&
\le A^2\lambda_1\omega^{\iota}(t)+A^2\lambda_1+A^2\lambda_2\omega^{\iota}(s)+A^2\lambda_2=A^2\lambda_1\omega^{\iota}(t)+A^2\lambda_2\omega^{\iota}(s)+A^2.
\end{align*}
And if there exists some $t_0>0$ such that $\omega(t)=0$ for all $t\in[0,t_0]$, then $\omega^{\iota}(t)=0$ for all $t\ge t_0^{-1}$. Applying the previous estimate to $s:=t_0^{-1}$ gives $\omega^{\iota}(\lambda_1t+\lambda_2t_0^{-1})\le A^2\lambda_1\omega^{\iota}(t)+A^2$ for all $t\in(0,t_0]$ and $\lambda_1,\lambda_2\ge 0$ satisfying $\lambda_1+\lambda_2=1$.

Remark that the proof and the conclusion from \cite[Lemma 3.8.1 $(2)$]{dissertation} fails since super-additivity for $\omega^{\iota}$ would imply $2\omega^{\iota}(t)\le\omega^{\iota}(2t)\le\omega^{\iota}(t)$ for all $t\in(0,+\infty)$, a contradiction.

\subsection{Comment on weights being defined on $\RR^d$ or $\CC^d$}
Usually, when being defined on $\RR^d$, $d\ge 2$, it is assumed that $\omega$ is extended \emph{radially} and so $\omega(x):=\omega(|x|)$ for any $x\in\RR^d$ \emph{(isotropic setting).} On the other hand, for $z\in\CC^d$, in the literature one can find the definition $\omega(z):=\omega(|z|_1)$ with $|z|_1:=\sum_{j=1}^d|z_j|$ for $z=(z_1,\dots,z_d)$; see e.g. \cite[Def. 1.1]{BraunMeiseTaylor90}.

Write $z_j=x_j+iy_j$, $1\le j\le d$, then $|z|_1=\sum_{j=1}^d(x_j^2+y_j^2)^{1/2}$ and $|z|=\left(\sum_{j=1}^d|z_j|^2\right)^{1/2}=\left(\sum_{j=1}^dx_j^2+y_j^2\right)^{1/2}$, so $|z|\le|z|_1\le d\max\{|z_j|: 1\le j\le d\}\le d|z|$ holds for all $z\in\CC^d$.

Using these estimates we see that for all non-decreasing $\omega$ satisfying \hyperlink{om1}{$(\omega_1)$} both notions are equivalent in the sense that the mappings $z\mapsto\omega(|z|)$ and $z\mapsto\omega(|z|_1)$ are equivalent w.r.t. relation \hyperlink{sim}{$\sim$}.\vspace{6pt}

More recently also non-radial weights have been considered \emph{(anisotropic setting);} we refer e.g. to \cite{NeytVindas23} and \cite{DebrouwereNeytVindas24}.

\subsection{The growth indices $\gamma(\omega)$ and $\overline{\gamma}(\omega)$}\label{growthindexsection}
Let in this section $\omega:[0,+\infty)\rightarrow[0,+\infty)$ be non-decreasing and such that $\lim_{t\rightarrow+\infty}\omega(t)=+\infty$; i.e. $\omega$ is a weight function as in \cite[Sect. 2.3]{index}. We briefly recall the definition of the growth index $\gamma(\omega)$, see \cite[Sect. 2.3]{index} and the references therein. For $\gamma>0$ we say that $\omega$ has property $(P_{\omega,\gamma})$ if
\begin{equation*}\label{newindex1}
\exists\;K>1:\;\;\;\limsup_{t\rightarrow+\infty}\frac{\omega(K^{\gamma}t)}{\omega(t)}<K.
\end{equation*}
If $(P_{\omega,\gamma})$ holds for some $K>1$, then also $(P_{\omega,\gamma'})$ is satisfied for all $\gamma'\le\gamma$ with the same $K$ since $\omega$ is non-decreasing. Moreover, we can restrict to $\gamma>0$, because for $\gamma\le 0$ condition $(P_{\omega,\gamma})$ is automatically satisfied (again since $\omega$ is non-decreasing and $K>1$). Then we put
\begin{equation*}\label{newindex2}
\gamma(\omega):=\sup\{\gamma>0: (P_{\omega,\gamma})\;\;\text{is satisfied}\},
\end{equation*}
and if none condition $(P_{\omega,\gamma})$ holds then set $\gamma(\omega):=0$.

\begin{remark}\label{secondcomprem}
Let $\omega$ be as before. We recall some connections between $\gamma(\omega)$ and sub-additivity-like conditions for $\omega$, see also \cite[Sect. 6]{subaddlike}:

\begin{itemize}
\item[$(*)$] By \cite[Cor. 2.14]{index} we get $\gamma(\omega)>0$ if and only if \hyperlink{om1}{$(\omega_1)$} is valid.

\item[$(*)$] By \cite[Thm. 2.11, Cor. 2.13]{index} we have that $\gamma(\omega)>1$ if and only if $\omega$ has \hyperlink{omsnq}{$(\omega_{\on{snq}})$}.

\item[$(*)$] If $\omega$ is also continuous, then by (the proof of) \cite[Prop. 1.3]{MeiseTaylor88} it is known that \hyperlink{omsnq}{$(\omega_{\on{snq}})$} is equivalent to the fact that $\omega\hyperlink{sim}{\sim}\kappa_{\omega}$, where
    \begin{equation}\label{fctkappa}
    \kappa_{\omega}(y):=\int_1^{+\infty}\frac{\omega(yt)}{t^2}dt.
    \end{equation}
    $\kappa_{\omega}$ is concave, continuous and we have $\kappa_{\omega}(0)=\omega(0)\ge 0$ (i.e. $\kappa_{\omega}(0)=0$ if and only if $\omega(0)=0$). Hence it is known that also \hyperlink{sub}{$(\omega_{\on{sub}})$} holds for $\kappa_{\omega}$; see e.g. \cite[Lemma 3.8.1 $(1)$]{dissertation}.

    \emph{Note:} If continuity for $\omega$ fails, then the proof of \cite[Thm. 2.11 $(i)\Rightarrow(ii)$]{index} yields the fact that we get all desired properties for $y\mapsto\widetilde{\kappa}_{\omega}(y):=\int_1^{+\infty}\frac{\kappa_{\omega}(yt)}{t^2}dt$
and so $\omega$ in \eqref{fctkappa} is replaced by $\kappa_{\omega}$. Consequently, in any case one has that $\omega$ is equivalent to a continuous weight satisfying \hyperlink{sub}{$(\omega_{\on{sub}})$}.
\end{itemize}
\end{remark}

The comments listed in Section \ref{relevantgrowthsect}, Remark \ref{secondcomprem} and Lemma \ref{alpha0equivalence} give the following result:

\begin{lemma}\label{chaimimplemma}
Let $\omega:[0,+\infty)\rightarrow[0,+\infty)$ be non-decreasing and such that $\lim_{t\rightarrow+\infty}\omega(t)=+\infty$. Then the following implications hold
\begin{equation}\label{chaimofimpl}
\gamma(\omega)=+\infty\Rightarrow\gamma(\omega)>1\Leftrightarrow\hyperlink{omsnq}{(\omega_{\on{snq}})}\Rightarrow\hyperlink{omnq}{(\omega_{\on{nq}})}\Rightarrow\hyperlink{om5}{(\omega_5)}\Rightarrow\hyperlink{om2}{(\omega_2)},
\end{equation}
and, moreover,
\begin{equation}\label{chaimofimpl1}
\gamma(\omega)>1\Rightarrow\hyperlink{alpha0}{(\alpha_0)}.
\end{equation}
\end{lemma}

Analogously, for $\gamma>0$ we say that $\omega$ has property $(\overline{P}_{\omega,\gamma})$ if
\begin{equation}\label{newindex3}
\exists\;A>1:\;\;\;\liminf_{t\rightarrow+\infty}\frac{\omega(A^{\gamma}t)}{\omega(t)}>A.
\end{equation}
If $(\overline{P}_{\omega,\gamma})$ holds for some $A>1$, then $(\overline{P}_{\omega,\gamma'})$ is satisfied for all $\gamma'\ge\gamma$ with the same $A$ since $\omega$ is non-decreasing. Moreover, we can restrict to $\gamma>0$ because for $\gamma\le 0$ condition $(\overline{P}_{\omega,\gamma})$ is never satisfied for any weight function: $\omega$ is assumed to be non-decreasing and $A>1$. Then set
\begin{equation}\label{newindex4}
\overline{\gamma}(\omega):=\inf\{\gamma>0: \;\;(\overline{P}_{\omega,\gamma})\;\;\text{is satisfied}\}.
\end{equation}
We obtain $(0\le)\gamma(\omega)\le\overline{\gamma}(\omega)$, see again \cite[Sect. 2.3 \& 2.4]{index}, and via \cite[Thm. 2.16 \& Cor. 2.17]{index} we get $\overline{\gamma}(\omega)<+\infty$ if and only if $\omega$ satisfies \hyperlink{om6}{$(\omega_6)$}.

\subsection{Comments on different notions of weight functions}\label{weightfctnotionsect}
In this section we comment in more detail on the different notions of weight functions appearing in the literature; recall the definitions given in the introduction (Section \ref{Intro}).\vspace{6pt}

For the abstract study of weight functions, the assumptions continuity and $\omega(0)=0$ in Definition \ref{defweightfct} can be removed; e.g. when investigating the growth indices $\gamma(\omega)$ and $\overline{\gamma}(\omega)$; see \cite[Sect. 2.2-2.5]{index} and Section \ref{growthindexsection}.

Moreover, for a BMT-weight function (see Definition \ref{defBMTweightfct}) we have the following terminology:

\begin{itemize}
\item[$(*)$] If $\omega$ satisfies in addition $\omega(t)=0$ for all $t\in[0,1]$, then $\omega$ is called \emph{normalized.}

\item[$(*)$] If $\omega$ has in addition \hyperlink{omnq}{$(\omega_{\on{nq}})$}, then $\omega$ is called \emph{non-quasianalytic} and \emph{quasianalytic} otherwise.

\item[$(*)$] If \hyperlink{omsnq}{$(\omega_{\on{snq}})$} holds then $\omega$ is called \emph{strong non-quasianalytic} or even \emph{strong} for short. This notion is inspired by \cite[Def. 1.8]{BonetBraunMeiseTaylorWhitneyextension}.
\end{itemize}

\begin{remark}\label{BMT-remark}
We gather some comments on the assumptions of BMT-weights:
\begin{itemize}
\item[$(*)$] In \cite[Def. 1.1]{BraunMeiseTaylor90} property \hyperlink{omnq}{$(\omega_{\on{nq}})$} and normalization have been standard assumptions.

\item[$(*)$] It is also known that normalization can be assumed w.l.o.g. since the corresponding weighted classes are preserved under equivalence; for the ultradifferentiable classes $\mathcal{E}_{\{\omega\}}$ resp. $\mathcal{E}_{(\omega)}$ see \cite[Rem. 1.2 $(b)$]{BonetBraunMeiseTaylorWhitneyextension} and \cite[Cor. 5.17]{compositionpaper}.

\item[$(*)$] Formally, the basic requirement $\lim_{t\rightarrow+\infty}\omega(t)=+\infty$ can be skipped since it follows in any case by \hyperlink{om3}{$(\omega_3)$}.

\item[$(*)$] Based on the notation \cite[Def. 1.1]{BraunMeiseTaylor90}, in the literature when dealing with BMT-weight functions \hyperlink{om1}{$(\omega_1)$} is frequently denoted by $(\alpha)$, \hyperlink{omnq}{$(\omega_{\on{nq}})$} by $(\beta)$, \hyperlink{om3}{$(\omega_3)$} by $(\gamma)$ and \hyperlink{om4}{$(\omega_4)$} by $(\delta)$.

\item[$(*)$] According to \cite[Sect. 2.1]{optimalRoumieu} and \cite[Sect. 2.2]{optimalBeurling} a (normalized) BMT-weight $\omega$ is called a (normalized) weight function and when only \hyperlink{om1}{$(\omega_1)$} fails then it is called a (normalized) \emph{pre-weight function.}
\end{itemize}
\end{remark}

Concerning BB-weights in \cite{Bjorck66}
\begin{itemize}
\item[$(*)$] the set of functions from Definition \ref{defBBweightfct} is denoted by $\mathfrak{M}$ and $\omega$ has been radially extended to $\RR^d$;

\item[$(*)$] requirements $(ii)$ and $(iii)$ together are abbreviated by $(\alpha)$, $(iv)$ by $(\beta)$ and $(v)$ by $(\gamma)$;

\item[$(*)$] according to \cite[Def. 1.2.5]{Bjorck66} the set $\mathfrak{M}_0$ denotes functions satisfying $(i)-(iv)$ in Definition \ref{defBBweightfct}.
\end{itemize}

For the sake of completeness let us mention that also in \cite[Sect. 8]{BraunMeiseTaylor90} the notion of BB-weights has been considered for a comparison of different theories. More precisely, in \cite[8.4 $(2)$]{BraunMeiseTaylor90} for weights in the sense of Beurling-Bj\"{o}rck all requirements from Definition \ref{defBBweightfct} except (formally) continuity have been assumed.\vspace{6pt}

Finally, concerning PV-weights in \cite{PetzscheVogt} (recall Definition \ref{defPVweightfct}) occasionally the following conditions appear:

\begin{itemize}
\item[$(5)$] \hyperlink{om3}{$(\omega_3)$} (abbreviated by $(\gamma_1)$),

\item[$(6)$] \hyperlink{om2}{$(\omega_2)$} (abbreviated by $(\beta_0)$), and

\item[$(7)$] $(\delta)$ which is
$$\exists\;H\ge 1:\;\;\;\liminf_{t\rightarrow+\infty}\frac{\omega(Ht)}{\omega(t)}>1.$$
\end{itemize}

Thus, $(\delta)$ in \cite[Sect. 5]{PetzscheVogt} does not mean \hyperlink{om4}{$(\omega_4)$} (i.e. the convexity of $\varphi_{\omega}$) and $(\gamma_1)$ must not be mixed and confused with Petzsche's condition for weight sequences introduced in \cite{petzsche}. Since $\omega$ is non-decreasing in $(\delta)$ we necessarily have to choose $H>1$ and, indeed, via \eqref{newindex3} and \eqref{newindex4} condition $(\delta)$ does precisely mean $\overline{\gamma}(\omega)<+\infty$; i.e. \hyperlink{om6}{$(\omega_6)$} for $\omega$. And this observation (for the associated weight function case $\omega=\omega_{\mathbf{M}}$) is consistent with \cite[Lemma 5.3]{PetzscheVogt} and the known characterization of  \hyperlink{om6}{$(\omega_6)$} in terms of $\mathbf{M}$; see e.g. \cite[Thm. 3.1]{modgrowthstrange} and the references there.

\subsection{Matrix admissible weight functions}\label{LFYconjgsection}
In \cite{dissertation} and \cite{compositionpaper} we have developed the idea that to each BMT-weight function $\omega$ we can associate the \emph{weight matrix}
\begin{equation}\label{assoweightmatrix}
\mathcal{M}_{\omega}:=\{\mathbf{W}^{(\ell)}: \ell>0\},\hspace{15pt}W^{(\ell)}_j:=\exp(\frac{1}{\ell}\varphi^{*}_{\omega}(\ell j)),\;j\in\NN.
\end{equation}
Recall that $\varphi^{*}_{\omega}$ denotes the Legendre-Fenchel-Young conjugate of $\varphi_{\omega}=\omega\circ\exp$; see \eqref{LFYconj} and \hyperlink{om4}{$(\omega_4)$}. However, from the results and ideas in \cite{dissertation} and \cite{compositionpaper} and also from the given basic definition in \cite{BraunMeiseTaylor90}, it turns out that basically for the definition of $\mathcal{M}_{\omega}$ it is convenient and sufficient to work with normalized, continuous and non-decreasing $\omega:[0,+\infty)\rightarrow[0,+\infty)$ satisfying \hyperlink{om3}{$(\omega_3)$} and \hyperlink{om4}{$(\omega_4)$}. (Frequently, we have denoted this set of functions by $\hypertarget{omset0}{\mathcal{W}_0}$ and occasionally more properties on $\omega$ are required; see Definition \ref{defBMTweightfct}.)

Since the definition of the weighted spaces and of $\mathcal{M}_{\omega}$ involves the conjugate $\varphi^{*}_{\omega}$ from \eqref{LFYconj} we comment now on properties of this conjugate for arbitrary $\omega:[0,+\infty)\rightarrow[0,+\infty)$ in detail; see also \cite[Rem. 1.3, Lemma 1.5]{BraunMeiseTaylor90}. It turns out that even a more general setting is sufficient to introduce $\mathcal{M}_{\omega}$. Indeed, it suffices to assume that $\omega$ is non-decreasing and $\log(t)=o(\omega(t))$ as $t\rightarrow+\infty$ (i.e. \hyperlink{om3}{$(\omega_3)$}) in order to ensure well-definedness of the Young conjugate and hence of each sequence $\mathbf{W}^{(\ell)}$; we refer also to the comments in \cite[Sect. 3.1]{dissertation}, \cite[Sect. 2.7, $(a)-(c)$]{genLegendreconjBMT}, \cite[Sect. 6 $(III)$]{modgrowthstrangeII}. For future reference we call any non-decreasing function $\omega$ satisfying \hyperlink{om3}{$(\omega_3)$} \emph{matrix admissible.}

\begin{itemize}
\item[$(a)$] If $\omega$ is normalized, then $\sup_{y\in\RR}\{xy-\varphi_{\omega}(y)\}=\sup_{y\ge 0}\{xy-\varphi_{\omega}(y)\}$ for all $x\ge 0$. Note that in concrete applications, e.g. in the BMT-setting, it is natural to consider $x\in[0,+\infty)$ in \eqref{LFYconj} but, however, naturally one should consider $\sup_{y\in\RR}$ in the definition since $\RR$ is the natural domain of definition for $\varphi_{\omega}$.

\item[$(b)$] By definition $\varphi^{*}_{\omega}$ is convex and $\varphi^{*}_{\omega}(0)=\sup_{y\in\RR}\{-\varphi_{\omega}(y)\}$ and so, if $\inf_{t\ge 0}\omega(t)=0$, then $\varphi^{*}_{\omega}(0)=0$ holds and if $\omega$ is non-decreasing, then $\varphi^{*}_{\omega}(0)=-\omega(0)$.

    Thus, if $\omega$ is normalized, then $\varphi^{*}_{\omega}(0)=0=\omega(0)$ and, moreover, $\varphi^{*}_{\omega}$ is non-decreasing on $[0,+\infty)$ (see $(a)$). But, however, in general this property is not satisfied.

\item[$(c)$] Naturally, one should restrict $\varphi^{*}_{\omega}$ to $[0,+\infty)$: Let $x<0$, then
$$\varphi^{*}_{\omega}(x)=\sup_{y\in\RR}\{xy-\varphi_{\omega}(y)\}\ge\sup_{y\in(-\infty,0]}\{xy-\varphi_{\omega}(y)\}\ge\sup_{y\in(-\infty,0]}\{xy\}-\sup_{y\in(-\infty,0]}\{\varphi_{\omega}(y)\}=+\infty.$$
For this note that $\omega(t)\in[0,+\infty)$ for all $t\in[0,1]$ and so $\varphi^{*}_{\omega}=+\infty$ on $(-\infty,0)$.

\item[$(d)$] $x\mapsto\frac{\varphi^{*}_{\omega}(x)}{x}$ is non-decreasing on $(0,+\infty)$ and $\lim_{x\rightarrow+\infty}\frac{\varphi^{*}_{\omega}(x)}{x}=+\infty$: The first property follows because $y-\frac{\varphi_{\omega}(y)}{x_1}\le y-\frac{\varphi_{\omega}(y)}{x_2}\Leftrightarrow\frac{1}{x_2}\le\frac{1}{x_1}$ for all $0<x_1\le x_2$ and $y\in\RR$. Concerning the second one assume that $\frac{\varphi^{*}_{\omega}(x)}{x}\le C$ for some $C>0$ and all $x>0$, then by definition $x(y-C)\le\varphi_{\omega}(y)$ for all $x\ge 0$ and $y\in\RR$. But this is impossible: Choose a fixed $y>C$ and then let $x\rightarrow+\infty$.

\item[$(e)$] \hyperlink{om3}{$(\omega_3)$} for $\omega$ precisely means $\lim_{y\rightarrow+\infty}\frac{\varphi_{\omega}(y)}{y}=+\infty$ and this property characterizes $\varphi^{*}_{\omega}(x)<+\infty$ for all $x\in[0,+\infty)$: Indeed, if \hyperlink{om3}{$(\omega_3)$} is valid then for any $C>0$ (large) we can find $d\in\RR$, w.l.o.g. $d<0$, such that $\varphi_{\omega}(y)\ge yC+d\Leftrightarrow-d\ge yC-\varphi_{\omega}(y)$ for all $y\ge 0$. Since $\varphi_{\omega}(y)\ge 0$ for any $y\in\RR$ and because $C>0$ the same estimate holds for any $y\in\RR$. Hence $\varphi^{*}_{\omega}(C)\le -d$ follows. Since $C>0$ is arbitrary $\varphi^{*}_{\omega}$ is well-defined on $[0,+\infty)$. Conversely, if $\varphi^{*}_{\omega}$ is well-defined on $[0,+\infty)$, then for all $x\ge 0$ there exists $C_x>0$ such that $xy-\varphi_{\omega}(y)\le C_x$ for all $y\in\RR$ and so $x-\frac{C_x}{y}\le\frac{\varphi_{\omega}(y)}{y}$ for all $y>0$. When $x\ge 0$ is fixed, then this estimate verifies $\liminf_{y\rightarrow+\infty}\frac{\varphi_{\omega}(y)}{y}\ge x$ and as $x\rightarrow+\infty$ we infer \hyperlink{om3}{$(\omega_3)$} for $\omega$.

\item[$(f)$] Now consider the double conjugate $\varphi^{**}_{\omega}$: This function is by definition also convex and, by $(c)$, one has $\varphi^{**}_{\omega}(x)=\sup_{y\in\RR}\{xy-\varphi^{*}_{\omega}(y)\}=\sup_{y\ge 0}\{xy-\varphi^{*}_{\omega}(y)\}$ for all $x\in\RR$. Comments $(d)$ and $(e)$ together imply that $\varphi^{**}_{\omega}$ is always well-defined (on $\RR$): Indeed, $(d)$ gives \hyperlink{om3}{$(\omega_3)$} for $\omega^{*}:=\varphi^{*}_{\omega}\circ\log$ and so well-definedness on $[0,+\infty)$ whereas on $(-\infty,0)$ this follows since $\varphi^{**}_{\omega}$ is non-decreasing.

    This fact does also follow from $\varphi^{**}_{\omega}(x)=\sup_{y\in\RR}\{xy-\varphi^{*}_{\omega}(y)\}=\sup_{y\in\RR}\{xy-\sup_{z\in\RR}\{zy-\varphi_{\omega}(z)\}\}\le\varphi_{\omega}(x)$, $x\in\RR$ arbitrary.\vspace{6pt}

Concerning the converse, note that $\varphi^{**}_{\omega}=\varphi_{\omega}$ clearly implies \hyperlink{om4}{$(\omega_4)$} for $\omega$. If $\omega$ is in addition non-decreasing, then we show that \hyperlink{om4}{$(\omega_4)$} implies $\varphi^{**}_{\omega}\ge\varphi_{\omega}$ and so $\varphi^{**}_{\omega}=\varphi_{\omega}$. The following arguments are inspired by (the proof of) \cite[Prop. 1.6]{PetzscheVogt}:

Let $t\in\RR$ be fixed and since $\varphi_{\omega}$ is convex one has that $u\mapsto\frac{\varphi_{\omega}(u)-\varphi_{\omega}(t)}{u-t}$ is non-decreasing, $u\in\RR$ such that $u\neq t$. Since $\varphi_{\omega}$ is non-decreasing this difference quotient is non-negative. So there exists $s=s(t)\ge 0$ such that
$$\inf_{u>t}\frac{\varphi_{\omega}(u)-\varphi_{\omega}(t)}{u-t}=(\varphi_{\omega})'_{+}(t)\ge s\ge(\varphi_{\omega})'_{-}(t)=\sup_{u<t}\frac{\varphi_{\omega}(u)-\varphi_{\omega}(t)}{u-t}.$$
Multiplying these estimates with $u-t$ we get that $s(u-t)\le\varphi_{\omega}(u)-\varphi_{\omega}(t)$ for all $u\in\RR$ and note that for $u=t$ this estimate is becoming trivial. Consequently, $\inf_{s\in\RR}\sup_{u\in\RR}\{s(u-t)-\varphi_{\omega}(u)\}\le-\varphi_{\omega}(t)$ holds which is equivalent to having $-\inf_{s\in\RR}\sup_{u\in\RR}\{s(u-t)-\varphi_{\omega}(u)\}\ge\varphi_{\omega}(t)$ for all $t\in\RR$. The left-hand side is equal to $\sup_{s\in\RR}\{-\sup_{u\in\RR}\{s(u-t)-\varphi_{\omega}(u)\}\}=\sup_{s\in\RR}\{st-\sup_{u\in\RR}\{su-\varphi_{\omega}(u)\}\}=\varphi_{\omega}^{**}(t)$ and we are done.

\emph{Note:} If $\omega$ is normalized, non-decreasing and satisfies \hyperlink{om4}{$(\omega_4)$}, then by the above and $(b)$ we get $(0=)-\varphi^{*}_{\omega}(0)=\omega(0)=\omega(1)=\varphi_{\omega}(0)=\varphi^{**}_{\omega}(0)=\sup_{y\in\RR}\{-\varphi^{*}_{\omega}(y)\}=\sup_{y\ge 0}\{-\varphi^{*}_{\omega}(y)\}$ and for having $(0=)-\varphi^{*}_{\omega}(0)=\sup_{y\ge 0}\{-\varphi^{*}_{\omega}(y)\}$ it suffices to assume normalization since then $\varphi^{*}_{\omega}$ is non-decreasing on $[0,+\infty)$; see again $(a)$ and $(b)$.
\end{itemize}

Summarizing, assume that $\omega$ is \emph{non-decreasing} and such that \emph{\hyperlink{om3}{$(\omega_3)$}} holds. Then $\lim_{t\rightarrow+\infty}\omega(t)=+\infty$, so $\omega$ is a weight function in the notion of \cite{index} (and \cite{genLegendreconj}, \cite{genLegendreconjBMT}), and by $(a)-(e)$ it follows that each $\mathbf{W}^{(\ell)}$ is a well-defined sequence of non-negative real numbers, each $\mathbf{W}^{(\ell)}$ is \emph{log-convex} (i.e. $(M.1)$ in \cite{Komatsu73} is valid), $\lim_{j\rightarrow+\infty}(W^{(\ell)}_j)^{1/j}=+\infty$ for all $\ell>0$, and finally $W^{(\ell_1)}_j\le W^{(\ell_2)}_j$ for all $j\in\NN$, $0<\ell_1\le\ell_2$.  Thus basically ``non-decreasing'' and \hyperlink{om3}{$(\omega_3)$} for $\omega$ are sufficient to introduce $\mathcal{M}_{\omega}$ according to \eqref{assoweightmatrix} and this explains the notion \emph{``matrix admissible''.}

If $\omega$ is also normalized, then each $\mathbf{W}^{(\ell)}$ is normalized; i.e. $1=W^{(\ell)}_0\le W^{(\ell)}_1$. Note that, when $\omega$ is not normalized, then one can always switch to an equivalent weight $\sigma$ which is normalized and this notion of equivalence transfers in a precise manner to the corresponding associated weight matrices and which ensures the equality of the corresponding (matrix) weighted spaces; we refer to \cite[Sect. 6 $(III)$, $(6.6)$]{modgrowthstrangeII}.

Finally, \hyperlink{om4}{$(\omega_4)$} gives the equality $\varphi^{**}_{\omega}=\varphi_{\omega}$ and ensures the equivalence between $\omega$ and (each) $\omega_{\mathbf{W}^{(\ell)}}$, see \cite[Lemma 5.7]{compositionpaper}. When $\omega$ is matrix admissible but not having \hyperlink{om4}{$(\omega_4)$} then this equivalence is not true in general: This follows by the example constructed in Section \ref{failuresection} and recall that each associated weight function satisfies \hyperlink{om4}{$(\omega_4)$}, see also Remark \ref{nonassorem}.\vspace{6pt}

Moreover, if $\omega$ is matrix admissible and satisfies in addition \hyperlink{om1}{$(\omega_1)$}, then also the equalities $\mathcal{E}_{\{\mathcal{M}_{\omega}\}}=\mathcal{E}_{\{\omega\}}$ and $\mathcal{E}_{(\mathcal{M}_{\omega})}=\mathcal{E}_{(\omega)}$ (as l.c.v.s.) hold. Here the classes are introduced directly (by definition) via weighting the derivatives in terms of $\varphi^{*}_{\omega}$ (cf. \cite[Sect. 4 \& 5]{compositionpaper}). For this fact estimate \cite[$(5.10)$]{compositionpaper} is required and for its proof \hyperlink{om4}{$(\omega_4)$} is not needed. The same comment applies to analogously defined weighted spaces.

If in addition \hyperlink{om6}{$(\omega_6)$} is valid, then by \cite[Lemma 5.9, $(5.11)$]{compositionpaper} all $\mathbf{W}^{(\ell)}$ are equivalent; i.e. $0<\inf_{j\in\NN_{>0}}(W^{(\ell_1)}_j/W^{(\ell_2)}_j)^{1/j}\le\sup_{j\in\NN_{>0}}(W^{(\ell_1)}_j/W^{(\ell_2)}_j)^{1/j}<+\infty$ for all $\ell_1,\ell_2>0$. However, for the converse implication which is also shown in \cite[Lemma 5.9]{compositionpaper} the double-conjugate and hence \hyperlink{om4}{$(\omega_4)$} is required. And also for the characterization in \cite[Cor 5.8 (2)]{compositionpaper}, which establishes a connection to condition \emph{``moderate growth''} (i.e. $(M.2)$ from \cite{Komatsu73}) for some/each $\mathbf{W}^{(\ell)}$, the convexity is used by involving the equivalence between $\omega$ and (each) $\omega_{\mathbf{W}^{(\ell)}}$ (see \cite[Lemma 5.7]{compositionpaper}).

\subsection{The (counter-)example by U. Franken}
In \cite[Prop. 3]{Franken94} the following statement has been shown:

\begin{theorem}\label{Frankenthm}
There exists a non-quasianalytic BMT-weight function $\omega$ such that each continuous increasing function $\sigma$ with $\sigma\ge\omega$ cannot have simultaneously \hyperlink{sub}{$(\omega_{\on{sub}})$} and \hyperlink{omnq}{$(\omega_{\on{nq}})$} and hence cannot be a weight in the sense of Beurling-Bj\"{o}rck.
\end{theorem}

Recall for arbitrary functions $\sigma,\omega:[0,+\infty)\rightarrow[0,+\infty)$ the notation $\sigma\ge\omega$ which means $\sigma(t)\ge\omega(t)$ for all $t\ge 0$ (see Section \ref{growthrelsection}).

\section{Weighted $L^p$-type spaces}\label{generalclasssect}

\subsection{Basic definitions}\label{weightedbasicdefsection}
For functions $\sigma,\tau:\RR^d\rightarrow[0,+\infty)$ we write $\sigma\le\tau$ if $\sigma(x)\le\tau(x)$ for all $x\in\RR^d$ and similarly all further relations from Section \ref{growthrelsection} transfer. We give the following general definition:

\begin{definition}\label{functionmatrixdef}
Consider the set of functions
$$\mathcal{W}:=\{\omega^{\ell}: \RR^d\rightarrow[0,+\infty),\;\;\;\ell>0:\;\;\;\omega^{\ell_1}\le\omega^{\ell_2},\;\forall\;0<\ell_1\le\ell_2\}.$$
$\mathcal{W}$ is denoted as a \emph{weight function matrix} and $\ell>0$ is the crucial \emph{matrix index parameter.} $\mathcal{W}$ is called
\begin{itemize}
\item[$(*)$] \emph{simple} if $\mathcal{W}=\{\omega\}$ which means that $\omega^{\ell}\equiv\omega$ for any $\ell>0$;

\item[$(*)$] \emph{radial} or \emph{isotropic} if $\omega^{\ell}(x)=\omega^{\ell}(|x|)$ for $x\in\RR^d$ and $\ell>0$;

\item[$(*)$] \emph{non-decreasing} if all $\omega^{\ell}$ are non-decreasing which means $\omega^{\ell}(x)\le\omega^{\ell}(y)$ for all $x,y\in\RR^d$ such that $|x|\le|y|$;

\item[$(*)$] \emph{continuous} if all $\omega^{\ell}$ are continuous,

\item[$(*)$] \emph{unbounded} if
\begin{equation}\label{unboundedweight}
\forall\;\ell>0:\;\;\;\sup_{x\in\RR^d}\omega^{\ell}(x)=+\infty.
\end{equation}
\end{itemize}
\end{definition}
If $\mathcal{W}$ is non-decreasing, then $\sup_{x\in\RR^d: |x|=t}\omega^{\ell}(x)\le\inf_{\epsilon>0}\inf_{x\in\RR^d: |x|=t+\epsilon}\omega^{\ell}(x)$ for any $t\ge 0$ fixed (and any $\ell>0$). Therefore, if $\mathcal{W}$ is in addition continuous, then $\mathcal{W}$ has to be radial; i.e. each $\omega^{\ell}$ is constant on each sphere $|x|=t$.

If $\mathcal{W}$ is radial, non-decreasing and unbounded, then $\lim_{|x|\rightarrow+\infty}\omega^{\ell}(x)=+\infty$ for all $\ell>0$.\vspace{6pt}

For any $1\le p\le\infty$ set $L^p=L^p(\RR^d,\CC)$. Let $\mathcal{W}$ be a weight function matrix and consider the weighted $L^{\infty}$-type space by
$$\mathcal{L}^{\infty}_{\omega^{\ell}}:=\left\{f\in L^{\infty}:\;\;\;\|f\|_{\infty,\omega^{\ell}}:=\ess\sup_{x\in\RR^d}|f(x)|e^{\omega^{\ell} (x)}<+\infty\right\},$$
and define
\begin{equation*}\label{weightedef1}
\mathcal{L}^{\infty}_{\{\mathcal{W}\}}:=\bigcup_{\ell>0}\mathcal{L}^{\infty}_{\omega^{\ell}},\hspace{15pt}\mathcal{L}^{\infty}_{(\mathcal{W})}:=\bigcap_{\ell>0}\mathcal{L}^{\infty}_{\omega^{\ell}},
\end{equation*}
endowed with their natural topologies. Here, $\ess\sup$ denotes the \emph{essential supremum.} Moreover, for any $1\le p<\infty$ we assume that $\mathcal{W}$ is \emph{continuous} and for all (parameters) $\ell>0$ let us set
$$\mathcal{L}^p_{\omega^{\ell}}:=\left\{f\in L^p:\;\;\;\|f\|_{p,\omega^{\ell}}:=\left(\int_{\RR^d}|f(x)|^pe^{\omega^{\ell}(x)}dx\right)^{1/p}<+\infty\right\},$$
and put
\begin{equation*}\label{weightedef2}
\mathcal{L}^p_{\{\mathcal{W}\}}:=\bigcup_{\ell>0}\mathcal{L}^p_{\omega^{\ell}},\hspace{15pt}\mathcal{L}^p_{(\mathcal{W})}:=\bigcap_{\ell>0}\mathcal{L}^p_{\omega^{\ell}}.
\end{equation*}
Again, the spaces $\mathcal{L}^p_{\{\mathcal{W}\}}$ and $\mathcal{L}^p_{(\mathcal{W})}$ are endowed with their natural topologies. Here and throughout the whole paper to lighten notation we use the convention to write $f$ instead of $[f]$; i.e. identify a representative of the equivalence class with the whole class.

Moreover, from now on write $[\cdot]$ as a joint notation if we mean either $\{\cdot\}$ or $(\cdot)$. A similar convention is also applied to the appearing growth conditions for weight function matrices in the next sections. $\mathcal{L}^p_{\{\mathcal{W}\}}$ are called classes of \emph{Roumieu-type} and $\mathcal{L}^p_{(\mathcal{W})}$ of \emph{Beurling-type.} Note that for the Beurling-type all large parameters $\ell>0$ are becoming crucial, whereas for the Roumieu-type all small parameters $\ell$ are important.

The following obvious and continuous inclusions are valid for any (continuous) weight function matrix:
$$\mathcal{L}^p_{(\mathcal{W})}\subseteq\mathcal{L}^p_{\{\mathcal{W}\}}\subseteq L^p,\;\;\;1\le p\le\infty.$$
Note that for any simple weight function matrix the Roumieu- and the Beurling-type coincide. On the other hand, if $\mathcal{W}=\{\omega^{\ell}>0\}$ \emph{stabilizes above,} i.e. there exists $\ell_0>0$ such that $\omega^{\ell}\equiv\omega^{\ell_0}$ for all $\ell\ge\ell_0$, then $\mathcal{L}^p_{(\mathcal{W})}=\mathcal{L}^p_{(\omega^{\ell_0})}$ when identifying $\omega^{\ell_0}$ with the simple matrix $\mathcal{W}=\{\omega^{\ell_0}\}$. Similarly, if $\mathcal{W}$ \emph{stabilizes below,} i.e. there exists $\ell_0>0$ such that $\omega^{\ell}\equiv\omega^{\ell_0}$ for all $\ell\le\ell_0$, then $\mathcal{L}^p_{\{\mathcal{W}\}}=\mathcal{L}^p_{\{\omega^{\ell_0}\}}$. Therefore, in these cases one can omit all superfluous weight functions $\omega^{\ell}$, $\ell\neq\ell_0$, but the Roumieu- and the Beurling-type will not coincide in general.

\subsection{On the unboundedness of weight function matrices}\label{unboundedsect}
Assume that $\omega^{\ell_0}=0$ for some $\ell_0>0$, then $L^p$ coincides as set and l.c.v.s. with $\mathcal{L}^{p}_{\{\mathcal{W}\}}$ and, similarly, when $\omega^{\ell}=0$ for all $\ell>0$ then $L^p$ coincides with $\mathcal{L}^{p}_{(\mathcal{W})}$ (and with $\mathcal{L}^{p}_{\{\mathcal{W}\}}$).\vspace{6pt}

More generally, if
\begin{equation}\label{boundedweightequroum}
\exists\;\ell_0>0:\;\;\;\sup_{x\in\RR^d}\omega^{\ell_0}(x)<+\infty,
\end{equation}
then $\mathcal{L}^p_{\{\mathcal{W}\}}=L^p$ (as sets and l.c.v.s.) and, analogously, if
\begin{equation}\label{boundedweightequbeur}
\forall\;\ell>0:\;\;\;\sup_{x\in\RR^d}\omega^{\ell}(x)<+\infty,
\end{equation}
then $\mathcal{L}^p_{(\mathcal{W})}=L^p$ (as sets and l.c.v.s.), $1\le p\le\infty$. Now let us prove the converse implication.

\begin{proposition}\label{thirdweightcomppropzero}
Let $\mathcal{W}=\{\omega^{\ell}: \ell>0\}$ be a weight function matrix.

\begin{itemize}
\item[$(i)$] The inclusion (as sets) $L^{\infty}\subseteq\mathcal{L}^{\infty}_{\{\mathcal{W}\}}$  implies \eqref{boundedweightequroum} whereas $L^{\infty}\subseteq\mathcal{L}^{\infty}_{(\mathcal{W})}$  implies \eqref{boundedweightequbeur} for $\mathcal{W}$.

\item[$(ii)$] If $\mathcal{W}$ is continuous, then the inclusion (as sets) $L^p\subseteq\mathcal{L}^{p}_{\{\mathcal{W}\}}$  for some/any $1\le p<\infty$ implies \eqref{boundedweightequroum} and $L^p\subseteq\mathcal{L}^{p}_{(\mathcal{W})}$ implies \eqref{boundedweightequbeur} for $\mathcal{W}$.
\end{itemize}
\end{proposition}

\demo{Proof}
$(i)$ By assumption the equality $L^{\infty}=\mathcal{L}^{\infty}_{[\mathcal{W}]}$ holds (as sets) and we show that this already implies \eqref{boundedweightequroum} resp. \eqref{boundedweightequbeur} for $\mathcal{W}$: Otherwise, take the constant function $f\equiv1$, and then, if \eqref{boundedweightequroum} resp. \eqref{boundedweightequbeur} is violated, we get $\sup_{x\in\RR^d}|f(x)|e^{\omega^{\ell}(x)}=+\infty$ for all $\ell>0$ resp. for all sufficiently large $\ell>0$. But this implies $f\notin\mathcal{L}^{\infty}_{[\mathcal{W}]}=L^{\infty}$, a contradiction.\vspace{6pt}

$(ii)$ By assumption $L^p=\mathcal{L}^{p}_{[\mathcal{W}]}$ holds (as sets) and again we proceed by contradiction.

Assume now that \eqref{boundedweightequroum} is violated for $\mathcal{W}$. Then we can find a sequence $(x_n)_{n\in\NN_{>0}}$ in $\RR^d$ such that $\omega^{1/n}(x_n)>n^2$ for all $n\in\NN_{>0}$. By the continuity of each weight function and the pointwise order of all weight functions in the matrix one can assume that $|x_n|\rightarrow+\infty$ and, moreover, the continuity of (each) $\omega^{1/n}$ also implies that
\begin{equation}\label{thirdweightcomppropzeroequ}
\forall\;n\in\NN_{>0}\;\exists\;\delta_n>0\;\forall\;x\in I_n:=\{y\in\RR^d: |y-x_n|\le\delta_n\}:\;\;\;\omega^{1/n}(x)\ge n^2.
\end{equation}
W.l.o.g. we can assume $\delta_n\le 1$ and that the sets (balls) $I_n$ are pairwise disjoint. Then set $J_n:=\int_{x\in\RR^d: x\in I_n}1dx$ and let $g^p:\RR^d\rightarrow\RR$ be defined as follows:
\begin{equation}\label{thirdweightcomppropzeroequ1}
g^p(x):=\left(\frac{1}{2^nJ_n}\right)^{1/p},\;\;\;\text{for}\;x\in I_n,\hspace{15pt}g^p(x):=0,\;\;\;\text{for}\;x\notin I_n.
\end{equation}
First, $g^p\in L^p$ follows because $g^p$ is continuous on $\RR^d$ almost everywhere and
$$\int_{\RR^d}|g^p(x)|^pdx=\sum_{n=1}^{+\infty}\int_{x\in\RR^d: x\in I_n}|g^p(x)|^pdx=\sum_{n=1}^{+\infty}\frac{1}{2^nJ_n}\int_{x\in\RR^d: x\in I_n}1dx=\sum_{n=1}^{+\infty}\frac{1}{2^n}=1.$$
On the other hand let $\ell>0$ be given (small), arbitrary but fixed. There exists $n_{\ell}\in\NN_{>0}$ such that $\frac{1}{n}\le\ell$ for all $n\ge n_{\ell}$ and so, by \eqref{thirdweightcomppropzeroequ} and \eqref{thirdweightcomppropzeroequ1} and the pointwise order of the weight functions, for all $n\in\NN_{>0}$ with $n\ge n_{\ell}$ and $x\in\RR^d$ with $x\in I_n$ we get $|g^p(x)|^pe^{\omega^{\ell}(x)}\ge |g^p(x)|^pe^{\omega^{1/n}(x)}\ge\frac{1}{2^nJ_n}e^{n^2}$. Therefore, we can estimate as follows:
\begin{align*}
&\int_{\RR^d}|g^p(x)|^pe^{\omega^{\ell}(x)}dx=\sum_{n=1}^{+\infty}\int_{x\in\RR^d: x\in I_n}|g^p(x)|^pe^{\omega^{\ell}(x)}dx
\\&
=\sum_{n=1}^{n_{\ell}-1}\int_{x\in\RR^d: x\in I_n}|g^p(x)|^pe^{\omega^{\ell}(x)}dx+\sum_{n=n_{\ell}}^{+\infty}\int_{x\in\RR^d: x\in I_n}|g(x)|^pe^{\omega^{\ell}(x)}dx
\\&
\ge\sum_{n=1}^{n_{\ell}-1}\int_{x\in\RR^d: x\in I_n}|g^p(x)|^pe^{\omega^{\ell}(x)}dx+\sum_{n=n_{\ell}}^{+\infty}\frac{1}{2^nJ_n}e^{n^2}\int_{x\in\RR^d: x\in I_n}1dx
\\&
=\sum_{n=1}^{n_{\ell}-1}\int_{x\in\RR^d: x\in I_n}|g^p(x)|^pe^{\omega^{\ell}(x)}dx+\sum_{n=n_{\ell}}^{+\infty}\frac{e^{n^2}}{2^n}=+\infty.
\end{align*}
Since $\ell>0$ was arbitrary, this computation verifies $g^p\notin\mathcal{L}^{p}_{\{\mathcal{W}\}}=L^p$, a contradiction.\vspace{6pt}

If \eqref{boundedweightequbeur} is violated, then there exists some $\ell_0>0$ such that $\sup_{x\in\RR^d}\omega^{\ell}(x)=+\infty$ for all $\ell\ge\ell_0$ and so we can find a sequence $(x_n)_{n\in\NN_{>0}}$ in $\RR^d$ with $|x_n|\rightarrow+\infty$ and such that $\omega^{\ell_0+1/n}(x_n)>n^2$ for all $n\in\NN_{>0}$. Moreover, analogously to \eqref{thirdweightcomppropzeroequ}, we get
\begin{equation}\label{thirdweightcomppropzeroequ2}
\forall\;n\in\NN_{>0}\;\exists\;\delta_n>0\;\forall\;x\in\RR^d,\;x\in I_n:=\{y\in\RR^d: |y-x_n|\le\delta_n\}:\;\;\;\omega^{\ell_0+1/n}(x)\ge n^2.
\end{equation}
Again w.l.o.g. assume $\delta_n\le 1$ and that the sets $I_n$ are pairwise disjoint. Let $g^p:\RR^d\rightarrow\RR$ be defined as in \eqref{thirdweightcomppropzeroequ1} and so $g^p\in L^p$ follows as before. Now let $\ell>\ell_0$ be given (large), arbitrary but fixed. There exists $n_{\ell}\in\NN_{>0}$ such that $\ell_0+\frac{1}{n}\le\ell$ for all $n\ge n_{\ell}$ and so, by \eqref{thirdweightcomppropzeroequ2}, for all $n\in\NN_{>0}$ with $n\ge n_{\ell}$ and $x\in\RR^d$ with $x\in I_n$ we get $|g^p(x)|^pe^{\omega^{\ell}(x)}\ge |g^p(x)|^pe^{\omega^{\ell_0+1/n}(x)}\ge\frac{1}{2^nJ_n}e^{n^2}$. Using this and the same estimate of the crucial integral as before one infers
$$\int_{\RR^d}|g^p(x)|^pe^{\omega^{\ell}(x)}dx=+\infty,$$
and so, since $\ell>\ell_0$ was arbitrary, $g^p\notin\mathcal{L}^{p}_{(\mathcal{W})}$ which yields a contradiction.
\qed\enddemo

Proposition \ref{thirdweightcomppropzero} and the previous comments immediately imply the following:

\begin{theorem}\label{unboundednesscor}
Let $\mathcal{W}=\{\omega^{\ell}: \ell>0\}$ be a weight function matrix.

\begin{itemize}
\item[$(i)$] The following are equivalent:
\begin{itemize}
\item[$(*)$] Condition \eqref{boundedweightequroum} holds,

\item[$(*)$] $L^{\infty}=\mathcal{L}^{\infty}_{\{\mathcal{W}\}}$ holds as sets,

\item[$(*)$] $L^{\infty}=\mathcal{L}^{\infty}_{\{\mathcal{W}\}}$ holds as l.c.v.s.
\end{itemize}

\item[$(ii)$] The following are equivalent:
\begin{itemize}
\item[$(*)$] Condition \eqref{boundedweightequbeur} holds,

\item[$(*)$] $L^{\infty}=\mathcal{L}^{\infty}_{(\mathcal{W})}$ holds as sets,

\item[$(*)$] $L^{\infty}=\mathcal{L}^{\infty}_{(\mathcal{W})}$ holds as l.c.v.s.
\end{itemize}

\item[$(iii)$] If $\mathcal{W}$ is in addition continuous, then the following are equivalent:
\begin{itemize}
\item[$(*)$] Condition \eqref{boundedweightequroum} holds,

\item[$(*)$] $L^p=\mathcal{L}^{p}_{\{\mathcal{W}\}}$ holds as sets for some/any $1\le p<\infty$,

\item[$(*)$] $L^p=\mathcal{L}^{p}_{\{\mathcal{W}\}}$ holds as l.c.v.s. for some/any $1\le p<\infty$.
\end{itemize}

\item[$(iv)$] If $\mathcal{W}$ is in addition continuous, then the following are equivalent:
\begin{itemize}
\item[$(*)$] Condition \eqref{boundedweightequbeur} holds,

\item[$(*)$] $L^p=\mathcal{L}^{p}_{(\mathcal{W})}$ holds as sets for some/any $1\le p<\infty$,

\item[$(*)$] $L^p=\mathcal{L}^{p}_{(\mathcal{W})}$ holds as l.c.v.s. for some/any $1\le p<\infty$.
\end{itemize}
\end{itemize}
\end{theorem}

Therefore, in order to treat proper (new) weighted subclasses of $L^p$ it is natural to assume always that $\mathcal{W}$ is \emph{unbounded;} i.e. that \eqref{unboundedweight} holds. Note that in the Beurling-type case, if \eqref{boundedweightequbeur} formally fails and \eqref{unboundedweight} is (only) valid for all $\ell\ge\ell_0$ (or $\ell>\ell_0$) large, then by definition of the spaces we can skip all $0<\ell<\ell_0$ (or $0<\ell\le\ell_0$) without changing the weighted class since we are interested in all large $\ell$. Hence \eqref{unboundedweight} is natural for both types and analogous comments explain the notions for $\mathcal{W}$ being radial, non-decreasing or continuous.

However, note that in the forthcoming proofs it is not required necessarily that the weight function matrices are unbounded. We finish this section by commenting on the (formal) opposite extreme.

\begin{remark}\label{infinityremark}
Let $\mathcal{W}:=\{\omega^{\ell}: \ell>0\}$ be given and assume that there exists $\ell_0>0$ such that $\omega^{\ell}\equiv+\infty$ for all $\ell\ge\ell_0$ resp. such that $\omega^{\ell}\equiv+\infty$ for all $\ell>0$ (this notion means $\omega^{\ell}(x)=+\infty$ for all $x\in\RR^d$). In this case, formally, it makes sense to set $\mathcal{L}^{p}_{(\mathcal{W})}=\{0\}$ resp. $\mathcal{L}^{p}_{\{\mathcal{W}\}}=\{0\}$, $1\le p\le\infty$, and this convention does not change when (only) assuming that for any $\ell\ge\ell_0$ (resp. any $\ell>0$) one can find some $x_{\ell}\in\RR^d$ such that $\omega^{\ell}(x)=+\infty$ for all $x\in\RR^d$ such that $|x|\ge|x_{\ell}|$.

Conversely, we use the convention to identify $\{0\}$ with the weighted $L^p$-type spaces when being formally defined via the previously described ``non-standard weight matrices''.
\end{remark}

\subsection{On the non-triviality of weighted $L^p$-type spaces}\label{nontrivialsect}
We verify that the defined weighted spaces are \emph{non-trivial} under mild and natural assumptions on the (radial) weight function matrix $\mathcal{W}$.

\begin{lemma}\label{nontriviallemma}
Let $\mathcal{W}$ be non-decreasing and radial.
\begin{itemize}
\item[$(i)$] It holds that $\{0\}\subsetneq\mathcal{L}^{\infty}_{(\mathcal{W})}\subseteq\mathcal{L}^{\infty}_{\{\mathcal{W}\}}$.

\item[$(ii)$] If $\mathcal{W}$ is in addition continuous, then $\{0\}\subsetneq\mathcal{L}^{p}_{(\mathcal{W})}\subseteq\mathcal{L}^{p}_{\{\mathcal{W}\}}$ for all $1\le p<\infty$.
\end{itemize}
\end{lemma}

\demo{Proof}
$(i)$ If \eqref{boundedweightequbeur} is valid, then $\mathcal{L}^{\infty}_{(\mathcal{W})}=L^{\infty}\neq\{0\}$ by $(ii)$ in Theorem \ref{unboundednesscor}.

If \eqref{boundedweightequbeur} is violated and since $\mathcal{W}$ is non-decreasing we can find $\ell_0>0$ such that $\lim_{|x|\rightarrow\infty}\omega^{\ell}(x)=+\infty$ for all $\ell\ge\ell_0$.

We introduce the following functions:
$$\sigma^{\infty}(x):=\omega^{n+1}(x),\;\;\;\forall\;x\in\RR^d,\;|x|\in[n,n+1),\;\;\;n\in\NN,$$
and
\begin{equation*}\label{psiinf}
\psi^{\infty}_{\mathcal{W}}(x):=e^{-\sigma^{\infty}(x)^2}.
\end{equation*}
So for any $\ell>0$ and all $x\in\RR^d$ such that $|x|\in[n,n+1)$ one has $|\psi^{\infty}_{\mathcal{W}}(x)|e^{\omega^{\ell}(x)}=e^{\omega^{\ell}(x)-(\omega^{n+1}(x))^2}$. Let $\ell\ge\ell_0$ be fixed and hence for all $|x|\in[n,n+1)$ with $n+1\ge\ell$ we have $|\psi^{\infty}_{\mathcal{W}}(x)|e^{\omega^{\ell}(x)}\le e^{\omega^{\ell}(x)(1-\omega^{\ell}(x))}$ and this last expression tends to $0$ as $|x|\rightarrow+\infty$. Thus,
$$\forall\;\ell\ge\ell_0\;\exists\;C_{\ell}\ge 1\;\forall\;x\in\RR^d:\;\;\;|\psi^{\infty}_{\mathcal{W}}(x)|e^{\omega^{\ell}(x)}\le C_{\ell},$$
and so $\psi^{\infty}_{\mathcal{W}}\in\mathcal{L}^{\infty}_{(\mathcal{W})}$.\vspace{6pt}

$(ii)$ If \eqref{boundedweightequbeur} holds, then $\mathcal{L}^{p}_{(\mathcal{W})}=L^p\neq\{0\}$ by $(iv)$ in Theorem \ref{unboundednesscor}.

If \eqref{boundedweightequbeur} fails, then since $\mathcal{W}$ is non-decreasing there exists $\ell_0>0$ such that $\lim_{|x|\rightarrow\infty}\omega^{\ell}(x)=+\infty$ for all $\ell\ge\ell_0$. In particular, there exists $x_0\in\RR^d$, $|x_0|\ge 1$, such that $\omega^{\ell}(x)>1$ for all $\ell\ge\ell_0$ and all $x\in\RR^d$ with $|x|\ge|x_0|$. W.l.o.g. assume that $|x_0|=n_0=\ell_0\in\NN_{>0}$. Now, for any $n\in\NN_{>0}$ with $n\ge n_0$ we can find real numbers $a_n\ge 1$ such that
\begin{equation}\label{anchoice}
(\omega^{n_0}(x))^{a_n}\ge 2(d+1)\log(1+|x|),\;\;\;\forall\;x\in\RR^d,\;|x|\in[n,n+1),\;n\ge n_0.
\end{equation}
Indeed, for this it suffices to take $a_n:=\max\{1,\frac{\log(2(d+1)\log(2+n))}{\log(\omega^{n_0}(n))}\}$ and \eqref{anchoice} holds as well for all $\omega^{n}$, $n\ge n_0$. Moreover, for all $n\ge n_0$ choose positive reals $b_n\ge 1$ such that
\begin{equation}\label{bnchoice}
(\omega^n(x))^{b_n}\ge 2\omega^n(x),\;\;\;\forall\;x\in\RR^d,\;|x|\in[n,n+1),\;n\ge n_0.
\end{equation}
Indeed, it suffices to take
$$b_n:=\frac{\log(\omega^n(n+1))+\log(2)}{\log(\omega^n(n))}(>1).$$
Finally, in order to complete the definition one can set $a_n:=1$ and $b_n:=1$ for $0\le n<n_0$.

Now introduce the auxiliary weight function
$$\sigma^p(x):=(\omega^{n+1}(x))^{a_nb_n},\;\;\;\forall\;x\in\RR^d,\;|x|\in[n,n+1),\;n\in\NN,$$
and set
\begin{equation}\label{psip}
\psi^p_{\mathcal{W}}(x):=e^{-p^{-1}\sigma^p(x)}.
\end{equation}
Then, for any $\ell\ge\ell_0$ (fixed) we estimate as follows:
\begin{align*}
&\int_{\RR^d}|\psi^p_{\mathcal{W}}(x)|^pe^{\omega^{\ell}(x)}dx=\sum_{n=0}^{+\infty}\int_{x\in\RR^d,\;|x|\in[n,n+1)}e^{\omega^{\ell}(x)-(\omega^{n+1}(x))^{a_nb_n}}dx
\\&
=\sum_{n<\ell}\int_{x\in\RR^d,\;|x|\in[n,n+1)}e^{\omega^{\ell}(x)-(\omega^{n+1}(x))^{a_nb_n}}dx+\sum_{n\ge \ell}\int_{x\in\RR^d,\;|x|\in[n,n+1)}e^{\omega^{\ell}(x)-(\omega^{n+1}(x))^{a_nb_n}}dx
\\&
\le\sum_{n<\ell}\int_{x\in\RR^d,\;|x|\in[n,n+1)}e^{\omega^{\ell}(x)-(\omega^{n+1}(x))^{a_nb_n}}dx+\sum_{n\ge\ell}\int_{x\in\RR^d,\;|x|\in[n,n+1)}e^{\omega^{n}(x)-(\omega^{n}(x))^{a_nb_n}}dx
\\&
\le\sum_{n<\ell}\int_{x\in\RR^d,\;|x|\in[n,n+1)}e^{\omega^{\ell}(x)-(\omega^{n+1}(x))^{a_nb_n}}dx+\sum_{n\ge\ell}\int_{x\in\RR^d,\;|x|\in[n,n+1)}\frac{1}{(1+|x|^2)^{d+1}}dx<+\infty.
\end{align*}
The last inequality is valid since $\frac{1}{(1+|x|)^{2(d+1)}}\le\frac{1}{(1+|x|^2)^{(d+1)}}$ for all $x\in\RR^d$. For the second last one note that $e^{\omega^{n}(x)-(\omega^{n}(x))^{a_nb_n}}\le\frac{1}{(1+|x|)^{2(d+1)}}$ holds on the intervals under consideration because for $|x|\in[n,n+1)$, $n\ge\ell_0=n_0$, via the choices of $a_n$ and $b_n$ in \eqref{anchoice} and \eqref{bnchoice}: $$(\omega^{n}(x))^{a_nb_n}-\omega^{n}(x)\ge(\omega^{n}(x))^{a_{n}b_{n}}-(\omega^{n}(x))^{a_{n}}\ge(\omega^{n}(x))^{a_n}\ge(\omega^{\ell_0}(x))^{a_n}\ge 2(d+1)\log(1+|x|).$$

Consequently, since $\ell\ge\ell_0$ was arbitrary it follows that $\psi^p_{\mathcal{W}}\in\mathcal{L}^{p}_{(\mathcal{W})}$.
\qed\enddemo

\subsection{Stability under translation}\label{stabtranslationsect}
For concrete applications it is important to know whether weighted spaces satisfy certain stability properties. Concerning the classes introduced above we investigate \emph{invariance w.r.t. translations.} This property is somehow natural and crucially required in the proofs of the main characterizing results in the next section.

\emph{Notation:} For given $f:\RR^d\rightarrow\CC$ and $x_0\in\RR^d$ set $f_{x_0}:=f(x-x_0)$ and consider the \emph{translation operator} $T_{x_0}: f\mapsto f_{x_0}$.

The main statement in this context is the following result:

\begin{proposition}\label{translationlemma}
Let $\mathcal{S}=\{\sigma^{\ell}: \ell>0\}$ and $\mathcal{T}=\{\tau^{\ell}: \ell>0\}$  be weight function matrices.
\begin{itemize}
\item[$(i)$] If the matrices are related by
\begin{equation}\label{mixedom1refomequ2}
\forall\;\ell>0\;\exists\;n>0\;\exists\;L\ge 1\;\forall\;x,y\in\RR^d:\;\;\;\tau^{\ell}(x+y)\le\sigma^n(x)+\sigma^n(y)+L,
\end{equation}
then
\begin{align*}
&\forall\;\ell>0\;\exists\;n>0\;\exists\;L\ge 1\;\forall\;x_0\in\RR^d\;\forall\;f\in\mathcal{L}^{\infty}_{\sigma^n}:\;\;\;\|T_{x_0}(f)\|_{\infty,\tau^{\ell}}\le e^Le^{\sigma^n(x_0)}\|f\|_{\infty,\sigma^n}.
\end{align*}
If both matrices are in addition continuous, then for any $1\le p<\infty$ we get
\begin{align*}
&\forall\;\ell>0\;\exists\;n>0\;\exists\;L\ge 1\;\forall\;x_0\in\RR^d\;\forall\;f\in\mathcal{L}^{p}_{\sigma^n}:\;\;\;\|T_{x_0}(f)\|^p_{p,\tau^{\ell}}\le e^Le^{\sigma^n(x_0)}\|f\|^p_{p,\sigma^n}.
\end{align*}

\item[$(ii)$] If the matrices are related by
\begin{equation}\label{mixedom1refomequ5}
\forall\;n>0\;\exists\;\ell>0\;\exists\;L\ge 1\;\forall\;x,y\in\RR^d:\;\;\;\tau^{\ell}(x+y)\le\sigma^n(x)+\sigma^n(y)+L,
\end{equation}
then
\begin{align*}
&\forall\;n>0\;\exists\;\ell>0\;\exists\;L\ge 1\;\forall\;x_0\in\RR^d\;\forall\;f\in\mathcal{L}^{\infty}_{\sigma^{n}}:\;\;\;\|T_{x_0}(f)\|_{\infty,\tau^{\ell}}\le e^Le^{\sigma^{n}(x_0)}\|f\|_{\infty,\sigma^{n}}.
\end{align*}
\end{itemize}
If both matrices are in addition continuous, then for any $1\le p<\infty$ we get
\begin{align*}
&\forall\;n>0\;\exists\;\ell>0\;\exists\;L\ge 1\;\forall\;x_0\in\RR^d\;\forall\;f\in\mathcal{L}^{p}_{\sigma^{n}}:\;\;\;\|T_{x_0}(f)\|^p_{p,\tau^{\ell}}\le e^Le^{\sigma^{n}(x_0)}\|f\|^p_{p,\sigma^{n}}.
\end{align*}
\end{proposition}

Proposition \ref{translationlemma} implies the continuity of (each) operator $T_{x_0}$ between weighted $L^p$-type classes and note:

\begin{itemize}
\item[$(*)$] The choice for the new index $n$ is only depending on given $\ell$ via \eqref{mixedom1refomequ2} resp. the index $\ell$ is depending on given $n$ via \eqref{mixedom1refomequ5} and hence not on the translation parameter $x_0$. Moreover, it is also not depending on the choice of the representative of the equivalence class.

\item[$(*)$] The statement becomes unclear if $\RR^d$ is replaced by some open and bounded $U\subseteq\RR^d$.

\item[$(*)$] The special case $\mathcal{S}=\mathcal{T}$, i.e. \eqref{mixedom1refomequ2} for $\mathcal{S}$, implies that the classes $\mathcal{L}^{p}_{(\mathcal{S})}$, $1\le p\le\infty$, are translation invariant and each $T_{x_0}$ is acting continuously on the corresponding weighted space. And \eqref{mixedom1refomequ5} for $\mathcal{S}$ implies that $\mathcal{L}^{p}_{\{\mathcal{S}\}}$, $1\le p\le\infty$, is translation invariant and each $T_{x_0}$ is acting continuously.
\end{itemize}

\demo{Proof}
We only prove $(i)$ in detail since $(ii)$ is analogous. When taking into account \eqref{mixedom1refomequ2} we get for all $x_0\in\RR^d$:
\begin{align*}
\|T_{x_0}(f)\|_{\infty,\tau^{\ell}}&=\|f_{x_0}\|_{\infty,\tau^{\ell}}=\ess\sup_{x\in\RR^d}|f_{x_0}(x)|e^{\tau^{\ell}(x)}=\ess\sup_{x\in\RR^d}|f(x-x_0)|e^{\tau^{\ell}(x)}
\\&
=\ess\sup_{y\in\RR^d}|f(y)|e^{\tau^{\ell}(y+x_0)}\le e^{L}e^{\sigma^n(x_0)}\ess\sup_{y\in\RR^d}|f(y)|e^{\sigma^n(y)}=e^{L}e^{\sigma^n(x_0)}\|f\|_{\infty,\sigma^n}.
\end{align*}
Similarly, for all $1\le p<\infty$ and $x_0\in\RR^d$:
\begin{align*}
\|T_{x_0}(f)\|^p_{p,\tau^{\ell}}&=\|f_{x_0}\|^p_{p,\tau^{\ell}}=\int_{\RR^d}|f_{x_0}(x)|^pe^{\tau^{\ell}(x)}dx=\int_{\RR^d}|f(x-x_0)|^pe^{\tau^{\ell}(x)}dx
\\&
=\int_{\RR^d}|f(y)|^pe^{\tau^{\ell}(y+x_0)}dy\le e^{L}e^{\sigma^n(x_0)}\int_{\RR^d}|f(y)|^pe^{\sigma^n(y)}dy=e^{L}e^{\sigma^n(x_0)}\|f\|^p_{p,\sigma^n}.
\end{align*}
\qed\enddemo

We continue with an auxiliary technical result for radial weight matrices which should be compared with $(f)$ in Section \ref{relevantgrowthsect}.

\begin{lemma}\label{mixedom1refom}
Let $\mathcal{S}=\{\sigma^{\ell}: \ell>0\}$ and $\mathcal{T}=\{\tau^{\ell}: \ell>0\}$ be non-decreasing and radial weight function matrices. Consider the following (mixed) assertions:
\begin{itemize}
\item[$(i)$] The matrices are related by
\begin{equation}\label{mixedom1refomequ1}
\forall\;\ell>0\;\exists\;n>0\;\exists\;L\ge 1\;\forall\;t\ge 0:\;\;\;\tau^{\ell}(2t)\le\sigma^n(t)+L.
\end{equation}
\item[$(ii)$] The matrices are related by \eqref{mixedom1refomequ2}.

\item[$(iii)$] The matrices are related by
\begin{equation}\label{mixedom1refomequ3}
\forall\;\ell>0\;\exists\;n>0\;\exists\;L\ge 1\;\forall\;t\ge 0:\;\;\;\tau^{\ell}(2t)\le 2\sigma^n(t)+L.
\end{equation}
\end{itemize}
Then the implications $(i)\Rightarrow(ii)\Rightarrow(iii)$ are valid. The proof shows that we can choose for given $\ell$ in each assertion the same parameter $n$ and the assumption ``non-decreasing'' is only used to show the first implication. The analogous statement holds when ``$\forall\;\ell>0\;\exists\;n>0$'' is replaced by ``$\forall\;n>0\;\exists\;\ell>0$'' and involving \eqref{mixedom1refomequ5}; for future reference we only mention explicitly
\begin{equation}\label{mixedom1refomequ4}
\forall\;n>0\;\exists\;\ell>0\;\exists\;L\ge 1\;\forall\;t\ge 0:\;\;\;\tau^{\ell}(2t)\le\sigma^n(t)+L.
\end{equation}
\end{lemma}

\demo{Proof}
Since $\tau^{\ell}$ and $\sigma^n$ are non-decreasing, by \eqref{mixedom1refomequ1} one has $\tau^{\ell}(s+t)\le\tau^{\ell}(2\max\{s,t\})\le\sigma^n(\max\{s,t\})+L\le\sigma^n(s)+\sigma^n(t)+L$ for all $s,t\ge 0$. On the other hand, in \eqref{mixedom1refomequ2} we set $s=t$ and immediately get \eqref{mixedom1refomequ3}. The second part follows analogously.
\qed\enddemo

This technical lemma motivates the next comment:

\begin{remark}\label{translationlemmarem}
In view of Lemma \ref{mixedom1refom} the conclusions in Proposition \ref{translationlemma} hold when assuming for non-decreasing, radial (and continuous) matrices $\mathcal{S}$, $\mathcal{T}$, properties \eqref{mixedom1refomequ1} resp. \eqref{mixedom1refomequ4} which are ``more natural'' and easier to check.

Instead of \eqref{mixedom1refomequ1} one could expect the weaker mixed \hyperlink{om1}{$(\omega_1)$}-like condition
\begin{equation}\label{mixedom1beur}
\forall\;\ell>0\;\exists\;n>0:\;\;\;\tau^{\ell}(2t)=O(\sigma^n(t))\;\;\;t\rightarrow+\infty,
\end{equation}
and instead of \eqref{mixedom1refomequ4} the condition
\begin{equation}\label{mixedom1roum}
\forall\;n>0\;\exists\;\ell>0:\;\;\;\tau^{\ell}(2t)=O(\sigma^{n}(t))\;\;\;t\rightarrow+\infty.
\end{equation}
When both matrices are non-decreasing, then \eqref{mixedom1beur} is equivalent to $\tau^{\ell}(s+t)\le L(\sigma^n(s)+\sigma^n(t))+L$ for some $L\ge 1$ and all $s,t\ge 0$ and similarly for \eqref{mixedom1roum}; the relation between the indices $\ell$ and $n$ is preserved in this equivalence. This follows by generalizing directly the arguments in $(f)$ in Section \ref{relevantgrowthsect} to the mixed setting. However, when involving \eqref{mixedom1beur} and \eqref{mixedom1roum} the proof of Proposition \ref{translationlemma} yields terms $e^{L\sigma^n(y)}$. By the appearance of the additional multiplicative constant $L$ in the exponent the conclusion in Proposition \ref{translationlemma} becomes unclear in general. Similar comments apply to the growth relations introduced in the next section; we refer to Section \ref{specialcasesect1} for more explanations.
\end{remark}

\section{Characterization of inclusion relations for weighted $L^p$-type spaces}\label{charactinclsection}
The aim is to prove the characterization of inclusion relations between weighted $L^p$-type spaces defined by weight function matrices when assuming mild standard growth requirements. In order to proceed, first let us introduce the crucial growth relations. Let $\mathcal{S}=\{\sigma^{\ell}: \ell>0\}$, $\mathcal{T}=\{\tau^{\ell}: \ell>0\}$ be given and consider
\begin{equation}\label{growthrel1}
	\forall\;\ell>0\;\exists\;n>0\;\exists\;C\ge 1\;\forall\;x\in\RR^d:\;\;\;\tau^{\ell}(x)\le\sigma^n(x)+C,
\end{equation}
and
\begin{equation}\label{growthrel2}
	\forall\;n>0\;\exists\;\ell>0\;\exists\;C\ge 1\;\forall\;x\in\RR^d:\;\;\;\tau^{\ell}(x)\le\sigma^n(x)+C.
\end{equation}
We write $\mathcal{S}(\preceq)\mathcal{T}$ if \eqref{growthrel1} is valid and $\mathcal{S}\{\preceq\}\mathcal{T}$ if \eqref{growthrel2} holds and using this notation the goal is to verify the following main characterizing statements:

\begin{theorem}\label{maincharthmbeur}
Let $\mathcal{S}=\{\sigma^{\ell}: \ell>0\}$ and $\mathcal{T}=\{\tau^{\ell}: \ell>0\}$ be weight function matrices and assume that both matrices satisfy
\begin{equation}\label{maincharthmbeurequ}
\forall\;\ell>0\;\exists\;n>0\;\exists\;L\ge 1\;\forall\;x,y\in\RR^d:\;\;\;\omega^{\ell}(x+y)\le\omega^n(x)+\omega^n(y)+L;
\end{equation}
i.e. \eqref{mixedom1refomequ2} for the same matrix. Then the following are equivalent:

\begin{itemize}
\item[$(i)$] The matrices are related by $\mathcal{S}(\preceq)\mathcal{T}$.

\item[$(ii)$] The inclusion $\mathcal{L}^{\infty}_{(\mathcal{S})}\subseteq\mathcal{L}^{\infty}_{(\mathcal{T})}$  holds as sets.

\item[$(iii)$] The continuous inclusion $\mathcal{L}^{\infty}_{(\mathcal{S})}\subseteq\mathcal{L}^{\infty}_{(\mathcal{T})}$ holds.
\end{itemize}

If both matrices are also continuous, then analogous equivalences are valid when the symbol/functor $\mathcal{L}^{\infty}$ is replaced by $\mathcal{L}^p$ for some/any $1\le p<\infty$.
\end{theorem}

\begin{theorem}\label{maincharthmroum}
Let $\mathcal{S}=\{\sigma^{\ell}: \ell>0\}$ and $\mathcal{T}=\{\tau^{\ell}: \ell>0\}$ be weight function matrices and assume that both matrices satisfy
\begin{equation}\label{maincharthmroumequ}
\forall\;n>0\;\exists\;\ell>0\;\exists\;L\ge 1\;\forall\;x,y\in\RR^d:\;\;\;\omega^{\ell}(x+y)\le\omega^n(x)+\omega^n(y)+L;
\end{equation}
i.e. \eqref{mixedom1refomequ5} for the same matrix. Then the following are equivalent:

\begin{itemize}
\item[$(i)$] The matrices are related by $\mathcal{S}\{\preceq\}\mathcal{T}$.

\item[$(ii)$] The inclusion $\mathcal{L}^{\infty}_{\{\mathcal{S}\}}\subseteq\mathcal{L}^{\infty}_{\{\mathcal{T}\}}$  holds as sets.

\item[$(iii)$] The continuous inclusion $\mathcal{L}^{\infty}_{\{\mathcal{S}\}}\subseteq\mathcal{L}^{\infty}_{\{\mathcal{T}\}}$ holds.
\end{itemize}

If both matrices are also continuous, then analogous equivalences are valid when the symbol/functor $\mathcal{L}^{\infty}$ is replaced by $\mathcal{L}^p$ for some/any $1\le p<\infty$.
\end{theorem}

\emph{Note:} \eqref{maincharthmbeurequ} resp. \eqref{maincharthmroumequ} yield the fact that all weighted spaces under consideration in the above results are \emph{translation invariant.}\vspace{6pt}

First, the next statement is immediate by the defining system of seminorms and the introduced growth relations before and it verifies $(i)\Rightarrow(ii),(iii)$ in Theorems \ref{maincharthmbeur} and \ref{maincharthmroum}.

\begin{proposition}\label{firstweightcompprop}
Let $\mathcal{S}=\{\sigma^{\ell}: \ell>0\}$ and $\mathcal{T}=\{\tau^{\ell}: \ell>0\}$ be weight function matrices.

\begin{itemize}
\item[$(i)$] If $\mathcal{S}(\preceq)\mathcal{T}$, then
$$\forall\;\ell>0\;\exists\;n>0\;\exists\;C\ge 1\;\forall\;f\in\mathcal{L}^{\infty}_{\sigma^n}:\;\;\;\|f\|_{\infty,\tau^{\ell}}\le e^C\|f\|_{\infty,\sigma^n},$$
and if $\mathcal{S}\{\preceq\}\mathcal{T}$, then
$$\forall\;n>0\;\exists\;\ell>0\;\exists\;C\ge 1\;\forall\;f\in\mathcal{L}^{\infty}_{\sigma^{n}}:\;\;\;\|f\|_{\infty,\tau^{\ell}}\le e^C\|f\|_{\infty,\sigma^{n}}.$$
Thus $\mathcal{S}[\preceq]\mathcal{T}$ implies $\mathcal{L}^{\infty}_{[\mathcal{S}]}\subseteq\mathcal{L}^{\infty}_{[\mathcal{T}]}$ with continuous inclusion.

\item[$(ii)$] If both $\mathcal{S}$ and $\mathcal{T}$ are in addition continuous, then analogously $\mathcal{S}[\preceq]\mathcal{T}$ implies $\mathcal{L}^p_{[\mathcal{S}]}\subseteq\mathcal{L}^p_{[\mathcal{T}]}$ with continuous inclusion for some/any $1\le p<\infty$. (The only difference compared with $(i)$ is that in the estimates $e^C$ has to be replaced by $e^{C/p}$.)
\end{itemize}
\end{proposition}

The aim now is to verify the (expected) characterization of inclusion relations in terms of $(\preceq)$ and $\{\preceq\}$ between weight function matrices. For this we have to impose mild extra conditions on the weights; in Section \ref{generacharactsect} we present for both the Beurling- and the Roumieu-type an abstract functional analytic argument (see Propositions \ref{Beurlingcompprop} and \ref{Roumcompprop}). This approach, in particular the Beurling-type result, has been inspired by \cite[Thm. 1.3.18]{Bjorck66}. Afterwards, in Section \ref{alternativecharactsect}, we give for the Roumieu-type a different constructive approach involving special (optimal) functions.\vspace{6pt}

In the proofs of Propositions \ref{Beurlingcompprop} and \ref{Roumcompprop} the \emph{non-triviality} of $\mathcal{L}^{p}_{[\mathcal{S}]}$ is required and in the isotropic setting in view of Lemma \ref{nontriviallemma} it suffices to assume in addition that $\mathcal{S}$ is non-decreasing. However, by taking into account the convention from Remark \ref{infinityremark} we comment now on the fact that, if these spaces are trivial, then the implications follow automatically.

\begin{remark}\label{infinityremark1}
When $\{0\}=\mathcal{L}^{p}_{[\mathcal{S}]}\subseteq\mathcal{L}^{p}_{[\mathcal{T}]}$, $1\le p\le\infty$, then the inclusion is trivial and the convention from Remark \ref{infinityremark} gives that $\mathcal{S}=\{\sigma^{\ell}: \ell>0\}$ corresponds to the weight function matrix satisfying $\sigma^{\ell}\equiv+\infty$ for all $\ell>0$ sufficiently large in the Beurling-case resp. $\sigma^{\ell}\equiv+\infty$ for any $\ell>0$ in the Roumieu-case. And therefore, by recalling \eqref{growthrel1} resp. \eqref{growthrel2}, the relation $\mathcal{S}(\preceq)\mathcal{T}$ resp. $\mathcal{S}\{\preceq\}\mathcal{T}$ is formally trivial.

Finally, by using the convention $+\infty=+\infty$, the relations formally even follow when in addition $\mathcal{L}^{p}_{[\mathcal{T}]}=\{0\}$.
\end{remark}

\subsection{Proof of the main characterizing results}\label{generacharactsect}
We start with the \emph{Beurling-type.}

\begin{proposition}\label{Beurlingcompprop}
Let $\mathcal{S}=\{\sigma^{\ell}: \ell>0\}$ and $\mathcal{T}=\{\tau^{\ell}: \ell>0\}$ be weight function matrices and assume that both matrices satisfy \eqref{maincharthmbeurequ}.

\begin{itemize}
\item[$(i)$] If both matrices are also continuous, then the (continuous) inclusion $\mathcal{L}^{p}_{(\mathcal{S})}\subseteq\mathcal{L}^{p}_{(\mathcal{T})}$ (for some/any $1\le p<\infty$) implies $\mathcal{S}(\preceq)\mathcal{T}$.

\item[$(ii)$] The (continuous) inclusion $\mathcal{L}^{\infty}_{(\mathcal{S})}\subseteq\mathcal{L}^{\infty}_{(\mathcal{T})}$ implies $\mathcal{S}(\preceq)\mathcal{T}$.
\end{itemize}
\end{proposition}

\demo{Proof}
$(i)$ First, by \eqref{maincharthmbeurequ} applied to the matrix $\mathcal{T}$ we get:
\begin{align*}
&\forall\;\ell>0\;\exists\;\ell'>0\;\exists\;L\ge 1\;\forall\;y,x_0\in\RR^d:
\\&
\tau^{\ell}(x_0)=\tau^{\ell}(x_0+y-y)\le\tau^{\ell'}(x_0+y)+\tau^{\ell'}(-y)+L.
\end{align*}
Thus, for any $f\in\mathcal{L}^{p}_{(\mathcal{S})}\subseteq\mathcal{L}^{p}_{(\mathcal{T})}\subseteq L^p$ and $x_0\in\RR^d$ it holds that
$$\|f_{x_0}\|^p_{p,\tau^{\ell'}}=\int_{\RR^d}|f(y)|^pe^{\tau^{\ell'}(y+x_0)}dy\ge e^{-L}e^{\tau^{\ell}(x_0)}\int_{\RR^d}|f(y)|^pe^{-\tau^{\ell'}(-y)}dy,$$
and clearly $A_1:=\int_{\RR^d}|f(y)|^pe^{-\tau^{\ell'}(-y)}dy\le\int_{\RR^d}|f(y)|^pdy<+\infty$ by assumption.

Second, since both $\mathcal{L}^{p}_{(\mathcal{S})}$ and $\mathcal{L}^{p}_{(\mathcal{T})}$ are Fr\'{e}chet spaces (w.r.t. their natural topologies) we have that the inclusion $\mathcal{L}^{p}_{(\mathcal{S})}\subseteq\mathcal{L}^{p}_{(\mathcal{T})}$ is continuous (in any case); see e.g. \cite[Prop. 4.5 \& Rem. 4.6]{PTTvsmatrix}. Consequently,
\begin{equation}\label{Beurlingcomppropequ}
\forall\;\ell'>0\;\exists\;\ell''>0\;\exists\;A\ge 1\;\forall\;f\in\mathcal{L}^{p}_{(\mathcal{S})}:\;\;\;\|f\|_{p,\tau^{\ell'}}\le A\|f\|_{p,\sigma^{\ell''}}.
\end{equation}
Now, in view of Remark \ref{infinityremark1} we can assume that there exists $0\neq f\in\mathcal{L}^{p}_{(\mathcal{S})}$ and \eqref{maincharthmbeurequ} for $\mathcal{S}$ gives via $(i)$ in Proposition \ref{translationlemma} that $f_{x_0}\in\mathcal{L}^{p}_{(\mathcal{S})}$ for all $x_0\in\RR^d$. The idea is to apply \eqref{Beurlingcomppropequ} to the family $\{f_{x_0}: x_0\in\RR^d\}$.

Let $\ell>0$ be given, arbitrary but from now on fixed. Let $\ell'>0$ be the index appearing in the estimate in the first step and being related to $\ell$ via \eqref{maincharthmbeurequ} for the matrix $\mathcal{T}$. Then apply \eqref{Beurlingcomppropequ} to $\ell'$ and let $\ell''$ be the index subject to $\ell'$ in this property.

Finally, when applying \eqref{maincharthmbeurequ} to $\mathcal{S}$ and to the index $\ell''$ then there exist some $\ell'''>0$ and $L_1\ge 1$ such that for all $x_0\in\RR^d$:
\begin{align*}
\|f_{x_0}\|^p_{p,\sigma^{\ell''}}=\int_{\RR^d}|f(y)|^pe^{\sigma^{\ell''}(y+x_0)}dy\le e^{L_1}e^{\sigma^{\ell'''}(x_0)}\int_{\RR^d}|f(y)|^pe^{\sigma^{\ell'''}(y)}dy.
\end{align*}
Note that $A_2:=\int_{\RR^d}|f(y)|^pe^{\sigma^{\ell'''}(y)}dy<+\infty$ since $f\in\mathcal{L}^{p}_{(\mathcal{S})}$.

Combining everything, so far we have shown
\begin{align*}
&\forall\;\ell>0\;\exists\;\ell'''>0\;\exists\;A,L,L_1\ge 1\;\exists\;A_1,A_2>0\;\forall\;x_0\in\RR^d:
\\&
\tau^{\ell}(x_0)\le L-\log(A_1)+p\log(A)+L_1+\log(A_2)+\sigma^{\ell'''}(x_0);
\end{align*}
so \eqref{growthrel1} is verified between the indices $\ell$ and $\ell'''$ and with $C:=L-\log(A_1)+p\log(A)+L_1+\log(A_2)$; i.e. relation $\mathcal{S}(\preceq)\mathcal{T}$.\vspace{6pt}

$(ii)$ Similarly, since $\mathcal{L}^{\infty}_{(\sigma)}$ and $\mathcal{L}^{\infty}_{(\tau)}$ are Fr\'{e}chet spaces (w.r.t. their natural topology), by \cite[Prop. 4.5 \& Rem. 4.6]{PTTvsmatrix} we have again that the inclusion $\mathcal{L}^{\infty}_{(\mathcal{S})}\subseteq\mathcal{L}^{\infty}_{(\mathcal{T})}$ is continuous in any case. Thus
\begin{equation}\label{fourthweightcomppropequ1}
\forall\;\ell'>0\;\exists\;\ell''>0\;\exists\;A\ge 1\;\forall\;f\in\mathcal{L}^{\infty}_{(\mathcal{S})}:\;\;\;\|f\|_{\infty,\tau^{\ell'}}\le A\|f\|_{\infty,\sigma^{\ell''}},
\end{equation}
and we apply this crucial estimate again to the family $\{f_{x_0}: x_0\in\RR^d\}$ corresponding to a fixed $0\neq f\in\mathcal{L}^{\infty}_{(\sigma)}$; recall Remark \ref{infinityremark1}. The rest follows analogously as in $(i)$ by involving \eqref{maincharthmbeurequ} and $(i)$ in Proposition \ref{translationlemma} for both matrices and note that $A_1:=\ess\sup_{y\in\RR^d}|f(y)|e^{-\tau^{\ell}(-y)}\le\ess\sup_{y\in\RR^d}|f(y)|<+\infty$ for all $\ell>0$ because $f\in\mathcal{L}^{\infty}_{(\mathcal{T})}$ and $A_2:=\ess\sup_{x\in\RR^d}|f(y)|e^{\sigma^{\ell'''}(y)}<+\infty$ for any $\ell'''>0$ because even $f\in\mathcal{L}^{\infty}_{(\mathcal{S})}$ (by assumption). Therefore, $\mathcal{S}(\preceq)\mathcal{T}$ follows as before.
\qed\enddemo

The next result deals with the \emph{Roumieu case.}

\begin{proposition}\label{Roumcompprop}
Let $\mathcal{S}=\{\sigma^{\ell}: \ell>0\}$ and $\mathcal{T}=\{\tau^{\ell}: \ell>0\}$ be given weight function matrices and assume that both matrices satisfy \eqref{maincharthmroumequ}.
	
	\begin{itemize}
		\item[$(i)$] If both matrices are also continuous, then the (continuous) inclusion $\mathcal{L}^{p}_{\{\mathcal{S}\}}\subseteq\mathcal{L}^{p}_{\{\mathcal{T}\}}$ for some/any $1\le p<\infty$ implies $\mathcal{S}\{\preceq\}\mathcal{T}$.
		
		\item[$(ii)$] The (continuous) inclusion $\mathcal{L}^{\infty}_{\{\mathcal{S}\}}\subseteq\mathcal{L}^{\infty}_{\{\mathcal{T}\}}$ implies $\mathcal{S}\{\preceq\}\mathcal{T}$.
	\end{itemize}
\end{proposition}

\demo{Proof}
$(i)$ We proceed analogously as in Proposition \ref{Beurlingcompprop}: By the (continuous) inclusion $\mathcal{L}^{p}_{\{\mathcal{S}\}}\subseteq\mathcal{L}^{p}_{\{\mathcal{T}\}}$ and since we are dealing with (LB)-spaces, by \emph{de Wilde's closed graph theorem and Grothendieck's factorization theorem} we obtain
\begin{equation}\label{fourthweightcomppropequroum}
\forall\;\ell>0\;\exists\;\ell'>0\;\exists\;A\ge 1\;\forall\;f\in\mathcal{L}^{p}_{\sigma^{\ell}}:\;\;\;\|f\|_{p,\tau^{\ell'}}\le A\|f\|_{p,\sigma^{\ell}}.
\end{equation}
Fix an index $\ell_0>0$ and $0\neq f\in\mathcal{L}^{p}_{\sigma^{\ell_0}}$ (recall Remark \ref{infinityremark1}). Moreover, assumption \eqref{maincharthmroumequ} for $\mathcal{S}$ and $(ii)$ in Proposition \ref{translationlemma} yield the existence of some $\ell>0$ such that $f_{x_0}\in\mathcal{L}^{p}_{\sigma^{\ell}}$ for any $x_0\in\RR^d$. (The indices $\ell_0$ and $\ell$ are related by \eqref{maincharthmroumequ}.) Indeed, the proof of $(ii)$ in Proposition \ref{translationlemma} gives
\begin{align*}
\|f_{x_0}\|^p_{p,\sigma^{\ell}}\le e^{L}e^{\sigma^{\ell_0}(x_0)}\int_{\RR^d}|f(y)|^pe^{\sigma^{\ell_0}(y)}dy.
\end{align*}
Note that $A_1:=\int_{\RR^d}|f(y)|^pe^{\sigma^{\ell_0}(y)}dy<+\infty$ since $f\in\mathcal{L}^{p}_{\sigma^{\ell_0}}$. We apply \eqref{fourthweightcomppropequroum} to the family $\{f_{x_0}: x_0\in\RR^d\}$ and to this index $\ell$. Moreover, when using \eqref{maincharthmroumequ} for $\mathcal{T}$, then analogously as in the proof above for some index $\ell'$ (being related to $\ell$ via \eqref{fourthweightcomppropequroum}) we can find some index $\ell''>0$ and $L_1\ge 1$ and get for all $x_0\in\RR^d$:
$$\|f_{x_0}\|^p_{p,\tau^{\ell'}}=\int_{\RR^d}|f(y)|^pe^{\tau^{\ell'}(y+x_0)}dy\ge e^{-L_1}e^{\tau^{\ell''}(x_0)}\int_{\RR^d}|f(y)|^pe^{-\tau^{\ell'}(-y)}dy.$$
Clearly $A_2:=\int_{\RR^d}|f(y)|^pe^{-\tau^{\ell'}(-y)}dy\le\int_{\RR^d}|f(y)|^pdy<+\infty$ since $f\in\mathcal{L}^{p}_{\{\mathcal{T}\}}$ (by assumption). Combining everything, so far we have shown
\begin{align*}
	&\forall\;\ell_0>0\;\exists\;\ell''>0\,\exists\;A,L,L_1\ge 1\;\exists\;A_1,A_2>0\;\forall\;x_0\in\RR^d:
	\\&
	\tau^{\ell''}(x_0)\le L_1-\log(A_2)+p\log(A)+L+\log(A_1)+\sigma^{\ell_0}(x_0),
\end{align*}
which yields \eqref{growthrel2} and so $\mathcal{S}\{\preceq\}\mathcal{T}$.\vspace{6pt}

$(ii)$ Similarly, by the continuous inclusion $\mathcal{L}^{\infty}_{\{\mathcal{S}\}}\subseteq\mathcal{L}^{\infty}_{\{\mathcal{T}\}}$ and since we are dealing with (LB)-spaces we obtain
\begin{equation}\label{fourthweightcomppropequ1roum}
	\forall\;\ell>0\;\exists\;\ell'>0\;\exists\;A\ge 1\;\forall\;f\in\mathcal{L}^{\infty}_{\sigma^{\ell}}:\;\;\;\|f\|_{\infty,\tau^{\ell'}}\le A\|f\|_{\infty,\sigma^{\ell}}.
\end{equation}
Again we fix an index $\ell_0>0$ and $0\neq f\in\mathcal{L}^{\infty}_{\sigma^{\ell_0}}$ (recall Remark \ref{infinityremark1}) and $(ii)$ in Proposition \ref{translationlemma} yields the existence of some $\ell>0$ such that $f_{x_0}\in\mathcal{L}^{p}_{\sigma^{\ell}}$ for any $x_0\in\RR^d$. (Again the indices $\ell_0$ and $\ell$ are related by \eqref{maincharthmroumequ}.)

We apply \eqref{fourthweightcomppropequ1roum} to the family $\{f_{x_0}: x_0\in\RR^d\}$ and to this index $\ell$. Then follow the estimates in $(i)$ using the same notation (for the indices) and by involving $(ii)$ in Proposition \ref{translationlemma} and \eqref{maincharthmroumequ}; see also $(ii)$ in Proposition \ref{Beurlingcompprop}. Therefore $\mathcal{S}\{\preceq\}\mathcal{T}$ holds and note that here $A_1:=\ess\sup_{x\in\RR^d}|f(y)|e^{\sigma^{\ell_0}(y)}<+\infty$ because $f\in\mathcal{L}^{\infty}_{\sigma^{\ell_0}}$ and $A_2:=\ess\sup_{y\in\RR^d}|f(y)|e^{-\tau^{\ell'}(-y)}\le\ess\sup_{y\in\RR^d}|f(y)|<+\infty$ since $f\in\mathcal{L}^{\infty}_{\{\mathcal{T}\}}$ (by assumption).
\qed\enddemo

\subsection{Alternative proofs for the Roumieu-type}\label{alternativecharactsect}
We present now different techniques and ideas for the Roumieu case to show $(ii)\Rightarrow(i)$ in Theorem \ref{maincharthmroum}. Here we are not using functional analytic tools but involve special (optimal) functions. One has to work with different assumptions for the weight function matrices but the new introduced optimal functions can be used in other related contexts as well.

Let $\mathcal{W}=\{\omega^{\ell}: \ell>0\}$ be a weight function matrix. Then set
\begin{equation}\label{charactfctdef}
\theta^p_{\omega^{\ell}}(x):=e^{-p^{-1}\omega^{\ell}(x)},\;\;\;x\in\RR^d,\;\ell>0,\;1\le p<\infty,
\end{equation}
and
\begin{equation}\label{charactfctdef1}
\theta^{\infty}_{\omega^{\ell}}(x):=e^{-\omega^{\ell}(x)},\;\;\;x\in\RR^d,\;\ell>0.
\end{equation}
The extreme case $\omega^{\ell}=0$ gives $\theta^p_{\omega^{\ell}}=1=\theta^{\infty}_{\omega^{\ell}}$. In order to proceed, we require the following relations between the functions $\omega^{\ell}$:
\begin{equation}\label{weakdifferentgrowth}
\exists\;\ell'_0>\ell_0>0:\;\;\;\sup_{x\in\RR^d}\omega^{\ell'_0}(x)-\omega^{\ell_0}(x)=+\infty,
\end{equation}
\begin{equation}\label{weakdifferentgrowth1}
\exists\;\ell'_0>\ell_0>0:\;\;\;\liminf_{|x|\rightarrow+\infty}\omega^{\ell'_0}(x)-\omega^{\ell_0}(x)>0,
\end{equation}
and
\begin{equation}\label{strongdifferentgrowth}
\exists\;a>1\;\forall\;\ell>0\;\exists\;\ell'>0\;\exists\;b\in\RR\;\forall\;x\in\RR^d:\;\;\;\omega^{\ell}(x)-\omega^{\ell'}(x)\ge a\log(1+|x|)+b.
\end{equation}

\begin{remark}\label{strongdifferentgrowthrem}
We point out:
\begin{itemize}
	\item[$(*)$] If \eqref{weakdifferentgrowth}, \eqref{weakdifferentgrowth1} holds for some indices $\ell'_0>\ell_0$, then by the order of the elements in the matrix also for all $\ell',\ell$ satisfying $\ell'\ge\ell'_0>\ell_0\ge\ell>0$.
	
	\item[$(*)$] By iteration it is clear that \eqref{strongdifferentgrowth} is equivalent to
\begin{equation}\label{strongdifferentgrowthvar}	
\forall\;a>1\;\forall\;\ell>0\;\exists\;\ell'>0\;\exists\;b\in\RR\;\forall\;x\in\RR^d:\;\;\;\omega^{\ell}(x)-\omega^{\ell'}(x)\ge a\log(1+|x|)+b.
\end{equation}
Note that in \eqref{strongdifferentgrowth}, \eqref{strongdifferentgrowthvar} necessarily the indices have to satisfy the relation $\ell>\ell'$.
\item[$(*)$] One shall note that \eqref{weakdifferentgrowth}, \eqref{weakdifferentgrowth1} and \eqref{strongdifferentgrowth} are crucial for the techniques in this section but they never can be satisfied for any \emph{simple} weight function matrix $\mathcal{W}:=\{\omega\}$. On the other hand, \eqref{maincharthmbeurequ} and \eqref{maincharthmroumequ} hold for any sub-additive $\omega$ (and coincide in the simple case) and hence the main results Theorems \ref{maincharthmbeur}, \ref{maincharthmroum} and Proposition \ref{translationlemma} can be applied, too.
\end{itemize}
\end{remark}

\begin{lemma}\label{thetafctlemma}
Let $\mathcal{W}=\{\omega^{\ell}: \ell>0\}$ be a weight function matrix.
\begin{itemize}
\item[$(i)$] We have
$$\forall\;\ell\ge\ell'>0:\;\;\;\theta^{\infty}_{\omega^{\ell}}\in\mathcal{L}^{\infty}_{\omega^{\ell'}}\subseteq\mathcal{L}^{\infty}_{\{\mathcal{W}\}}.$$

\item[$(ii)$] If $\mathcal{W}$ satisfies \eqref{weakdifferentgrowth} with indices $\ell'_0>\ell_0$, then
$$\forall\;\ell'\ge\ell'_0>\ell_0\ge\ell>0:\;\;\;\theta^{\infty}_{\omega^{\ell}}\notin\mathcal{L}^{\infty}_{\omega^{\ell'}}\supseteq\mathcal{L}^{\infty}_{(\mathcal{W})}.$$

\item[$(iii)$] Assume that $\mathcal{W}$ is continuous and \eqref{strongdifferentgrowth} holds. Then for any $1\le p<\infty$ we get:
$$\forall\;\ell>0\;\exists\;\ell'>0:\;\;\;\theta^p_{\omega^{\ell}}\in\mathcal{L}^p_{\omega^{\ell'}}\subseteq\mathcal{L}^p_{\{\mathcal{W}\}}.$$

\item[$(iv)$] Assume that $\mathcal{W}$ is continuous and \eqref{weakdifferentgrowth1} holds. Then for any $1\le p<\infty$ we get:
$$\exists\;\ell_0'>\ell_0\;\forall\;\ell'\ge\ell'_0>\ell_0\ge\ell>0:\;\;\;\theta^p_{\omega^{\ell}}\notin\mathcal{L}^p_{\omega^{\ell'}}\supseteq\mathcal{L}^p_{(\mathcal{W})}.$$
\end{itemize}
\end{lemma}

\demo{Proof}
$(i)$ This is clear since $\sup_{x\in\RR^d}|\theta^{\infty}_{\omega^{\ell}}(x)|e^{\omega^{\ell'}(x)}=\sup_{x\in\RR^d}e^{\omega^{\ell'}(x)-\omega^{\ell}(x)}\le 1$.

$(ii)$ According to the given indices and the assumption it holds that $\sup_{x\in\RR^d}|\theta^{\infty}_{\omega^{\ell}}(x)|e^{\omega^{\ell'}(x)}=\sup_{x\in\RR^d}e^{\omega^{\ell'}(x)-\omega^{\ell}(x)}\ge\sup_{x\in\RR^d}e^{\omega^{\ell'_0}(x)-\omega^{\ell_0}(x)}=+\infty$; see Remark \ref{strongdifferentgrowthrem}.

$(iii)$ We apply \eqref{strongdifferentgrowthvar} to (some fixed) $a\ge 2(d+1)$ and hence get for each $\ell>0$ some index $\ell'>0$ such that
$$\|\theta^p_{\omega^{\ell}}\|^p_{p,\omega^{\ell'}}=\int_{\RR^d}e^{\omega^{\ell'}(x)-\omega^{\ell}(x)}dx\le e^{-b}\int_{\RR^d}\frac{1}{(1+|x|)^{a}}dx.$$
The last integral is finite by the choice $a\ge 2(d+1)$.

$(iv)$ Let $\ell'\ge\ell'_0>\ell_0\ge\ell>0$ be given, with $\ell'_0$ and $\ell_0$ denoting the parameters appearing in \eqref{weakdifferentgrowth1} and then
$$\|\theta^p_{\omega^{\ell}}\|^p_{p,\omega^{\ell'}}=\int_{\RR^d}e^{\omega^{\ell'}(x)-\omega^{\ell}(x)}dx\ge\int_{\RR^d}e^{\omega^{\ell'_0}(x)-\omega^{\ell_0}(x)}dx=+\infty,$$
and the last equality is valid since by assumption there exists some $\epsilon>0$ and $t_{\epsilon}>0$ such that $\omega^{\ell'_0}(x)-\omega^{\ell_0}(x)\ge\epsilon$ for all $x\in\RR^d$ satisfying $|x|\ge t_{\epsilon}$. However, in general it is not clear if already \eqref{weakdifferentgrowth} implies this estimate; see Remark \ref{peakremark}.
\qed\enddemo

\emph{Note:}

\begin{itemize}
\item[$(*)$] $(i)$ and $(iii)$ in Lemma \ref{thetafctlemma} yield, in particular, the non-triviality of the corresponding weighted Roumieu-type spaces even in the anisotropic setting. This should be compared with Lemma \ref{nontriviallemma}: There it has been assumed that $\mathcal{W}$ is non-decreasing, whereas $(i)$ in Lemma \ref{thetafctlemma} requires no further assumption on $\mathcal{W}$ and $(iii)$ there requires information on the difference between the growth of the weight functions which fails in the simple case; recall Remark \ref{strongdifferentgrowthrem}.

\item[$(*)$] By $(ii)$ and $(iv)$ in Lemma \ref{thetafctlemma} we see that in general $\theta^{\infty}_{\omega^{\ell}}$ resp. $\theta^p_{\omega^{\ell}}$ does not belong to the corresponding Beurling-type spaces and so the proofs of the next results are becoming unclear in this setting.
\end{itemize}

Using this preparation we prove the following variant of $(ii)$ in Proposition \ref{Roumcompprop}.

\begin{proposition}\label{thirdweightcompprop}
Let $\mathcal{S}=\{\sigma^{\ell}: \ell>0\}$ and $\mathcal{T}=\{\tau^{\ell}: \ell>0\}$ be weight function matrices.
	
Then the inclusion $\mathcal{L}^{\infty}_{\{\mathcal{S}\}}\subseteq\mathcal{L}^{\infty}_{\{\mathcal{T}\}}$ (as sets) implies $\mathcal{S}\{\preceq\}\mathcal{T}$.
\end{proposition}

\demo{Proof}
By $(i)$ in Lemma \ref{thetafctlemma} and $\mathcal{L}^{\infty}_{\{\mathcal{S}\}}\subseteq\mathcal{L}^{\infty}_{\{\mathcal{T}\}}$ we get $\theta^{\infty}_{\sigma^{\ell}}\in\mathcal{L}^{\infty}_{\{\mathcal{T}\}}$ for any $\ell>0$. Hence
$$\forall\;\ell>0\;\exists\;\ell'>0:\;\;\;\sup_{x\in\RR^d}e^{\tau^{\ell'}(x)-\sigma^{\ell}(x)}<+\infty,$$
and so
$$\forall\;\ell>0\;\exists\;\ell'>0\;\exists\;D\ge 0\;\forall\;x\in\RR^d:\;\;\;\tau^{\ell'}(x)\le\sigma^{\ell}(x)+D,$$
which is precisely $\mathcal{S}\{\preceq\}\mathcal{T}$.
\qed\enddemo

\emph{Note:} This statement is more general than $(ii)$ in Proposition \ref{Roumcompprop} since no stability under translation (i.e. condition \eqref{maincharthmroumequ}) is required.

Moreover, the established technique immediately gives the next result:

\begin{proposition}\label{thirdstrongweightcompprop}
Let $\mathcal{S}=\{\sigma^{\ell}: \ell>0\}$ and $\mathcal{T}=\{\tau^{\ell}: \ell>0\}$ be weight function matrices. Then the following are equivalent:
\begin{itemize}
\item[$(i)$] We have the continuous inclusion $\mathcal{L}^{\infty}_{\{\mathcal{S}\}}\subseteq\mathcal{L}^{\infty}_{(\mathcal{T})}$.

\item[$(ii)$] We have the inclusion $\mathcal{L}^{\infty}_{\{\mathcal{S}\}}\subseteq\mathcal{L}^{\infty}_{(\mathcal{T})}$ as sets.

\item[$(iii)$] The matrices satisfy
\begin{equation}\label{thirdstrongweightcomppropequ}
\forall\;\ell>0\;\forall\;\ell'>0\;\exists\;D\ge 1\;\forall\;x\in\RR^d:\;\;\;\tau^{\ell}(x)\le\sigma^{\ell'}(x)+D;
\end{equation}
and let us abbreviate this relation by $\mathcal{S}\vartriangleleft\mathcal{T}$.
\end{itemize}
\end{proposition}

\demo{Proof}
$(iii)\Rightarrow(i)$ follows immediately by definition of the spaces and $(i)\Rightarrow(ii)$ is clear. For $(ii)\Rightarrow(iii)$ follow the proof of Proposition \ref{thirdweightcompprop} and replace ``$\exists\;\ell'>0$'' by ``$\forall\;\ell'>0$''.
\qed\enddemo

Finally, we treat the $L^p$-Roumieu-type spaces and so a variant of $(i)$ in Proposition \ref{Roumcompprop}. Here, due to technical reasons we focus on isotropic weight matrices and, before stating the next result, we introduce the following growth relations for $\mathcal{W}=\{\omega^{\ell}: \ell>0\}$:

\begin{equation}\label{weakom1}
\forall\;\ell>0\;\exists\;\ell_1>0\;\exists\;C\ge 1\;\forall\;t\in\NN:\;\;\;\omega^{\ell_1}(t+1)\le\omega^{\ell}(t)+C,
\end{equation}
\begin{equation}\label{weakom1var}
\forall\;\ell_1>0\;\exists\;\ell>0\;\exists\;C\ge 1\;\forall\;t\in\NN:\;\;\;\omega^{\ell_1}(t+1)\le\omega^{\ell}(t)+C.
\end{equation}

\begin{itemize}
\item[$(*)$] Similarly, as pointed out in Remark \ref{translationlemmarem} (see \eqref{mixedom1roum} resp. \eqref{mixedom1beur}), one could consider the expression $C\omega^{\ell}(t)+C$ on the right-hand side. However, the additional appearance of the multiplicative constant $C$ is not appropriate in general; see again Section \ref{specialcasesect1}.

\item[$(*)$] \eqref{maincharthmroumequ} obviously implies \eqref{weakom1} and \eqref{maincharthmbeurequ} implies \eqref{weakom1var}.
	
\item[$(*)$] If $\mathcal{W}$ is non-decreasing, then equivalently in \eqref{weakom1} and \eqref{weakom1var} we can replace ``$\forall\;t\in\NN$'' by ``$\forall\;t\ge 0$''; we only treat \eqref{weakom1}:

     When $t\ge 0$ with $t\notin\NN$ is given we chose $s_t\in\NN$ such that $s_t<t<s_t+1$ and get by iteration $\omega^{\ell_2}(s_t+1)\le\omega^{\ell_2}(s_t+2)\le\omega^{\ell_1}(s_t+1)+C\le\omega^{\ell}(s_t)+C_1+C\le\omega^{\ell}(t)+C_1+C$. Hence it suffices to replace $C$ by $C_1+C$ and $\ell_1$ by $\ell_2(\le\ell_1)$ to ensure \eqref{weakom1} for all $t\ge 0$.
\end{itemize}

\begin{proposition}\label{secondweightcompprop}
Let $\mathcal{S}=\{\sigma^{\ell}: \ell>0\}$ and $\mathcal{T}=\{\tau^{\ell}: \ell>0\}$ be continuous and non-decreasing (and hence radial) weight function matrices. Moreover, assume that $\mathcal{S}$ satisfies \eqref{strongdifferentgrowth} and that either $\mathcal{S}$ or $\mathcal{T}$ has \eqref{weakom1}.

Then the inclusion $\mathcal{L}^p_{\{\mathcal{S}\}}\subseteq\mathcal{L}^p_{\{\mathcal{T}\}}$ (as sets) for some/any $1\le p<\infty$ implies $\mathcal{S}\{\preceq\}\mathcal{T}$.
\end{proposition}

\demo{Proof}
By assumption \eqref{strongdifferentgrowth}, $(iii)$ in Lemma \ref{thetafctlemma} and the inclusion $\mathcal{L}^p_{\{\mathcal{S}\}}\subseteq\mathcal{L}^p_{\{\mathcal{T}\}}$ (as sets) we obtain $\theta^p_{\sigma^{\ell}}\in\mathcal{L}^p_{\{\mathcal{T}\}}$ for all $\ell>0$ and so
\begin{equation}\label{secondweightcomppropequ}
\forall\;\ell>0\;\exists\;\ell_1>0:\;\;\;\int_{\RR^d}|\theta^p_{\sigma^{\ell}}(x)|^pe^{\tau^{\ell_1}(x)}dx=\int_{\RR^d}e^{\tau^{\ell_1}(x)-\sigma^{\ell}(x)}dx<+\infty.
\end{equation}
Let now $x\in\RR^d$ with $n\le|x|<n+1$ for some $n\in\NN$. Then $e^{\tau^{\ell_1}(x)-\sigma^{\ell}(x)}\ge e^{\tau^{\ell_1}(n)-\sigma^{\ell}(n+1)}$ by using the assumption that both matrices are non-decreasing. Next recall that the volume of the $d$-dimensional sphere with radius $r$ is given by $V_d(r)=\frac{\pi^{d/2}r^d}{\Gamma(1+d/2)}$ (with $\Gamma$ denoting the Gamma-function) and therefore,
\begin{align*}
\int_{x\in\RR^d: n\le|x|<n+1}1dx&=\int_{x\in\RR^d: |x|<n+1}1dx-\int_{x\in\RR^d: |x|\le n}1dx=\frac{\pi^{d/2}((n+1)^d-n^d)}{\Gamma(1+d/2)}
\\&
\ge\frac{\pi^{d/2}}{\Gamma(1+d/2)}=:C_d,
\end{align*}
because $(n+1)^d-n^d\ge 1$ for all $n\in\NN$ and $d\in\NN_{>0}$. Then we estimate as follows:
\begin{align*}
\int_{\RR^d}e^{\tau^{\ell_1}(x)-\sigma^{\ell}(x)}dx&=\sum_{n=0}^{+\infty}\int_{x\in\RR^d: n\le|x|<n+1}e^{\tau^{\ell_1}(x)-\sigma^{\ell}(x)}dx\ge\sum_{n=0}^{+\infty}e^{\tau^{\ell_1}(n)-\sigma^{\ell}(n+1)}\int_{x\in\RR^d: n\le|x|<n+1}1dx
\\&
\ge C_d\sum_{n=0}^{+\infty}e^{\tau^{\ell_1}(n)-\sigma^{\ell}(n+1)}.
\end{align*}
Combining this estimate with \eqref{secondweightcomppropequ} we get, in particular, that $\lim_{n\rightarrow+\infty}e^{\tau^{\ell_1}(n)-\sigma^{\ell}(n+1)}=0$ or equivalently $\lim_{n\rightarrow+\infty}\tau^{\ell_1}(n)-\sigma^{\ell}(n+1)=-\infty$. Summarizing, so far we have shown
\begin{equation}\label{secondweightcomppropequ1}
\forall\;\ell>0\;\exists\;\ell_1>0\;\exists\;D\ge 0\;\forall\;n\in\NN:\;\;\;\tau^{\ell_1}(n)\le\sigma^{\ell}(n+1)+D.
\end{equation}
Now assume that \eqref{weakom1} for $\mathcal{S}$ holds and let $x\in\RR^d$ with $n\le|x|<n+1$ for some $n\in\NN$. Then when iterating this property we estimate by
\begin{align*}
\tau^{\ell_1}(x)&\le\tau^{\ell_1}(n+1)\le\sigma^{\ell}(n+2)+D\le\sigma^{\ell'}(n+1)+C+D\le\sigma^{\ell''}(n)+C_1+C+D
\\&
\le\sigma^{\ell''}(x)+C_1+C+D,
\end{align*}
thus $\mathcal{S}\{\preceq\}\mathcal{T}$ is verified between the indices $\ell''(\ge\ell)$ and $\ell_1$.

If $\mathcal{T}$ satisfies \eqref{weakom1}, then similarly by \eqref{secondweightcomppropequ1} for any $x\in\RR^d$ with $n\le|x|<n+1$ for some $n\in\NN_{>0}$ we have
\begin{align*}
\tau^{\ell''_1}(x)&\le\tau^{\ell''_1}(n+1)\le\tau^{\ell'_1}(n)+C\le\tau^{\ell_1}(n-1)+C_1+C\le\sigma^{\ell}(n)+C_1+C+D
\\&
\le\sigma^{\ell}(x)+C_1+C+D.
\end{align*}
And for $x\in\RR^d$ with $0\le|x|<1$ one has $\tau^{\ell''_1}(x)\le\tau^{\ell''_1}(1)\le\sigma^{\ell}(x)+A$ with $A:=\tau^{\ell''_1}(1)$. Summarizing, again $\mathcal{S}\{\preceq\}\mathcal{T}$ is verified between the indices $\ell$ and $\ell''_1(\le\ell_1)$.
\qed\enddemo

\begin{remark}\label{peakremark}
Without the additional assumptions on the weight matrices it is not clear that \eqref{secondweightcomppropequ} implies
$\lim_{|x|\rightarrow+\infty}e^{\tau^{\ell_1}(x)-\sigma^{\ell}(x)}=0$ or, what is enough to conclude, $\sup_{x\in\RR^d}e^{\tau^{\ell_1}(x)-\sigma^{\ell}(x)}<+\infty$. This is due to the fact that in general $t\mapsto e^{\tau^{\ell_1}(t)-\sigma^{\ell}(t)}$ might behave very irregular: This function can have (infinitely many) ``exploding peaks'' but being sufficiently small ``nearly everywhere on $\RR$'' and hence \eqref{secondweightcomppropequ} still can be valid. In view of Remark \ref{strongdifferentgrowthrem} this comment becomes relevant in the simple case when no further index is appearing.
\end{remark}

\begin{proposition}\label{secondstrongweightcompprop}
Let $\mathcal{S}=\{\sigma^{\ell}: \ell>0\}$ and $\mathcal{T}=\{\tau^{\ell}: \ell>0\}$ be continuous and non-decreasing (and hence radial) weight function matrices. Moreover, assume that $\mathcal{S}$ satisfies \eqref{strongdifferentgrowth} and \eqref{weakom1}. Then the following are equivalent:
\begin{itemize}
\item[$(i)$] We have the continuous inclusion $\mathcal{L}^p_{\{\mathcal{S}\}}\subseteq\mathcal{L}^p_{(\mathcal{T})}$ for some/any $1\le p<\infty$.

\item[$(ii)$] We have the inclusion $\mathcal{L}^p_{\{\mathcal{S}\}}\subseteq\mathcal{L}^p_{(\mathcal{T})}$ as sets for some/any $1\le p<\infty$.

\item[$(iii)$] The matrices satisfy $\mathcal{S}\vartriangleleft\mathcal{T}$.
\end{itemize}
The same equivalences are valid if $\mathcal{S}$ is assumed as above, but not necessarily having \eqref{weakom1}, and $\mathcal{T}$ satisfies \eqref{weakom1var} instead.
\end{proposition}

\demo{Proof}
$(iii)\Rightarrow(i)$ follows by definition of the spaces and $(i)\Rightarrow(ii)$ is clear. For $(ii)\Rightarrow(iii)$ replace ``$\exists\;\ell_1>0$'' by ``$\forall\;\ell_1>0$'' in the proof of Proposition \ref{secondweightcompprop} and so get the stronger variant of \eqref{secondweightcomppropequ1}. Then involve in the last step \eqref{weakom1} for $\mathcal{S}$ resp. \eqref{weakom1var} for $\mathcal{T}$.
\qed\enddemo

\section{Special cases of weight function matrices}\label{specialcasesect}
We investigate two special but crucial situations for the defining (non-simple and non-stabilizing) weight function matrix $\mathcal{W}:=\{\omega^{\ell}: \ell>0\}$. Let $\omega: \RR^d\rightarrow[0,+\infty)$ be a fixed function and consider
\begin{equation}\label{expotype}
	\omega^{\ell}(x):=\ell\omega(x),\;\;\;\ell>0,\;x\in\RR^d,
\end{equation}
resp.
\begin{equation}\label{dilatype}
	\omega^{\ell}(x):=\omega(\ell x),\;\;\;\ell>0,\;x\in\RR^d.
\end{equation}
If \eqref{expotype} is valid, then $\mathcal{W}$ is called to be of \emph{exponential-type} and write $\mathcal{W}^{\mathfrak{c}}_{\omega}$ instead of $\mathcal{W}$. If \eqref{dilatype} is valid, then $\mathcal{W}$ is called to be of \emph{dilatation-type} and we write $\mathcal{W}_{\omega,\mathfrak{c}}$.

This terminology is inspired by the notation introduced in \cite{weightedentireinclusion1} and \cite{weightedentireinclusion2} and by taking into account the fact that the crucial seminorms are expressed by the function $\exp\circ\omega^{\ell}$. Note that known examples in the literature are given in terms of an exponential-type (and isotropic) weight matrix; see Section \ref{examplesect} for more details.

\subsection{On the exponential-type case}\label{specialcasesect1}

$\mathcal{W}^{\mathfrak{c}}_{\omega}$ is a weight function matrix according to Definition \ref{functionmatrixdef}. We gather some immediate facts and conclusions:

\begin{itemize}
	\item[$(a)$] $\mathcal{W}^{\mathfrak{c}}_{\omega}$ is non-decreasing/continuous/radial if and only if $\omega$ is so.
	
	\item[$(b)$] The following are equivalent:
	\begin{itemize}
		\item[$(*)$] $\sup_{x\in\RR^d}\omega(x)<+\infty$,
		
		\item[$(*)$] \eqref{boundedweightequroum} holds,
		
		\item[$(*)$] \eqref{boundedweightequbeur} holds.
	\end{itemize}
	Therefore, $\mathcal{W}^{\mathfrak{c}}_{\omega}$ is unbounded if and only if $\sup_{x\in\RR^d}\omega(x)=+\infty$.
	
	\item[$(c)$] Let $\omega$ be non-decreasing and radial, then the following are equivalent:
	\begin{itemize}
		\item[$(*)$] $\omega$ satisfies \hyperlink{om1}{$(\omega_1)$}.
		
        \item[$(*)$] $\mathcal{W}^{\mathfrak{c}}_{\omega}$ satisfies
        $$\exists\;n_0,\ell_0>0\;\exists\;C\ge 1\;\forall\;t\ge 0:\;\;\;\omega^{\ell_0}(2t)\le\omega^{n_0}(t)+C.$$

        \item[$(*)$] $\mathcal{W}^{\mathfrak{c}}_{\omega}$ satisfies \eqref{mixedom1refomequ2}.

        \item[$(*)$] $\mathcal{W}^{\mathfrak{c}}_{\omega}$ satisfies \eqref{mixedom1refomequ5}.

		\item[$(*)$] $\mathcal{W}^{\mathfrak{c}}_{\omega}$ satisfies \eqref{mixedom1refomequ1}.

        \item[$(*)$] $\mathcal{W}^{\mathfrak{c}}_{\omega}$ satisfies \eqref{mixedom1refomequ3}.

		\item[$(*)$] $\mathcal{W}^{\mathfrak{c}}_{\omega}$ satisfies \eqref{mixedom1refomequ4}.

        \item[$(*)$] $\mathcal{W}^{\mathfrak{c}}_{\omega}$ satisfies \eqref{mixedom1beur}.

        \item[$(*)$] $\mathcal{W}^{\mathfrak{c}}_{\omega}$ satisfies \eqref{mixedom1roum}.
	\end{itemize}
	 These equivalences follow by definition and in view of Lemma \ref{mixedom1refom}; recall $(f)$ in Section \ref{relevantgrowthsect}.

	\item[$(d)$] The following are equivalent:
	\begin{itemize}
		\item[$(*)$] $\mathcal{W}^{\mathfrak{c}}_{\sigma}(\preceq)\mathcal{W}^{\mathfrak{c}}_{\tau}$ holds, i.e. \eqref{growthrel1}.
		
		\item[$(*)$]  $\mathcal{W}^{\mathfrak{c}}_{\sigma}\{\preceq\}\mathcal{W}^{\mathfrak{c}}_{\tau}$ holds, i.e. \eqref{growthrel2}.
		
\item[$(*)$] The matrices $\mathcal{W}^{\mathfrak{c}}_{\sigma}$ and $\mathcal{W}^{\mathfrak{c}}_{\tau}$ are related by
\begin{equation}\label{matrixweakrelation}
\exists\;\ell_0>0\;\exists\;n_0>0\;\exists\;C\ge 1\;\forall\;x\in\RR^d:\;\;\;\tau^{\ell_0}(x)\le\sigma^{n_0}(x)+C.
\end{equation}

\item[$(*)$] $\sigma\hyperlink{ompreceq}{\preceq}\tau$ holds.
	\end{itemize}
For this note that obviously both \eqref{growthrel1} and \eqref{growthrel2} imply \eqref{matrixweakrelation}, whereas $\sigma\hyperlink{ompreceq}{\preceq}\tau$ implies that the indices satisfy the relation $n=\ell C$ for any given $\ell>0$ resp. $\ell=\frac{n}{C}$ for any given $n$, $C\ge 1$ denoting the constant from relation \hyperlink{ompreceq}{$\preceq$}, and hence both $\mathcal{W}^{\mathfrak{c}}_{\sigma}(\preceq)\mathcal{W}^{\mathfrak{c}}_{\tau}$ and $\mathcal{W}^{\mathfrak{c}}_{\sigma}\{\preceq\}\mathcal{W}^{\mathfrak{c}}_{\tau}$.

Moreover, by analogous reasons, the following are equivalent:
\begin{itemize}
\item[$(*)$] $\mathcal{W}^{\mathfrak{c}}_{\sigma}\vartriangleleft\mathcal{W}^{\mathfrak{c}}_{\tau}$,

\item[$(*)$] $$\exists\;n_0>0\;\forall\;\ell>0\;\exists\;C\ge 1\;\forall\;x\in\RR^d:\;\;\;\tau^{\ell}(x)\le\sigma^{n_0}(x)+C,$$

\item[$(*)$] $$\exists\;\ell_0>0\;\forall\;n>0\;\exists\;C\ge 1\;\forall\;x\in\RR^d:\;\;\;\tau^{\ell_0}(x)\le\sigma^{n}(x)+C,$$

\item[$(*)$] $\sigma\hyperlink{omtriangle}{\vartriangleleft}\tau$.
\end{itemize}

	\item[$(e)$] \eqref{weakdifferentgrowth} holds if and only if $\sup_{x\in\RR}\omega(x)=+\infty$ and \eqref{weakdifferentgrowth1} if and only if $\liminf_{|x|\rightarrow+\infty}\omega(x)>0$.
	
	\item[$(f)$] \hyperlink{om3w}{$(\omega_{3,w})$} for $\omega$ in the anisotropic setting precisely means
	\begin{equation}\label{om3wreform}
		\exists\;a>0\;\exists\;b\in\RR\;\forall\;x\in\RR^d:\;\;\;\omega(x)\ge a\log(1+|x|)+b,
	\end{equation}
	and \hyperlink{om3}{$(\omega_3)$} is
	\begin{equation}\label{om3reform}
		\forall\;a\ge 1\;\exists\;b\in\RR\;\forall\;x\in\RR^d:\;\;\;\omega(x)\ge a\log(1+|x|)+b.
	\end{equation}
	We verify that \eqref{strongdifferentgrowthvar}, which is a reformulation of \eqref{strongdifferentgrowth}, is equivalent to \eqref{om3reform}: On the one hand, when $a\ge 1$ is given and \eqref{strongdifferentgrowthvar} holds then apply this condition to $a'>1$ (large) and $\ell>\ell'>0$ (both sufficiently small) such that $\frac{a'}{\ell-\ell'}\ge a$ is valid and which verifies \eqref{om3reform} for this particular $a$. Conversely, if \eqref{om3reform} holds and if $0<\ell'<\ell$ (both small) and $a\ge 1$ are given, then there exists $a'>1$ such that $a'\ge\frac{a}{\ell-\ell'}$ and when applying \eqref{om3reform} to $a'$ property \eqref{strongdifferentgrowthvar} follows with these choices for $a,\ell,\ell'$.
	
	\item[$(g)$] The following are equivalent for radial $\omega$:
	\begin{itemize}
		\item[$(*)$] \eqref{weakom1} for $\mathcal{W}^{\mathfrak{c}}_{\omega}$ is valid.
		
		\item[$(*)$] \eqref{weakom1var} for $\mathcal{W}^{\mathfrak{c}}_{\omega}$ is valid.
		
		\item[$(*)$] $\omega$ satisfies
		\begin{equation}\label{weakom1var1}
			\exists\;C\ge 1\;\forall\;n\in\NN:\;\;\;\omega(n+1)\le C\omega(n)+C.
		\end{equation}
	\end{itemize}
	If $\omega$ is non-decreasing and has \hyperlink{om1}{$(\omega_1)$}, then \eqref{weakom1var1} follows but the latter condition is much weaker: For example $\omega(t):=e^t$ clearly has \eqref{weakom1var1} with $C:=e$ but \hyperlink{om1}{$(\omega_1)$} fails. And, analogously as in Section \ref{alternativecharactsect}, in \eqref{weakom1var1} equivalently we can replace $n$ by any $t\ge 0$ when modifying $C$.
\end{itemize}
Summarizing, the main statements Theorems \ref{maincharthmbeur} and \ref{maincharthmroum} can be applied to exponential-type weight matrices and the characterizing relation of the inclusions of the weighted spaces is given by \hyperlink{ompreceq}{$\preceq$} between the defining weight functions. In the radial setting property \hyperlink{om1}{$(\omega_1)$} is crucial to ensure stability under translation.

\subsection{On the dilatation-type case}\label{specialcasesect2}
We give now the analogous comments for the weight function matrix $\mathcal{W}_{\omega,\mathfrak{c}}$. First note that $\mathcal{W}_{\omega,\mathfrak{c}}$ is formally a weight matrix according to Definition \ref{functionmatrixdef} if $\omega$ is non-decreasing. Moreover, we get:

\begin{itemize}
	\item[$(a)$] $\mathcal{W}_{\omega,\mathfrak{c}}$ is non-decreasing/continuous/radial if and only if $\omega$ is so.
	
	\item[$(b)$] The following are equivalent:
	\begin{itemize}
		\item[$(*)$] $\sup_{x\in\RR^d}\omega(x)<+\infty$,
		
		\item[$(*)$] \eqref{boundedweightequroum} holds,
		
		\item[$(*)$] \eqref{boundedweightequbeur} holds.
	\end{itemize}
	Therefore, $\mathcal{W}_{\omega,{\mathfrak{c}}}$ is unbounded if and only if $\sup_{x\in\RR^d}\omega(x)=+\infty$.
	
	\item[$(c)$] Let $\omega$ be non-decreasing and radial, then $\mathcal{W}_{\omega,\mathfrak{c}}$ satisfies \eqref{mixedom1refomequ1} and\eqref{mixedom1refomequ4} since both conditions amount to have $\omega(2\ell t)\le\omega(nt)+L$ and which is clear: choose $n\ge 2\ell$ when $\ell>0$ is given and $\ell\le\frac{n}{2}$ when $n>0$ is given. (Indeed, $L\ge 1$ is superfluous and can be taken arbitrarily.)
	
	Thus, also \eqref{mixedom1beur} and \eqref{mixedom1roum} follow.
	
	\item[$(d)$] By definition, $\mathcal{W}_{\sigma,\mathfrak{c}}(\preceq)\mathcal{W}_{\tau,\mathfrak{c}}$ and/or $\mathcal{W}_{\sigma,\mathfrak{c}}\{\preceq\}\mathcal{W}_{\tau,\mathfrak{c}}$ amounts to have $\sigma\hyperlink{ompreceqc}{\preceq_{\mathfrak{c}}}\tau$, again this relation yields a uniform scaling of the parameters like in $(d)$ in the previous section, and $\mathcal{W}_{\sigma,\mathfrak{c}}\vartriangleleft\mathcal{W}_{\tau,\mathfrak{c}}$ precisely means $\sigma\hyperlink{omtrianglec}{\vartriangleleft_{\mathfrak{c}}}\tau$. Let $\sigma$, $\tau$ be non-decreasing and radial. Then, by Proposition \ref{weightcomparisonprop} we get that when either $\sigma$ or $\tau$ satisfies in addition \hyperlink{om1}{$(\omega_1)$}, then $\mathcal{W}_{\sigma,\mathfrak{c}}(\preceq)\mathcal{W}_{\tau,\mathfrak{c}}$ and/or $\mathcal{W}_{\sigma,\mathfrak{c}}\{\preceq\}\mathcal{W}_{\tau,\mathfrak{c}}$ implies $\sigma\hyperlink{ompreceq}{\preceq}\tau$, whereas $\sigma\hyperlink{omtriangle}{\vartriangleleft}\tau$ implies $\mathcal{W}_{\sigma,\mathfrak{c}}\vartriangleleft\mathcal{W}_{\tau,\mathfrak{c}}$. And if either $\sigma$ or $\tau$ satisfies in addition \hyperlink{om6}{$(\omega_6)$}, then $\sigma\hyperlink{ompreceq}{\preceq}\tau$ implies $\mathcal{W}_{\sigma,\mathfrak{c}}\{\preceq\}\mathcal{W}_{\tau,\mathfrak{c}}$ and/or $\mathcal{W}_{\sigma,\mathfrak{c}}(\preceq)\mathcal{W}_{\tau,\mathfrak{c}}$, whereas $\mathcal{W}_{\sigma,\mathfrak{c}}\vartriangleleft\mathcal{W}_{\tau,\mathfrak{c}}$ implies $\sigma\hyperlink{omtriangle}{\vartriangleleft}\tau$.\vspace{6pt}
	
\item[$(e)$] \eqref{weakdifferentgrowth} holds if and only if $\sup_{x\in\RR}\omega(ax)-\omega(x)=+\infty$ for some $a>1$ and \eqref{weakdifferentgrowth1} if and only if $\liminf_{|x|\rightarrow+\infty}\omega(ax)-\omega(x)>0$ for some $a>1$.
	
	In general, both requirements are not clear but when $\omega\equiv\omega_{\mathbf{M}}$, $\omega_{\mathbf{M}}$ denoting the \emph{associated weight function,} see \eqref{assofunc}, and with $\mathbf{M}\in\RR_{>0}^{\NN}$ being a \emph{log-convex weight sequence} in the notion of \cite[Def. 2.4]{weightedentireinclusion1}, then both conditions are valid by recalling (the proof of) \cite[Lemma 2.9]{weightedentireinclusion1}. Indeed, in this case even any choice $a>1$ is sufficient to conclude.
	
	\item[$(f)$] $\omega$ satisfies \eqref{om3wreform} if and only if $\omega^{\ell}$ does so for some/any $\ell>0$ and similarly for \eqref{om3reform}. However, the equivalence to \eqref{strongdifferentgrowthvar} established in $(f)$ in the previous Section \ref{specialcasesect1} is not clear in general even if $\omega\equiv\omega_{\mathbf{M}}$ but in this situation the following is known: \eqref{strongdifferentgrowthvar} holds if and only if $\mathbf{M}$ is \emph{derivation closed} which is the classical condition $(M.2)'$ in \cite{Komatsu73}; see \cite[Lemma 2]{nuclearglobal2} even in a more general (mixed) situation.
	
	\item[$(g)$] If $\omega$ is radial, then \eqref{weakom1} resp. \eqref{weakom1var} for $\mathcal{W}_{\omega,\mathfrak{c}}$ amounts to have $\omega(\ell_1(n+1))\le\omega(\ell n)+C$ for some $C\ge 1$ and all $n\in\NN$. Thus, both properties hold automatically when $\omega$ is also assumed to be non-decreasing: In the first estimate choose $\ell_1:=\frac{\ell}{2}$ and in the second one $\ell:=2\ell_1$.
\end{itemize}

Summarizing, the main statements Theorems \ref{maincharthmbeur} and \ref{maincharthmroum} can be applied to dilatation-type weight matrices; here the inclusion is characterized in terms of relation \hyperlink{ompreceqc}{$\preceq_{\mathfrak{c}}$} between the defining weight functions. In view of comment $(d)$ under additional assumptions on the weights a connection to the corresponding exponential-type weight system and the (characterizing) relation \hyperlink{ompreceq}{$\preceq$} can be established. And by comment $(c)$ it holds that the associated weighted classes are invariant under translation for any non-decreasing and radial $\omega$.

\section{Weighted (test-)function classes of Beurling-Bj\"{o}rck-type}\label{BBclasssect}
We study related and more general cases and comment how and if the results from Sections \ref{generalclasssect} and \ref{charactinclsection} can be transferred.

\subsection{Weighted subspaces of a fixed Banach space}\label{subspacesect}
Let $\mathcal{W}=\{\omega^{\ell}: \ell>0\}$ be a weight function matrix and let $\mathcal{A}^{p}\subseteq L^p$ be a fixed Banach space, $1\le p\le\infty$. Then define $$\mathcal{A}^{\infty}_{\omega^{\ell}}:=\left\{f\in\mathcal{A}^{\infty}:\;\;\;\|f\|_{\infty,\omega^{\ell}}<+\infty\right\},$$
and
$$\mathcal{A}^{\infty}_{\{\mathcal{W}\}}:=\bigcup_{\ell>0}\mathcal{A}^{\infty}_{\omega^{\ell}},\hspace{15pt}\mathcal{A}^{\infty}_{(\mathcal{W})}:=\bigcap_{\ell>0}\mathcal{A}^{\infty}_{\omega^{\ell}},$$
endowed with their natural topologies. Moreover, if $\mathcal{W}$ is \emph{continuous} and $1\le p<\infty$, then for any parameter $\ell>0$ let us set
$$\mathcal{A}^p_{\omega^{\ell}}:=\left\{f\in\mathcal{A}^p:\;\;\;\|f\|_{p,\omega^{\ell}}<+\infty\right\},$$
and put accordingly
$$\mathcal{A}^p_{\{\mathcal{W}\}}:=\bigcup_{\ell>0}\mathcal{A}^p_{\omega^{\ell}},\hspace{15pt}\mathcal{A}^p_{(\mathcal{W})}:=\bigcap_{\ell>0}\mathcal{A}^p_{\omega^{\ell}}.$$
Again, $\mathcal{A}^p_{\{\mathcal{W}\}}$ and $\mathcal{A}^p_{(\mathcal{W})}$ are endowed with their natural topologies and the following continuous inclusions hold:
$$\mathcal{A}^p_{(\mathcal{W})}\subseteq\mathcal{A}^p_{\{\mathcal{W}\}}\subseteq\mathcal{A}^p,\;\;\;1\le p\le\infty.$$
Thus, the special case $\mathcal{A}^{p}=L^p$ corresponds to the situation studied in Sections \ref{generalclasssect} and \ref{charactinclsection} and obviously, by definition,
\begin{equation}\label{fixedBequ}
\mathcal{A}^{p}_{[\mathcal{W}]}=\mathcal{L}^{p}_{[\mathcal{W}]}\cap\mathcal{A}^p\subseteq\mathcal{L}^{p}_{[\mathcal{W}]},\;\;\;1\le p\le\infty.
\end{equation}
We comment in detail how the shown statements transfer if the symbol/functor $\mathcal{L}$ is replaced by $\mathcal{A}$. A natural basic requirement, to avoid trivial situations and which is assumed from now on, is $\mathcal{A}^p\neq\{0\}$; i.e. the fixed space is \emph{non-trivial.} Otherwise, of course $\mathcal{A}^p_{[\mathcal{W}]}=\{0\}$ holds.

\begin{itemize}
\item[$(a)$] Section \ref{unboundedsect}; Theorem \ref{unboundednesscor}: By definition, \eqref{boundedweightequroum} implies $\mathcal{A}^p_{\{\mathcal{W}\}}=\mathcal{A}^p$ (as sets and l.c.v.s.) and \eqref{boundedweightequbeur} yields the same equalities for the Beurling-types.

    Concerning the converse, the proof of $(i)$ in Proposition \ref{thirdweightcomppropzero} immediately transfers when we assume that $\mathcal{A}^{\infty}$ \emph{contains at least all (equivalence classes of) constant functions} and, indeed, it suffices to assume
    \begin{equation}\label{tildefunction}
    \exists\;(0\neq)\widetilde{f}\in\mathcal{C}(\RR^d):\;\;\;\widetilde{f}\in\mathcal{A}^{\infty},\;\;\;\inf_{x\in\RR^d}|\widetilde{f}(x)|>\epsilon>0.
    \end{equation}
    Note that even in the proof of $(i)$ in Proposition \ref{thirdweightcomppropzero} it is not sufficient to take any arbitrary (non-trivial) $f\in\mathcal{C}(\RR^d)\cap L^{\infty}$ since such $f$ can vanish identically near infinity. To transfer the proof of $(ii)$ we assume that $\mathcal{A}^p$ \emph{contains at least all (equivalence classes of) locally constant functions $f$ satisfying $\|f\|_{L^p}<+\infty$.} Thus $g^p\in\mathcal{A}^p$, $g^p$ the function defined in \eqref{thirdweightcomppropzeroequ1}. Therefore, in this case the characterizing Theorem \ref{unboundednesscor} is valid, too.

\item[$(b)$] Section \ref{nontrivialsect}: The \emph{non-triviality} of $\mathcal{A}^p_{[\mathcal{W}]}$ is \emph{not clear in general} even if $\mathcal{A}^p\neq\{0\}$ because the technical constructive proof of Lemma \ref{nontriviallemma} cannot be transferred automatically.

    In view of \eqref{fixedBequ} the non-triviality amounts to having
    $$\{0\}\neq\mathcal{L}^{p}_{[\mathcal{W}]}\cap\mathcal{A}^p.$$

\item[$(c)$] Section \ref{stabtranslationsect}: The proof of Proposition \ref{translationlemma} purely involves weights and their growth properties and transfers immediately. However, when $\mathcal{A}^p$ is \emph{not translation invariant,} then it is not clear in general that for given $f\in\mathcal{A}^p$ we get $T_{x_0}(f)\in\mathcal{A}^p$ and hence by definition the estimates in Proposition \ref{translationlemma} do not yield the fact that (each) $T_{x_0}$ is acting continuously between weighted spaces $\mathcal{A}^{p}_{[\mathcal{T}]}$ and $\mathcal{A}^{p}_{[\mathcal{S}]}$. Therefore, it is natural to assume from the beginning that $\mathcal{A}^p$ is \emph{translation invariant} and then this result immediately transfers to the weighted $\mathcal{A}^p$-setting.

\item[$(d)$] Section \ref{generacharactsect}; main Theorems \ref{maincharthmbeur} and \ref{maincharthmroum}: Proposition  \ref{firstweightcompprop} is clear by the analogous system of seminorms. And the proofs of the crucial Propositions \ref{Beurlingcompprop} and \ref{Roumcompprop} transfer when $\mathcal{A}^p$ is \emph{translation invariant} by taking into account the previous comment $(c)$ and this requirement is crucial in order to proceed. Indeed, Proposition \ref{translationlemma} transfers and shows that under \eqref{maincharthmbeurequ} resp. \eqref{maincharthmroumequ} the corresponding weighted class $\mathcal{A}^{p}_{(\mathcal{W})}$ resp. $\mathcal{A}^{p}_{\{\mathcal{W}\}}$ is translation invariant. Note (again) that the proofs of Propositions \ref{Beurlingcompprop} and \ref{Roumcompprop} require also the non-triviality of the corresponding (smaller) weighted classes and which is not clear by comment $(b)$ above. However, when taking into account the conventions mentioned in Remark \ref{infinityremark1} one is able to conclude in any case.

\item[$(e)$] Section \ref{alternativecharactsect}: Similarly as in comment $(b)$ before in general the methods involving the optimal functions $\theta^p_{\omega^{\ell}}$ from \eqref{charactfctdef} and $\theta^{\infty}_{\omega^{\ell}}$ from \eqref{charactfctdef1} fail since neither $\theta^p_{\omega^{\ell}}\in\mathcal{A}^{p}_{\{\mathcal{W}\}}$ nor $\theta^{\infty}_{\omega^{\ell}}\in\mathcal{A}^{\infty}_{\{\mathcal{W}\}}$ is clear. If these functions belong to the spaces, then all shown results transfer immediately.

    On the other hand, when one assumes $\widetilde{\theta}^{\infty}_{\omega^{\ell}}:=\widetilde{f}\cdot\theta^{\infty}_{\omega^{\ell}}\in\mathcal{A}^{\infty}_{\{\mathcal{W}\}}$, $\widetilde{f}\in\mathcal{C}(\RR^d)$ a function satisfying the infimum-condition from \eqref{tildefunction}, then the crucial estimate in the proof of Proposition \ref{thirdweightcompprop} reads as follows:
$$\forall\;\ell>0\;\exists\;\ell'>0\;\exists\;D\ge 0\;\forall\;x\in\RR^d:\;\;\;\tau^{\ell'}(x)\le\sigma^{\ell}(x)-\log(|\widetilde{f}(x)|)+D\le\sigma^{\ell}(x)-\log(\epsilon)+D,$$
which yields $\mathcal{S}\{\preceq\}\mathcal{T}$. Note that $-\log(\epsilon)>0$ because w.l.o.g. $\epsilon\in(0,1)$. Similarly Proposition \ref{thirdstrongweightcompprop} holds, too.

And assume that there exists a continuous function $(0\neq)\widetilde{f}\in\mathcal{C}(\RR^d)$ with $\inf_{x\in\RR^d}|\widetilde{f}(x)|>\epsilon^{1/p}$ for some $\epsilon>0$ and such that $\widetilde{\theta}^p_{\omega^{\ell}}:=\widetilde{f}\cdot\theta^p_{\omega^{\ell}}\in\mathcal{A}^{p}_{\{\mathcal{W}\}}$ for all $\ell>0$. Then also Propositions \ref{secondweightcompprop} and \ref{secondstrongweightcompprop} are valid (under the same assumptions on $\mathcal{S}$ and $\mathcal{T}$): First we have that \eqref{secondweightcomppropequ} turns into
\begin{equation}\label{secondweightcomppropequnew}
\forall\;\ell>0\;\exists\;\ell_1>0:\;\;\;\int_{\RR^d}|\widetilde{\theta}^p_{\sigma^{\ell}}|^pdx=\int_{\RR^d}|\widetilde{f}(x)|^pe^{\tau^{\ell_1}(x)-\sigma^{\ell}(x)}dx<+\infty.
\end{equation}
Second, the crucial estimate in the proof of Proposition \ref{secondweightcompprop} changes as follows:
\begin{align*}
&\int_{\RR^d}|\widetilde{f}(x)|^pe^{\tau^{\ell_1}(x)-\sigma^{\ell}(x)}dx=\sum_{n=0}^{+\infty}\int_{x\in\RR^d: n\le|x|<n+1}|\widetilde{f}(x)|^pe^{\tau^{\ell_1}(x)-\sigma^{\ell}(x)}dx
\\&
\ge\sum_{n=0}^{+\infty}e^{\tau^{\ell_1}(n)-\sigma^{\ell}(n+1)}\int_{x\in\RR^d: n\le|x|<n+1}|\widetilde{f}(x)|^pdx
\\&
\ge\sum_{n=0}^{+\infty}e^{\tau^{\ell_1}(n)-\sigma^{\ell}(n+1)}\inf_{x\in\RR^d}|\widetilde{f}(x)|^p\int_{x\in\RR^d: n\le|x|<n+1}1dx
\\&
=\sum_{n=0}^{+\infty}e^{\tau^{\ell_1}(n)-\sigma^{\ell}(n+1)}\inf_{x\in\RR^d}|\widetilde{f}(x)|^p\frac{\pi^{d/2}((n+1)^d-n^d)}{\Gamma(1+d/2)}\ge\sum_{n=0}^{+\infty}e^{\tau^{\ell_1}(n)-\sigma^{\ell}(n+1)}\frac{\pi^{d/2}\epsilon}{\Gamma(1+d/2)}.
\end{align*}
For the last estimate recall that $(n+1)^d-n^d\ge 1$ for all $n\in\NN$ and $d\in\NN_{>0}$.

Combining this estimate with \eqref{secondweightcomppropequnew} it follows again that $\sup_{n\in\NN}e^{\tau^{\ell_1}(n)-\sigma^{\ell}(n+1)}<+\infty$; i.e. \eqref{secondweightcomppropequ1} is verified and the rest follows as in the proof above.
\end{itemize}

\subsection{Weighted Fourier image spaces}\label{weightedFouriersection}
In this section we comment on weighted spaces which are defined analogously as the ones in Section \ref{generalclasssect} but the crucial difference is that we are weighting the Fourier transform $\widehat{f}$ of the functions under consideration:

For a fixed function space $\mathcal{A}$, a (continuous) weight function matrix $\mathcal{W}$ and $1\le p\le\infty$ we set $$\widehat{\mathcal{A}}^{p}_{(\mathcal{W})}:=\{f\in\mathcal{A}:\;\forall\;\ell>0:\; \|\widehat{f}\|_{p,\omega^{\ell}}<+\infty\},\;\;\;\widehat{\mathcal{A}}^{p}_{\{\mathcal{W}\}}:=\{f\in\mathcal{A}:\;\exists\;\ell>0:\; \|\widehat{f}\|_{p,\omega^{\ell}}<+\infty\},$$
and both spaces are endowed with their natural topologies. Of course, in order to make this definition useful one requires that the Fourier-transform is well-defined on $\mathcal{A}$ and that $\widehat{\mathcal{A}}^{p}_{[\mathcal{W}]}$ is non-trivial. We remark:

\begin{itemize}
\item[$(a)$] Proposition \ref{firstweightcompprop} is clear for this setting by the analogous definition of the seminorms but with $f$ being replaced by $\widehat{f}$.

\item[$(b)$] The non-triviality of $\widehat{\mathcal{A}}^p_{[\mathcal{W}]}$ is not guaranteed in general and has to be checked. This step requires technical investigations, in particular, when explicitly constructing (optimal) functions belonging to such classes.

\item[$(c)$] Assume that $\{0\}\neq\widehat{\mathcal{A}}^p_{[\mathcal{S}]}$ holds and that the underlying space $\mathcal{A}$ is \emph{modulation invariant;} i.e. if $f\in\mathcal{A}$ then also $g_{x_0}:=e^{i\langle\cdot,x_0\rangle}f\in\mathcal{A}$ for any $x_0\in\RR^d$. Then the proofs of Propositions \ref{Beurlingcompprop} and \ref{Roumcompprop} can be transferred as follows:

Take a (fixed) non-trivial $f\in\widehat{\mathcal{A}}^p_{[\mathcal{S}]}$ and for any $x_0\in\RR^d$ put $g_{x_0}(x):=e^{i\langle x,x_0\rangle}f(x)$. Then $\widehat{g}_{x_0}(x)=\widehat{f}(x-x_0)=T_{x_0}(\widehat{f})$ and $g_{x_0}\in\mathcal{A}$ for all $x_0\in\RR^d$ by modulation invariance of $\mathcal{A}$. We apply \eqref{Beurlingcomppropequ}, \eqref{fourthweightcomppropequ1}, \eqref{fourthweightcomppropequroum}, and \eqref{fourthweightcomppropequ1roum} to the family $\{g_{x_0}: x_0\in\RR^d\}$; see also the original proof of \cite[Thm. 1.3.18]{Bjorck66} (for the Beurling-type). Indeed, Proposition \ref{translationlemma} transfers (with $f$ being replaced by $\widehat{f}$) and shows that under \eqref{maincharthmbeurequ} resp. \eqref{maincharthmroumequ} the corresponding weighted class $\widehat{\mathcal{A}}^{p}_{(\mathcal{W})}$ resp. $\widehat{\mathcal{A}}^{p}_{\{\mathcal{W}\}}$ is \emph{modulation invariant.}

However, as mentioned in Remark \ref{infinityremark1} and in comment $(d)$ in Section \ref{subspacesect}, the same implication holds when $\{0\}=\widehat{\mathcal{A}}^p_{[\mathcal{S}]}$ (and $\{0\}=\widehat{\mathcal{A}}^p_{[\mathcal{T}]}$).

Summarizing, the analogues of the main results Theorems \ref{maincharthmbeur} and \ref{maincharthmroum} hold for this described ``Fourier-setting'' too when assuming the same conditions on the matrices $\mathcal{S}$ and $\mathcal{T}$.
\end{itemize}

\subsection{Test function spaces appearing in the literature}\label{examplesect}
We revisit now the different types and notions of prominent weighted (test-)function spaces which have been considered in the literature. All known examples are expressed in terms of radial/isotropic \emph{exponential-type weight function matrices} and by the Fourier framework described in the previous Section \ref{weightedFouriersection} when taking spaces $\mathcal{A}\subseteq L^1$. Finally, $(A)$ and $(B)$ correspond to $p=1$ and $(C)$ to $p=\infty$.

\begin{itemize}
	\item[$(A)$] Let $\omega$ be a \emph{non-quasianalytic} BMT-weight except necessarily satisfying the convexity condition \hyperlink{om4}{$(\omega_4)$}. Let $U\subseteq\RR^d$ be non-empty open and then set
	$$\mathcal{D}_{\{\omega\}}(U):=\{f\in\mathcal{D}(U),\;\,\exists\;\ell>0:\;\int_{\RR^d}|\widehat{f}(x)|e^{\ell\omega(x)}dx<+\infty\},$$
	and
	$$\mathcal{D}_{(\omega)}(U):=\{f\in\mathcal{D}(U),\;\;\forall\;\ell>0:\;\int_{\RR^d}|\widehat{f}(x)|e^{\ell\omega(x)}dx<+\infty\}.$$
	For this setting see \cite[Def. 3.1]{BraunMeiseTaylor90} and in this work, as described in the introduction, with the help of \hyperlink{om4}{$(\omega_4)$} those spaces have been transferred into the nowadays frequently used definition involving $\varphi^{*}_{\omega}$; they are usually denoted by $\mathcal{E}_{\{\omega\}}(U)$, $\mathcal{E}_{(\omega)}(U)$ (also for quasianalytic weights).
	
	\item[$(B)$] In \cite[Def. 1.3.1, 1.3.2 \& 1.3.4]{Bjorck66} for any $\omega:[0,+\infty)\rightarrow[0,+\infty)$ which satisfies \hyperlink{sub}{$(\omega_{\on{sub}})$}, i.e. $(\alpha)$ in \cite{Bjorck66}, the following (global) test function classes (of Beurling-type) have been introduced:
	$$(\mathfrak{D}_{(\omega)}(U)=)\mathfrak{D}_{\omega}(U):=\{f\in L^1(\RR^d): \exists\;K\subset\subset U,\;\supp(f)\subseteq K,\;\forall\;\ell>0:\;\int_{\RR^d}|\widehat{f}(x)|e^{\ell\omega(x)}dx<+\infty\}.$$
	The non-triviality of this class has been characterized in \cite[Thm. 1.3.7]{Bjorck66} in terms of the non-quasianalyticity condition for weight functions.
	
\item[$(C)$] Finally, in \cite{PetzscheVogt} for $\omega$ satisfying the conditions listed in Definition \ref{defPVweightfct} (i.e. PV-weights) the following test function classes have been studied for $U\subseteq\RR^d$ (see \cite[p. 18-19]{PetzscheVogt}):
	$$\widetilde{\mathcal{D}}_{\{\omega\}}(U):=\{f\in L^1(\RR^d): \exists\;K\subset\subset U,\;\supp(f)\subseteq K,\;\exists\;\ell>0:\;\sup_{x\in\RR^d}|\widehat{f}(x)|e^{\ell\omega(x)}<+\infty\},$$
	and
	$$\widetilde{\mathcal{D}}_{(\omega)}(U):=\{f\in L^1(\RR^d): \exists\;K\subset\subset U,\;\supp(f)\subseteq K,\;\forall\;\ell>0:\;\sup_{x\in\RR^d}|\widehat{f}(x)|e^{\ell\omega(x)}<+\infty\}.$$
\end{itemize}

Consequently, the established results can be applied to all these examples and since we deal with exponential-type weight matrices relation \hyperlink{ompreceq}{$\preceq$} between the defining weights is characterizing the inclusion(s) of the above weighted spaces and equivalence (as l.c.v.s.) is characterized in terms of \emph{equivalent weight functions,} i.e. relation \hyperlink{sim}{$\sim$} is relevant; see the comments in Sections \ref{specialcasesect1} and \ref{weightedFouriersection}.

Some parts of these characterizations have already been stated in the literature, for $(B)$ we refer to \cite[Thm. 1.3.18]{Bjorck66} concerning the Beurling-type, for $(A)$ in \cite[Rem. 3.2 $(2)$]{BraunMeiseTaylor90} it has been mentioned that U. Franken has also treated the Beurling-type in his Master thesis \cite{diplomafranken}. Note that \emph{non-quasianalyticity \hyperlink{omnq}{$(\omega_{\on{nq}})$}} for the defining weight $\omega$ is crucial in each setting in order to ensure \emph{non-triviality.}

\begin{remark}\label{truncationrem}
Indeed, in \cite{PetzscheVogt} the authors have even considered the non-standard situation when the (matrix) parameter $\ell$ is restricted to $(0,\Lambda)$ for $\widetilde{\mathcal{D}}_{(\omega)}(U)$ and to $(\Lambda,+\infty)$ for $\widetilde{\mathcal{D}}_{\{\omega\}}(U)$. Here, $\Lambda\in(0+\infty)$ denotes a \emph{``truncation parameter''.}

The established techniques transfer immediately with the obvious changes in the proofs of Proposition \ref{translationlemma}, Lemma \ref{mixedom1refom}, Theorems \ref{maincharthmbeur} and \ref{maincharthmroum} when restricting in the particular conditions the indices w.r.t. given truncation parameters $\Lambda'$ and $\Lambda$ for the matrices $\mathcal{T}$ and $\mathcal{S}$, respectively, and so also the above comments in Section \ref{BBclasssect} apply. More precisely, the characterizing growth relation $\mathcal{S}(\preceq)\mathcal{T}$ and hence \eqref{growthrel1} turns into
\begin{equation}\label{growthrel1trunc}
\forall\;\ell\in(0,\Lambda')\;\exists\;n\in(0,\Lambda)\;\exists\;C\ge 1\;\forall\;x\in\RR^d:\;\;\;\tau^{\ell}(x)\le\sigma^n(x)+C,
\end{equation}
whereas \eqref{growthrel2} turns into
\begin{equation}\label{growthrel2trunc}
\forall\;n\in(\Lambda,+\infty)\;\exists\;\ell\in(\Lambda',+\infty)\;\exists\;C\ge 1\;\forall\;x\in\RR^d:\;\;\;\tau^{\ell}(x)\le\sigma^n(x)+C.
\end{equation}
Now focus on the exponential-type system and first note that the equivalences stated in $(c)$ and $(d)$ in Section \ref{specialcasesect1} are not valid anymore in this situation: Recall that \hyperlink{om1}{$(\omega_1)$} for a non-decreasing radial $\omega$ means $\omega(2t)\le L\omega(t)+L$ for some $L\ge 1$ and all $t\ge 0$. This constant $L$ yields then a uniform scaling in the parameters of the matrix $\mathcal{W}^{\mathfrak{c}}_{\omega}$; i.e. one infers the relation $L\ell=n$ and which is a problem when restricting the parameters and when $\ell\rightarrow\Lambda'$ (resp. $n\rightarrow\Lambda$). Indeed, $L$ yields a transformation (dilation) of the truncation parameter.

Moreover, $\sigma\hyperlink{ompreceq}{\preceq}\tau$ means $\tau(x)\le C\sigma(x)+C$ for some $C\ge 1$ and all $x\in\RR^d$ and so \eqref{growthrel1trunc} and \eqref{growthrel2trunc} imply this relation as well (with no restrictions on $\Lambda'$, $\Lambda$). However, for the converse the constant $C$ is then required to be related to $\Lambda,\Lambda'$; more precisely for fixed $\Lambda,\Lambda'\in(0,+\infty)$ consider the following assertions:
\begin{itemize}
\item[$(i)$] One has
\begin{equation}\label{growthrel5trunc}
\exists\;C\ge 1\;\forall\;x\in\RR^d:\;\;\;\tau(x)\le\frac{\Lambda}{\Lambda'}\sigma(x)+C.
\end{equation}

\item[$(ii)$] \eqref{growthrel1trunc} holds.

\item[$(iii)$] \eqref{growthrel2trunc} holds.

\item[$(iv)$] One has
\begin{equation}\label{growthrel6trunc}
\forall\;0<K'<\Lambda'\;\forall\;K>\Lambda\;\exists\;D\ge 1\;\forall\;x\in\RR^d:\;\;\;\tau(x)\le\frac{K}{K'}\sigma(x)+D.
\end{equation}
\end{itemize}
Then the implications $(i)\Rightarrow(ii),(iii)\Rightarrow(iv)$ are valid:

If \eqref{growthrel5trunc} holds, then \eqref{growthrel1trunc} follows by taking $n:=\frac{\Lambda}{\Lambda'}\ell\in(0,\Lambda)$ for any given $\ell\in(0,\Lambda')$ and \eqref{growthrel2trunc} follows when taking $\ell:=n\frac{\Lambda'}{\Lambda}\in(\Lambda',+\infty)$ for any given $n\in(\Lambda,+\infty)$.

$(ii),(iii)\Rightarrow(iv)$: Fix $K'<\Lambda'$ and $K>\Lambda$ and then by \eqref{growthrel1trunc} one has with $\ell:=\Lambda'-\epsilon_1$ and $n:=\Lambda-\epsilon_2$, both $\epsilon_i$ chosen sufficiently small depending on given $K'$, that
$$\exists\;\epsilon_1,\epsilon_2>0:\;\;\;\frac{\Lambda-\epsilon_2}{\Lambda'-\epsilon_1}\le\frac{K}{K'}.$$
Note that by the point-wise order of the weight functions in the matrices in \eqref{growthrel1trunc} we are interested in $\ell\rightarrow\Lambda'$ and one can assume w.l.o.g. $n=n(\ell)\rightarrow\Lambda$. Similarly, \eqref{growthrel2trunc} implies with $n:=\Lambda+\epsilon_1$ and $\ell:=\Lambda'+\epsilon_2$, again both $\epsilon_i$ small enough depending on given $K$, that
$$\exists\;\epsilon_1,\epsilon_2>0:\;\;\;\frac{\Lambda+\epsilon_1}{\Lambda'+\epsilon_2}\le\frac{K}{K'}.$$
Thus \eqref{growthrel6trunc} holds in both cases.

Finally, \eqref{thirdstrongweightcomppropequ} and so $\mathcal{S}\vartriangleleft\mathcal{T}$ transfers in this setting to
\begin{equation}\label{thirdstrongweightcomppropequtrunc1}
\forall\;\ell\in(0,\Lambda')\;\forall\;n\in(\Lambda,+\infty)\;\exists\;D\ge 1\;\forall\;x\in\RR^d:\;\;\;\tau^{\ell}(x)\le\sigma^{n}(x)+D.
\end{equation}
Indeed \eqref{growthrel6trunc} and \eqref{thirdstrongweightcomppropequtrunc1} are equivalent: Let $\epsilon,\epsilon_1>0$, then when applying \eqref{growthrel6trunc} to $K:=\Lambda+\epsilon$ and $K':=\Lambda'-\epsilon_1$ we get \eqref{thirdstrongweightcomppropequtrunc1} with $\ell=K'$ and $n=K$. Conversely, we apply \eqref{thirdstrongweightcomppropequtrunc1} to $\ell=\Lambda'-\epsilon_1$ and $n=\Lambda-\epsilon$ and get \eqref{growthrel6trunc} with $K'=\ell$ and $K=n$. Since $\epsilon,\epsilon_1$ can be chosen arbitrarily small, the equivalence is established.

This also illustrates the difference compared with the standard setting when not involving a truncation parameter since in this case $\mathcal{S}\vartriangleleft\mathcal{T}$ obviously implies both $\mathcal{S}(\preceq)\mathcal{T}$ and $\mathcal{S}\{\preceq\}\mathcal{T}$ and hence the converse.
\end{remark}

\section{On the failure of the convexity condition $(\omega_4)$}\label{failuresection}
In this final section we provide the technical construction of a function $\omega:[0,+\infty)\rightarrow[0,+\infty)$ satisfying the following properties:

\begin{itemize}
	\item[$(I)$] $\omega$ is normalized, continuous, non-decreasing and $\lim_{t\rightarrow+\infty}\omega(t)=+\infty$.
	
	\item[$(II$)] $\omega$ satisfies \hyperlink{om3}{$(\omega_3)$}.
	
	\item[$(III)$] $\gamma(\omega)=+\infty$ holds. This implies that $\omega$ also has \hyperlink{omsnq}{$(\omega_{\on{snq}})$},  \hyperlink{om1}{$(\omega_1)$}, and that $\omega$ is equivalent to the continuous concave resp. subadditive function $\kappa_{\omega}$; see Lemma \ref{chaimimplemma} and Remark \ref{secondcomprem}. In fact, we extend \eqref{chaimofimpl} and even show that $\omega$ is \emph{slowly varying,} see \cite[$(1.2.1)$]{regularvariation} and \eqref{stepVIequ1}:
$$\forall\;u>0:\;\;\;\lim_{t\rightarrow+\infty}\frac{\omega(tu)}{\omega(t)}=1.$$

	\item[$(IV)$] There does not exist any function $\sigma:[0,+\infty)\rightarrow[0,+\infty)$ such that $\sigma$ and $\omega$ are equivalent and such that $\varphi_{\sigma}$ is convex. Thus, in particular, $\varphi_{\omega}$ is not convex. In view of $(I)$ we have that equivalence is described by $O$-growth relations as $t\rightarrow+\infty$; see Section \ref{growthrelsection}.

\item[$(V)$] There exist (uncountable) infinitely many pairwise non-equivalent functions satisfying $(I)-(IV)$.
\end{itemize}

These properties yield:

\begin{itemize}
\item[$(\mathfrak{C}_1)$] By $(I)-(III)$, and in view of Remark \ref{secondcomprem} and Lemmas \ref{alpha0equivalence} and \ref{chaimimplemma}, we have that $\omega$ is equivalent to a weight function in the sense of Beurling-Bj\"{o}rck (see Definition \ref{defBBweightfct}) and of Petzsche-Vogt (see Definition \ref{defPVweightfct}).

\item[$(\mathfrak{C}_2)$] By $(IV)$ it follows that in the whole equivalence class of $\omega$ we cannot find any BMT-weight function. In particular, $\omega$ itself is not a BMT-weight function and, indeed, even any arbitrary function being equivalent to $\omega$ violates \hyperlink{om4}{$(\omega_4)$}.

\item[$(\mathfrak{C}_3)$] Finally, $(V)$ implies that there exist ``quite many'' weight functions having this non-standard behavior.
\end{itemize}

When combining the information from $(\mathfrak{C}_1)$ and $(\mathfrak{C}_2)$ with the characterizations obtained in Theorems \ref{maincharthmbeur} and \ref{maincharthmroum} (see also Proposition \ref{thirdweightcompprop}), and by taking into account the comments in Section \ref{specialcasesect1}, especially the characterization of the growth conditions and relations in comments $(c)$ and $(d)$ there, within the exponential-type weight function matrix setting we obtain: There exist examples of weighted classes treated in the Beurling-Bj\"{o}rck-setting which cannot be described (alternatively) by any BMT-weight function and this statement applies to both Roumieu- and Beurling-type spaces.\vspace{6pt}

This (counter-)example should be compared with Theorem \ref{Frankenthm} and it is somehow converse to this result since we are constructing explicitly a BB-weight and not a BMT-weight function.

Finally, note that $\gamma(\omega)=+\infty$ implies $\overline{\gamma}(\omega)=+\infty$ and thus \hyperlink{om6}{$(\omega_6)$} fails (recall Section \ref{growthindexsection}).

\subsection{Preparatory results}\label{counterexprepsection}
We collect several preparatory results which are needed in the proof of the construction in the next section.

\begin{lemma}\label{counterlemma0}
Let $\omega:[0,+\infty)\rightarrow[0,+\infty)$ be arbitrary. Assume that there exists $\sigma:[0,+\infty)\rightarrow[0,+\infty)$ such that $\omega$ and $\sigma$ are equivalent and such that $\varphi_{\sigma}: t\mapsto\sigma(e^t)$ is convex (on $\RR$), i.e. $\sigma$ has \hyperlink{om4}{$(\omega_4)$}. Then $\omega$ has to satisfy:
\begin{equation}\label{nonconvexequ}
\exists\;A\ge 1\;\forall\;x,y\in\RR\;\forall\;0\le t\le 1:\;\;\;\varphi_{\omega}(ty+(1-t)x)\le At\varphi_{\omega}(y)+A(1-t)\varphi_{\omega}(x)+A.
\end{equation}
\end{lemma}

\demo{Proof}
The equivalence between $\omega$ and $\sigma$ gives that
$$\exists\;B\ge 1\;\forall\;t\ge 0:\;\;\;-1+B^{-1}\sigma(t)\le\omega(t)\le B\sigma(t)+B,$$
and by taking $s:=e^t$ these estimates transfer to the corresponding functions $\varphi_{\omega}$, $\varphi_{\sigma}$ as well (with the same choice for $B$ and for all $t\in\RR$).

Thus we can estimate as follows for all $x,y\in\RR$ and $0\le t\le 1$:
\begin{align*}
&\varphi_{\omega}(ty+(1-t)x)\le B\varphi_{\sigma}(ty+(1-t)x)+B\le Bt\varphi_{\sigma}(y)+B(1-t)\varphi_{\sigma}(x)+B
\\&
\le B^2t\varphi_{\omega}(y)+B^2t+B^2(1-t)\varphi_{\omega}(x)+B^2(1-t)+B=B^2t\varphi_{\omega}(y)+B^2(1-t)\varphi_{\omega}(x)+B^2+B.
\end{align*}
This verifies \eqref{nonconvexequ} for the choice $A:=B^2+B$.
\qed\enddemo

\begin{lemma}\label{counterlemma1}
Assume that $\omega:[0,+\infty)\rightarrow[0,+\infty)$ satisfies
\begin{equation}\label{counterlemma1equ1}
\liminf_{t\rightarrow+\infty}\omega(t)>1\Leftrightarrow\liminf_{t\rightarrow+\infty}\varphi_{\omega}(t)>1,
\end{equation}
and such that
\begin{equation}\label{counterlemma1equ2}
\forall\;j\in\NN_{>0}\;\exists\;0<t_j<\frac{1}{j}\;\exists\;a_j,b_j\ge 0:\;\;\;\varphi_{\omega}(t_jb_j+(1-t_j)a_j)=\varphi_{\omega}(b_j)=\frac{2j\varphi_{\omega}(a_j)}{1-jt_j}.
\end{equation}
Then there does not exist a function $\sigma:[0,+\infty)\rightarrow[0,+\infty)$ such that $\varphi_{\sigma}$ is convex and such that $\omega$ and $\sigma$ are equivalent.
\end{lemma}

Note that \eqref{counterlemma1equ1} is \eqref{liminfweighfct} and a mild requirement; it is clearly satisfied if $\lim_{t\rightarrow+\infty}\omega(t)=+\infty$.

\demo{Proof}
First, we claim that
$$\exists\;j_0\in\NN_{>0}\;\forall\;j\ge j_0:\;\;\;\varphi_{\omega}(b_j)>jt_j\varphi_{\omega}(b_j)+j(1-t_j)\varphi_{\omega}(a_j)+j.$$
For this note that
\begin{align*}
&\varphi_{\omega}(b_j)>jt_j\varphi_{\omega}(b_j)+j(1-t_j)\varphi_{\omega}(a_j)+j\Leftrightarrow\varphi_{\omega}(b_j)(1-jt_j)>j(1-t_j)\varphi_{\omega}(a_j)+j
\\&
\Leftrightarrow 2j\varphi_{\omega}(a_j)>j(1-t_j)\varphi_{\omega}(a_j)+j;
\end{align*}
we have used the second equality in \eqref{counterlemma1equ2} and the last estimate is valid since $1-t_j<1$ for all $j\in\NN_{>0}$ and since $\varphi_{\omega}(a_j)>1$ for all $j$ sufficiently large which holds by assumption \eqref{counterlemma1equ1}. By taking into account the first equality in \eqref{counterlemma1equ2} we have shown:
\begin{align*}
&\forall\;C\ge 1\;\exists\;j_C\in\NN_{>0}\;\forall\;j\ge j_C\;\exists\;0<t_j<\frac{1}{j}\;\exists\;a_j,b_j\ge 0:
\\&
\varphi_{\omega}(t_jb_j+(1-t_j)a_j)=\varphi_{\omega}(b_j)>Ct_j\varphi_{\omega}(b_j)+C(1-t_j)\varphi_{\omega}(a_j)+C.
\end{align*}
This violates \eqref{nonconvexequ} and finishes the proof by Lemma \ref{counterlemma1}.
\qed\enddemo

\begin{lemma}\label{counterlemma2}
Let $\omega:[0,+\infty)\rightarrow[0,+\infty)$ be non-decreasing and such that $\lim_{t\rightarrow+\infty}\omega(t)=+\infty$. Then $\omega$ satisfies $(P_{\omega,\gamma})$ if and only if
$$\exists\;K>1:\;\;\;\limsup_{t\rightarrow+\infty}\frac{\varphi_{\omega}(\gamma\log(K)+t)}{\varphi_{\omega}(t)}<K.$$
\end{lemma}

\demo{Proof}
We set $s:=e^t$ and the definition of $\varphi_{\omega}(t)=\omega(e^t)$ yields $\frac{\omega(K^{\gamma}s)}{\omega(s)}=\frac{\omega(e^{\gamma\log(K)}e^t)}{\omega(e^t)}=\frac{\varphi_{\omega}(\gamma\log(K)+t)}{\varphi_{\omega}(t)}$.
\qed\enddemo

\begin{corollary}\label{counterlemma3}
Let $\omega:[0,+\infty)\rightarrow[0,+\infty)$ be non-decreasing and such that $\lim_{t\rightarrow+\infty}\omega(t)=+\infty$. If
\begin{equation}\label{counterlemma3equ1}
\forall\;\gamma>0:\;\;\;\limsup_{t\rightarrow+\infty}\frac{\varphi_{\omega}(\gamma+t)}{\varphi_{\omega}(t)}\le 2
\end{equation}
is valid, then $\gamma(\omega)=+\infty$.
\end{corollary}

\demo{Proof}
In view of Lemma \ref{counterlemma2} property \eqref{counterlemma3equ1} means that $(P_{\omega,\gamma})$ is valid with the (uniform) choice $K:=e$ for any $\gamma>0$.
\qed\enddemo

\subsection{Construction of the (counter-)example}\label{constructioncountersection}
We start with the explicit construction of the weight function and involve the information obtained in Section \ref{counterexprepsection}. The idea is purely geometric and quite technical, therefore we split the proof into several steps.\vspace{6pt}

\emph{Step a - Comments on the geometric definition of $\omega$.}\vspace{6pt}

We introduce $\omega$ by defining explicitly $\varphi_{\omega}: t\mapsto\omega(e^t)$ on $[0,+\infty)$, i.e. $\omega$ on $[1,+\infty)$ via $\omega(s):=\varphi_{\omega}(\log(s))$. The graph of $\varphi_{\omega}$ consists of a set of straight lines and for this we have to introduce \emph{three auxiliary sequences} $(x_j)_{j\in\NN_{>0}}$, $(\overline{x}_j)_{j\in\NN_{>0}}$ and $(y_j)_{j\in\NN_{>0}}$. We assume that
$$\forall\;j\in\NN_{>0}:\;\;\;x_j<\overline{x}_j<y_j<x_{j+1},$$
and several more growth restrictions; see \emph{Step b} below. Indeed, $\overline{x}_j$ is given as the convex combination of the points $x_j$ and $y_j$ subject to a value $0<t_j<1$ which is coming from the \emph{fourth auxiliary sequence} $(t_j)_{j\in\NN_{>0}}$.

The graph of $\varphi_{\omega}$ is defined to be
\begin{itemize}
\item[$(*)$] the straight line connecting the points $(0,0)$ and $(x_1,\varphi_{\omega}(x_1))$ (since we assume $x_1>0$),

\item[$(*)$] the set of straight lines connecting the points $(x_j,\varphi_{\omega}(x_j))$ and $(\overline{x}_j,\varphi_{\omega}(\overline{x}_j))$ resp. $(\overline{x}_j,\varphi_{\omega}(\overline{x}_j))$ and $(y_j,\varphi_{\omega}(y_j))$ resp. $(y_j,\varphi_{\omega}(y_j))$ and $(x_{j+1},\varphi_{\omega}(x_{j+1}))$.

\item[$(*)$] The slope of the straight line connecting $(x_j,\varphi_{\omega}(x_j))$ and $(\overline{x}_j,\varphi_{\omega}(\overline{x}_j))$ is denoted by $k_j$; the slope of the line connecting $(y_j,\varphi_{\omega}(y_j))$ and $(x_{j+1},\varphi_{\omega}(x_{j+1}))$ by $\ell_j$.

\item[$(*)$] Moreover, the slope of the line connecting $(\overline{x}_j,\varphi_{\omega}(\overline{x}_j))$ and $(y_j,\varphi_{\omega}(y_j))$ is identically $0$ for all $j\in\NN_{>0}$, i.e. $\varphi_{\omega}(\overline{x}_j)=\varphi_{\omega}(y_j)$; see \eqref{stepIIIequ4} and \eqref{stepIIIequ6} below.

\item[$(*)$] This describes $\omega$ on $[1,+\infty)$ with $\omega(1)=\varphi_{\omega}(0)=0$ and in order to complete the definition we set $\omega(t):=0$ for $0\le t<1$ (i.e. normalization for $\omega$).
\end{itemize}

Finally, for treating $(V)$ we introduce the following notion: A sequence $\delta=(\delta_j)_{j\in\NN_{>0}}\in\RR_{>0}^{\NN}$ is called an \emph{admissible parameter sequence} if
\begin{equation}\label{step0equ1}
\forall\;j\in\NN_{>0}:\;\;\;\delta_{j+1}\le\delta_j\le(j+2)\delta_{j+1},\hspace{25pt}\lim_{j\rightarrow+\infty}\delta_jj=+\infty.
\end{equation}
An immediate example for such a sequence is $\delta_j:=\frac{1}{\log(e+j)}$.

Let $\delta=(\delta_j)_{j\in\NN_{>0}}$ be an admissible parameter sequence, arbitrary but from now on fixed. The definition of $\omega$ is depending on chosen $\delta$ which suggests the notation $\varphi_{\omega_{\delta}}$. However, since $\delta$ is now fixed, in order to lighten notation we skip $\delta$ and simply write $\varphi_{\omega}$ except in the very last part \emph{Step h:} There we verify $(V)$ and hence we have to emphasize the different parameter sequences.\vspace{6pt}

\emph{Step b - Definition of the auxiliary sequence $(t_j)_{j\in\NN_{>0}}$.}\vspace{6pt}

We choose $0<t_1<1$ arbitrary and then set
\begin{equation}\label{stepIequ1}
\forall\;j\in\NN_{>0}:\;\;\;t_{j+1}:=\frac{1}{2}\left(\frac{1}{j}-t_j\right).
\end{equation}
Let us show by induction that this sequence satisfies
\begin{equation}\label{stepIequ2}
\forall\;j\in\NN_{>0}:\;\;\;0<t_j<\frac{1}{j},
\end{equation}
which should be compared with \eqref{counterlemma1equ2}. The case $j=1$ is obvious and note that $t_{j+1}>0$ if and only if $t_j<\frac{1}{j}$ is verified, see \eqref{stepIequ1}. Moreover,
$$t_{j+1}<\frac{1}{j+1}\Leftrightarrow\frac{1}{j}-t_j<\frac{2}{j+1}\Leftrightarrow\frac{1-j}{j(j+1)}=\frac{1}{j}-\frac{2}{j+1}<t_j,$$
which clearly holds since $t_j>0$ by induction hypothesis. From now on let such a sequence $(t_j)_{j\in\NN_{>0}}$ be given and fixed and note that this sequence is not depending on the chosen admissible parameter sequence $\delta$.\vspace{6pt}

\emph{Step c - Definition of the required auxiliary sequences and basic growth properties.}\vspace{6pt}

We introduce the required auxiliary sequences $(x_j)_{j\in\NN_{>0}}$, $(\overline{x}_j)_{j\in\NN_{>0}}$ and $(y_j)_{j\in\NN_{>0}}$, collect the required growth properties and define $\varphi_{\omega}$ on these values.

Let us set
\begin{equation}\label{stepIIIequ1}
x_1:=\frac{2}{t_1},
\end{equation}
and when $x_j$, $j\in\NN_{>0}$, is given, then define
\begin{equation}\label{stepIIIequ2}
\forall\;j\in\NN_{>0}:\;\;\;\varphi_{\omega}(x_j):=\delta_jjx_j,
\end{equation}
\begin{equation}\label{stepIIIequ3}
\forall\;j\in\NN_{>0}:\;\;\;y_j:=\frac{2\delta_j^{-1}\varphi_{\omega}(x_j)}{1-jt_j}=\frac{2jx_j}{1-jt_j},
\end{equation}
\begin{equation}\label{stepIIIequ4}
\forall\;j\in\NN_{>0}:\;\;\;\varphi_{\omega}(y_j):=\frac{2j\varphi_{\omega}(x_j)}{1-jt_j}=\frac{2\delta_jj^2x_j}{1-jt_j}=\delta_jjy_j,
\end{equation}
\begin{equation}\label{stepIIIequ5}
\forall\;j\in\NN_{>0}:\;\;\;\overline{x}_j:=t_jy_j+(1-t_j)x_j,
\end{equation}
\begin{equation}\label{stepIIIequ6}
	\forall\;j\in\NN_{>0}:\;\;\;\varphi_{\omega}(\overline{x}_j):=\varphi_{\omega}(y_j)\left(=\frac{2\delta_jj^2x_j}{1-jt_j}\right),
\end{equation}
\begin{equation}\label{stepIIIequ7}
\forall\;j\in\NN_{>0}:\;\;\;x_{j+1}:=(j+1)y_j=\frac{2j(j+1)x_j}{1-jt_j}.
\end{equation}
Concerning the denominators in the defining expressions recall that by \eqref{stepIequ2} we get $0<1-jt_j<1$ for any $j\in\NN_{>0}$. Moreover, we immediately have $1<\frac{2}{t_1}=x_1$, that $x_j<\overline{x}_j<y_j<x_{j+1}$ for all $j\in\NN_{>0}$ and that all three sequences tend to infinity as $j\rightarrow+\infty$. By \eqref{stepIIIequ2} the slope of the straight line connecting $(0,0)$ and $(x_1,\varphi_{\omega}(x_1))$ is equal to $\delta_1$. Next we prove
\begin{equation}\label{stepIIIequ8}
\forall\;j\in\NN_{>0}:\;\;\;x_j\ge\frac{2j}{t_j}(>2j^2\ge 2j).
\end{equation}
The case $j=1$ is precisely \eqref{stepIIIequ1}. By \eqref{stepIIIequ7} we get
$$x_{j+1}\ge\frac{2(j+1)}{t_{j+1}}\Leftrightarrow\frac{2j(j+1)x_j}{1-jt_j}\ge\frac{2(j+1)}{t_{j+1}}\Leftrightarrow\frac{jx_j}{1-jt_j}\ge\frac{1}{t_{j+1}},$$
and \eqref{stepIequ1} yields
$$\frac{1}{t_{j+1}}\le\frac{2j^2}{t_j(1-jt_j)}\Leftrightarrow\frac{t_j(1-jt_j)}{2j^2}\le t_{j+1}\Leftrightarrow\frac{t_j}{2j}\left(\frac{1}{j}-t_j\right)\le\frac{1}{2}\left(\frac{1}{j}-t_j\right),$$
which is clearly satisfied (recall that $t_j<1$). Finally, the induction hypothesis implies $\frac{jx_j}{1-jt_j}\ge\frac{2j^2}{t_j(1-jt_j)}$ and by combining everything we have shown the desired estimate for $j+1$, thus \eqref{stepIIIequ8} is verified. Next we study the growth differences between the auxiliary sequences; these identities are then crucial for \emph{Step f} below.\vspace{6pt}

\emph{Claim: One has $\min\{y_j-\overline{x}_j, y_j-x_j, x_{j+1}-y_j, \overline{x}_j-x_j\}\ge j$ for all $j\in\NN_{>0}$.}\vspace{6pt}

First note that, since $x_j<\overline{x}_j$, for having $y_j-x_j\ge j$ it suffices to verify $y_j-\overline{x}_j\ge j$.

Concerning this estimate, by taking into account \eqref{stepIIIequ5} and \eqref{stepIIIequ3} we get
\begin{align*}
&y_j-\overline{x}_j=(1-t_j)y_j-(1-t_j)x_j=(1-t_j)(y_j-x_j)=(1-t_j)\left(\frac{2jx_j}{1-jt_j}-x_j\right)\ge j
\\&
\Leftrightarrow jx_j\ge j\frac{1-jt_j}{2(1-t_j)}+x_j\frac{1-jt_j}{2}.
\end{align*}
This last estimate holds since $\frac{1-jt_j}{2}\le\frac{1}{2}$ and $\frac{1-jt_j}{2(1-t_j)}\le\frac{1}{2}$, so we get $j\frac{1-jt_j}{2(1-t_j)}+x_j\frac{1-jt_j}{2}\le\frac{j+x_j}{2}$ and the estimate $j+x_j\le 2jx_j$ is clear because $x_j\ge j$ for any $j\in\NN_{>0}$ (see \eqref{stepIIIequ8}).

Next, we see that \eqref{stepIIIequ7} yields $x_{j+1}-y_j=jy_j>jx_j\ge j$.

Finally, by \eqref{stepIIIequ5} and \eqref{stepIIIequ3} we get $$\overline{x}_j-x_j=t_j(y_j-x_j)=t_j\left(\frac{2jx_j}{1-jt_j}-x_j\right)\ge j\Leftrightarrow jx_j\ge j\frac{1-jt_j}{2t_j}+x_j\frac{1-jt_j}{2},$$
and let us see that this estimate holds: First, $x_j\frac{1-jt_j}{2}<\frac{x_j}{2}$ is clear and so it suffices to have $\left(j-\frac{1}{2}\right)x_j\ge j\frac{1-jt_j}{2t_j}$. Second, by taking into account \eqref{stepIIIequ8} for this estimate it suffices to have $\left(j-\frac{1}{2}\right)\frac{2j}{t_j}\ge\frac{j}{2t_j}$ and which is clear for all $j\in\NN_{>0}$. Thus the claim is verified.\vspace{6pt}

Now let us compute and compare the slopes $k_j$ and $\ell_j$ introduced in \emph{Step a} explicitly. By definition we have
\begin{align*}
k_j&:=\frac{\varphi_{\omega}(\overline{x}_j)-\varphi_{\omega}(x_j)}{\overline{x}_j-x_j}=\delta_j\frac{2j^2x_j(1-jt_j)^{-1}-jx_j}{t_j(y_j-x_j)}=\delta_j\frac{2j^2x_j-jx_j(1-jt_j)}{t_j(2jx_j(1-jt_j)^{-1}-x_j)(1-jt_j)}
\\&
=\delta_j\frac{jx_j(2j-1+jt_j)}{t_j(2jx_j-x_j+jx_jt_j)}=\delta_j\frac{jx_j(2j-1+jt_j)}{t_jx_j(2j-1+jt_j)},
\end{align*}
and so
\begin{equation}\label{stepIIIequ9}
\forall\;j\in\NN_{>0}:\;\;\;k_j=\frac{\delta_jj}{t_j}.
\end{equation}
Since $k_j=\frac{\delta_jj}{t_j}>\delta_jj^2$ (see \eqref{stepIequ2}), the second part of \eqref{step0equ1} implies $\lim_{j\rightarrow+\infty}k_j=+\infty$.

Similarly, for $\ell_j$ one has by \eqref{stepIIIequ2}, \eqref{stepIIIequ4} and \eqref{stepIIIequ7} that
\begin{align*}
\ell_j&:=\frac{\varphi_{\omega}(x_{j+1})-\varphi_{\omega}(y_j)}{x_{j+1}-y_j}=\frac{\delta_{j+1}(j+1)x_{j+1}-\delta_jjy_j}{x_{j+1}-y_j}=\frac{\delta_{j+1}(j+1)^2-\delta_jj}{j}.
\end{align*}
The first part in \eqref{step0equ1} implies, on the one hand, that $\ell_j\le\frac{\delta_j(j^2+j+1)}{j}$ because $\delta_{j+1}\le\delta_j$, and, on the other hand, also $\ell_j>0\Leftrightarrow\delta_{j+1}\frac{(j+1)^2}{j}>\delta_j$ holds since $\frac{(j+1)^2}{j}=j+2+\frac{1}{j}>j+2$. Consequently, we have verified
\begin{equation}\label{stepIIIequ10}
\forall\;j\in\NN_{>0}:\;\;\;0<\ell_j\le\frac{\delta_j(j^2+j+1)}{j}.
\end{equation}
Moreover, it holds that
\begin{equation}\label{stepIIIequ11}
\forall\;j\in\NN_{>0}:\;\;\;\ell_j<k_{j+1}
\end{equation}
because, first, by \eqref{stepIIIequ9} we get $$k_{j+1}>\ell_j\Leftrightarrow\frac{\delta_{j+1}(j+1)}{t_{j+1}}>\frac{\delta_{j+1}(j+1)^2-\delta_jj}{j}\Leftrightarrow\delta_{j+1}(j+1)>t_{j+1}\left(\delta_{j+1}\left(j+2+\frac{1}{j}\right)-\delta_j\right).$$ This estimate is valid since $$t_{j+1}\left(\delta_{j+1}\left(j+2+\frac{1}{j}\right)-\delta_j\right)<t_{j+1}\delta_{j+1}\left(j+2+\frac{1}{j}\right)<\frac{1}{j+1}\delta_{j+1}\left(j+2+\frac{1}{j}\right),$$ which follows by taking into account \eqref{stepIequ2} and then, in order to conclude, note that $\frac{1}{j+1}\left(j+2+\frac{1}{j}\right)<j+1$ is obvious.\vspace{6pt}

Summarizing, from the given construction it follows that $\varphi_{\omega}:[0,+\infty)\rightarrow[0,+\infty)$ is continuous and non-decreasing with $\lim_{t\rightarrow+\infty}\varphi_{\omega}(t)=+\infty$. The last property follows, in particular, also by the next step.\vspace{6pt}

\emph{Step d - $\omega$ satisfies \hyperlink{om3}{$(\omega_3)$} (and hence \hyperlink{om3w}{$(\omega_{3,w})$}, too).}\vspace{6pt}

In order to conclude we have to show that $\lim_{t\rightarrow+\infty}\frac{\varphi_{\omega}(t)}{t}=+\infty$ and by construction and definition it suffices to investigate the ``corner points'' of the graph; i.e. that
\begin{equation}\label{stepdequ}
\lim_{j\rightarrow+\infty}\frac{\varphi_{\omega}(x_j)}{x_j}=\lim_{j\rightarrow+\infty}\frac{\varphi_{\omega}(y_j)}{y_j}=+\infty.
\end{equation}
For this note that $x\mapsto\frac{kx+d}{x}=k+\frac{d}{x}$ is non-increasing on $(0,+\infty)$ if and only if $d\ge 0$ and in view of the geometric construction recall that the slopes of the lines connecting $(\overline{x}_j,\varphi_{\omega}(\overline{x}_j))$ and $(y_j,\varphi_{\omega}(y_j))$ are identically $0$.

For all $j\in\NN_{>0}$ by \eqref{stepIIIequ2} we get
$$\frac{\varphi_{\omega}(x_j)}{x_j}=\frac{\delta_jjx_j}{x_j}=\delta_jj,$$
and \eqref{stepIIIequ4} implies
$$\frac{\varphi_{\omega}(y_j)}{y_j}=\delta_jj,$$
which proves \eqref{stepdequ} by taking into account the second part in \eqref{step0equ1}.\vspace{6pt}

\emph{Step e - There does not exist $\sigma:[0,+\infty)\rightarrow[0,+\infty)$ such that $\sigma\hyperlink{sim}{\sim}\omega$ and such that $\varphi_{\sigma}$ is convex (i.e. $\sigma$ satisfies in addition \hyperlink{om4}{$(\omega_4)$}).}\vspace{6pt}

We apply Lemma \ref{counterlemma1}. First note that as a consequence of {\itshape Step c} and {\itshape Step d} we have that $\lim_{t\rightarrow+\infty}\varphi_{\omega}(t)=+\infty$ and hence, in particular, \eqref{counterlemma1equ1} is valid. In order to conclude we prove
\eqref{counterlemma1equ2} and verify this condition for the choices $t_j$ from {\itshape Step b} (recall \eqref{stepIequ2}) and when taking $a_j:=x_j$, $b_j:=y_j$. By combining \eqref{stepIIIequ5}, \eqref{stepIIIequ6} and \eqref{stepIIIequ2} we get for all $j\in\NN_{>0}$
$$\varphi_{\omega}(t_jy_j+(1-t_j)x_j)=\varphi_{\omega}(\overline{x}_j)=\varphi_{\omega}(y_j)=\frac{2\delta_jj^2x_j}{1-jt_j}=\frac{2j\varphi_{\omega}(x_j)}{1-jt_j},$$
which is precisely \eqref{counterlemma1equ2} for the aforementioned choices.\vspace{6pt}

\emph{Step f - $\omega$ satisfies $\gamma(\omega)=+\infty$.}\vspace{6pt}

In order to prove $\gamma(\omega)=+\infty$ we want to apply Corollary \ref{counterlemma3} and so the aim is to verify \eqref{counterlemma3equ1}. In fact we even show
\begin{equation}\label{stepVIequ1}
\forall\;\gamma>0:\;\;\;\lim_{s\rightarrow+\infty}\frac{\varphi_{\omega}(\gamma+s)}{\varphi_{\omega}(s)}=1\Leftrightarrow\lim_{t\rightarrow+\infty}\frac{\omega(e^{\gamma}t)}{\omega(t)}=1,
\end{equation}
i.e. $\omega$ is \emph{slowly varying,} see \cite[$(1.2.1)$]{regularvariation}, \cite[Sect. 2.2]{index} and the comments at the end of Section \ref{growthrelsection}.\vspace{6pt}

Let $\gamma>0$ be given, arbitrary large but from now on fixed. By the previous \emph{Claim in Step c} it follows that for all sufficiently large $s\in[0,+\infty)$ (depending on given $\gamma$) we eventually have to distinguish between the following three different cases.\vspace{6pt}

{\itshape Case i} - Let $x_j\le s\le\overline{x}_j$, for $j\in\NN_{>0}$ sufficiently large depending on $\gamma$, then by the claim $\gamma+s\le y_j$.

In this situation, by definition of the graph and the slope $k_j$ and since $\varphi_{\omega}$ is non-decreasing, we have the following estimate:
$$\frac{\varphi_{\omega}(\gamma+s)}{\varphi_{\omega}(s)}\le\frac{\varphi_{\omega}(s)}{\varphi_{\omega}(s)}+\frac{\gamma k_j}{\varphi_{\omega}(s)}\le 1+\frac{\gamma k_j}{\varphi_{\omega}(x_j)}.$$
\eqref{stepIIIequ2}, \eqref{stepIIIequ9} and \eqref{stepIIIequ8} together imply
$$\frac{k_j}{\varphi_{\omega}(x_j)}=\frac{\delta_jj}{t_j\delta_jjx_j}=\frac{1}{t_jx_j}\le\frac{1}{2j},$$
which tends to $0$ as $j\rightarrow+\infty$.\vspace{6pt}

{\itshape Case ii} - Let $\overline{x}_j\le s\le y_j$, then the claim gives that $\gamma+s\le x_{j+1}$ for all $j\in\NN_{>0}$ sufficiently large.

In this situation, by definition of the graph and the slope $\ell_j$ and since $\varphi_{\omega}$ is non-decreasing, we have the following estimate:
$$\frac{\varphi_{\omega}(\gamma+s)}{\varphi_{\omega}(s)}\le\frac{\varphi_{\omega}(s)}{\varphi_{\omega}(s)}+\frac{\gamma\ell_j}{\varphi_{\omega}(s)}\le 1+\frac{\gamma\ell_j}{\varphi_{\omega}(\overline{x}_j)}.$$
\eqref{stepIIIequ6}, \eqref{stepIIIequ10} and \eqref{stepIIIequ8} together with \eqref{stepIequ2} imply
$$\frac{\ell_j}{\varphi_{\omega}(\overline{x}_j)}\le\frac{\delta_j(j^2+j+1)}{j2j^2\delta_jx_j(1-jt_j)^{-1}}\le\frac{(j^2+j+1)(1-jt_j)t_j}{4j^4}\le\frac{j^2+j+1}{4j^4},$$
which tends to $0$ as $j\rightarrow+\infty$.\vspace{6pt}

{\itshape Case iii} - Let $y_j\le s\le x_{j+1}$, then $\gamma+s\le\overline{x}_{j+1}$ for all $j\in\NN_{>0}$ sufficiently large.

In this situation we have, by taking into account \eqref{stepIIIequ11}, that
$$\frac{\varphi_{\omega}(\gamma+s)}{\varphi_{\omega}(s)}\le\frac{\varphi_{\omega}(s)}{\varphi_{\omega}(s)}+\frac{\gamma\max\{\ell_j,k_{j+1}\}}{\varphi_{\omega}(s)}\le 1+\frac{\gamma k_{j+1}}{\varphi_{\omega}(y_j)}.$$
Then note that by \eqref{stepIequ1} and the fact that $t_j<1$ for all $j\in\NN_{>0}$ we have $\frac{(1-jt_j)t_j}{t_{j+1}}\le 2j\Leftrightarrow(1-jt_j)t_j\le(1-jt_j)$. By combining this with \eqref{stepIIIequ4}, \eqref{stepIIIequ9}, \eqref{stepIIIequ8} and property $\delta_{j+1}\le\delta_j$ (see the first part in \eqref{step0equ1}) we obtain
$$\frac{k_{j+1}}{\varphi_{\omega}(y_j)}=\frac{\delta_{j+1}(j+1)}{t_{j+1}2\delta_jj^2x_j(1-jt_j)^{-1}}\le\frac{(j+1)(1-jt_j)t_j}{t_{j+1}4j^3}\le\frac{j+1}{2j^2},$$
which tends to $0$ as $j\rightarrow+\infty$.

Summarizing, {\itshape Cases i-iii} verify \eqref{stepVIequ1} and hence this part is finished.\vspace{6pt}

\emph{Step g - $\omega$ satisfies \hyperlink{omnq}{$(\omega_{\on{nq}})$} and \hyperlink{alpha0}{$(\alpha_0)$}; more precisely, $\omega$ is equivalent to $\kappa_{\omega}$ (see \eqref{fctkappa}).}\vspace{6pt}

The fact that $\omega$ has \hyperlink{omnq}{$(\omega_{\on{nq}})$} follows from \emph{Step f} and \eqref{chaimofimpl} whereas \hyperlink{alpha0}{$(\alpha_0)$} holds by \emph{Step f} and \eqref{chaimofimpl1}. By Remark \ref{secondcomprem} the function $\omega$ is equivalent to $\kappa_{\omega}$ which is continuous, sub-additive and $\kappa_{\omega}(0)=\omega(0)=0$.\vspace{6pt}

\emph{Step h - There exist (uncountable) infinitely many pairwise non-equivalent weights satisfying all the requirements shown above.}\vspace{6pt}

Indeed, we are going to verify that if $\delta=(\delta_j)_{j\in\NN_{>0}}$ and $\delta'=(\delta'_j)_{j\in\NN_{>0}}$ are two admissible parameter sequences satisfying
\begin{equation}\label{admissiblesequencequ}
\sup_{j\in\NN_{>0}}\frac{\delta_j}{\delta'_j}=+\infty\;\;\;\text{or}\;\;\;\sup_{j\in\NN_{>0}}\frac{\delta'_j}{\delta_j}=+\infty,
\end{equation}
then the corresponding weights $\omega_{\delta}$ and $\omega_{\delta'}$ are not equivalent. From this the assertion follows immediately; see also Remark \ref{admissiblerem} for an explicit one-parameter family.\vspace{6pt}

For this apply \eqref{stepIIIequ2} and when assuming the first condition in \eqref{admissiblesequencequ} we get
$$\sup_{x\ge 0}\frac{\varphi_{\omega_{\delta}}(x)}{\varphi_{\omega_{\delta'}}(x)}\ge\sup_{j\in\NN_{>0}}\frac{\varphi_{\omega_{\delta}}(x_j)}{\varphi_{\omega_{\delta'}}(x_j)}=\sup_{j\in\NN_{>0}}\frac{\delta_j jx_j}{\delta'_jjx_j}=\sup_{j\in\NN_{>0}}\frac{\delta_j}{\delta'_j}=+\infty.$$
This shows that $\varphi_{\omega_{\delta}}\hyperlink{sim}{\sim}\varphi_{\omega_{\delta'}}$, equivalently $\omega_{\delta}\hyperlink{sim}{\sim}\omega_{\delta'}$, cannot be satisfied. The argument for the second condition in \eqref{admissiblesequencequ} is analogous. Alternatively, instead of \eqref{stepIIIequ2} one can also use \eqref{stepIIIequ4} or \eqref{stepIIIequ6} and for this recall the following facts: Of course, the values of $\varphi_{\omega}$ are depending, by definition, on chosen $\delta$. On the other hand, note that the auxiliary sequence $(t_j)_{j\in\NN_{>0}}$ studied in \emph{Step b,} the ``corner point'' sequences $(x_j)_{j\in\NN_{>0}}$, $(y_j)_{j\in\NN_{>0}}$ and thus also the intermediate sequence $(\overline{x}_j)_{j\in\NN_{>0}}$ are \emph{not depending on the admissible parameter sequence.} In this sense the geometric construction is made uniformly w.r.t. the parameter sequences $\delta$.\vspace{6pt}

This finishes the construction of the desired (counter-)example resp. of the uncountable infinite family of (counter-)examples.

\begin{remark}\label{admissiblerem}
Let an admissible parameter sequence $\delta$ satisfying $\lim_{j\rightarrow+\infty}\delta_j=0$ be given; consider e.g. $\delta_j:=\frac{1}{\log(e+j)}$. However, note that this additional condition is not required necessarily in the above arguments and $\delta$ can be bounded by below by some $\epsilon>0$. In any case, under this assumption each power sequence $\delta^{\alpha}=(\delta_j^{\alpha})_{j\in\NN_{>0}}$, $0<\alpha<1$, is again an admissible parameter sequence since $\delta_{j+1}^{\alpha}\le\delta_j^{\alpha}\le(j+2)^{\alpha}\delta_{j+1}^{\alpha}\le(j+2)\delta_{j+1}^{\alpha}$ for all $j\in\NN_{>0}$ and $\delta_j^{\alpha}\ge\delta_j$ for all $j$ such that $\delta_j\le 1$; i.e. for all $j$ sufficiently large. Thus all requirements from \eqref{step0equ1} hold for $\delta^{\alpha}:=(\delta^{\alpha}_j)_{j\in\NN_{>0}}$ and (each) $\delta':=\delta^{\alpha}$ serves as an example for an admissible parameter sequence satisfying the second condition in \eqref{admissiblesequencequ}; more precisely $\sup_{j\in\NN_{>0}}\frac{\delta^{\alpha}_j}{\delta^{\beta}_j}=+\infty$ for all $0<\alpha<\beta\le 1$.

Note that in general it is not clear that $\delta^{\alpha}$ is admissible if $\alpha>1$.
\end{remark}

\begin{remark}\label{nonassorem}
Let $\mathbf{M}=(M_p)_{p\in\NN}\in\RR_{>0}^{\NN}$ be given, then the \emph{associated function} $\omega_{\mathbf{M}}: [0,+\infty)\rightarrow[0,+\infty)\cup\{+\infty\}$ is defined as follows:
\begin{equation}\label{assofunc}
\omega_{\mathbf{M}}(t):=\sup_{p\in\NN}\log\left(\frac{M_0t^p}{M_p}\right),\qquad t\ge 0,
\end{equation}
with the conventions $0^0:=1$ and $\log(0)=-\infty$. When treating this function, naturally one should assume that $\lim_{p\rightarrow+\infty}(M_p)^{1/p}=+\infty$ holds in order to ensure $\omega_{\mathbf{M}}: [0,+\infty)\rightarrow[0,+\infty)$ (well-definedness). For more details we refer to \cite[Chapitre I]{mandelbrojtbook}, the more recent work \cite{regularnew} and see also \cite[Sect. 3]{Komatsu73}. It is known and immediate by definition that \hyperlink{om4}{$(\omega_4)$} holds for $\omega_{\mathbf{M}}$ and consequently in the whole equivalence class (w.r.t. relation \hyperlink{sim}{$\sim$}) of the previously constructed function there does not exist any associated weight function. In particular, the obtained (counter-)example cannot be described within the weight sequence setting which means that $\omega\equiv\omega_{\mathbf{M}}$ fails for any weight sequence $\mathbf{M}$.

On the other hand the above constructed $\omega$ is matrix admissible since crucially \hyperlink{om3}{$(\omega_3)$} holds.
\end{remark}

We finish by mentioning the following \emph{question/goal:} Transfer the construction of $\varphi_{\omega}$ (resp. of $\omega$) in Section \ref{failuresection} to the higher dimensional case (anisotropic setting). We expect that the analogous statement holds as well but such an explicit construction is more technical and involved and also the arguments using the growth index $\gamma(\cdot)$ are becoming unclear in this general situation since so far this quantity is only defined and studied for $\omega: [0,+\infty)\rightarrow[0,+\infty)$.

\bibliographystyle{plain}
\bibliography{Bibliography}
\end{document}